\documentclass{amsart}
\usepackage{amssymb,epsf,amscd}
\allowdisplaybreaks
\newtheorem{theorem}[subsubsection]{Theorem}               
\newtheorem{acknowledgements}[subsubsection]{Acknowledgements}

\newtheorem{diagram}[subsubsection]{Diagram}        
\newtheorem{lemma}[subsubsection]{Lemma}                   
\newtheorem{notation}[subsubsection]{Notation}                   
\newtheorem{corollary}[subsubsection]{Corollary}           
                   
\newtheorem{conjecture}[subsubsection]{Conjecture}        
\newtheorem{exercize}[subsubsection]{Exercise}        
         
\newtheorem{note}[subsubsection]{Note}         
\theoremstyle{definition}                      
\newtheorem{definition}[subsubsection]{Definition}         
\newtheorem{altdefinition}[subsubsection]{Alternative Definition}         
\newtheorem{observation}[subsubsection]{Observation}         
\newtheorem{example}[subsubsection]{Example}                
               
\theoremstyle{remark}
\def\rtb{
\begin{flushright}
\boxed{}
\end{flushright}}
\newtheorem{remark}[subsubsection]{Remark}                 
               
\numberwithin{equation}{subsection}
\def\Arc[#1]{
\ifcase#1
\qbezier[25](0.966,-0.259)(1.04,0)(0.966,0.259)
\or
\qbezier[25](0.966,0.259)(0.897,0.518)(0.707,0.707)
\or
\qbezier[25](0.707,0.707)(0.518,0.897)(0.259,0.966)
\or
\qbezier[25](0.259,0.966)(0,1.04)(-0.259,0.966)
\or
\qbezier[25](-0.259,0.966)(-0.518,0.897)(-0.707,0.707)
\or
\qbezier[25](-0.707,0.707)(-0.897,0.518)(-0.966,0.259)
\or
\qbezier[25](-0.966,0.259)(-1.04,0)(-0.966,-0.259)
\or
\qbezier[25](-0.966,-0.259)(-0.897,-0.518)(-0.707,-0.707)
\or
\qbezier[25](-0.707,-0.707)(-0.518,-0.897)(-0.259,-0.966)
\or
\qbezier[25](-0.259,-0.966)(0,-1.04)(0.259,-0.966)
\or
\qbezier[25](0.259,-0.966)(0.518,-0.897)(0.707,-0.707)
\or
\qbezier[25](0.707,-0.707)(0.897,-0.518)(0.966,-0.259)
\fi}

\def\DottedCircle{
\qbezier[4](0.966,-0.259)(1.04,0)(0.966,0.259)
\qbezier[4](0.966,0.259)(0.897,0.518)(0.707,0.707)
\qbezier[4](0.707,0.707)(0.518,0.897)(0.259,0.966)
\qbezier[4](0.259,0.966)(0,1.04)(-0.259,0.966)
\qbezier[4](-0.259,0.966)(-0.518,0.897)(-0.707,0.707)
\qbezier[4](-0.707,0.707)(-0.897,0.518)(-0.966,0.259)
\qbezier[4](-0.966,0.259)(-1.04,0)(-0.966,-0.259)
\qbezier[4](-0.966,-0.259)(-0.897,-0.518)(-0.707,-0.707)
\qbezier[4](-0.707,-0.707)(-0.518,-0.897)(-0.259,-0.966)
\qbezier[4](-0.259,-0.966)(0,-1.04)(0.259,-0.966)
\qbezier[4](0.259,-0.966)(0.518,-0.897)(0.707,-0.707)
\qbezier[4](0.707,-0.707)(0.897,-0.518)(0.966,-0.259)
}

\newcommand\lbb[1]{\label{#1} 
                   }                                    
\def\Zset{{\mathbb Z}}       
\def\Rset{{\mathbb R}}       
\def\Qset{{\mathbb Q}}       

\def\ZHS{\mathbb{Z}HS^3}

\def\be{\beta}

\def\be{\begin{equation}}
\def\fe{\end{equation}}
\def\KI{Z}
\def\KIC{\check{Z}}
\def\KIH{\hat{Z}}
\def\wh{\omega}
\def\CA{A}
\def\CB{B}
\def\Picture#1{
\begin{picture}(2,1)(0,0)
#1
\end{picture}
}

\begin{document}
\title{The lines of the Kontsevich integral and Rozansky's rationality conjecture}
\author{Andrew Kricker}
\address{Department of Mathematical and Computing Sciences, 
Tokyo Institute of Technology}
\email{kricker@is.titech.ac.jp}
\date{First version: February , 2000. Current version: May , 2000.
Your comments are welcome.}

\begin{abstract}
This work develops some technology for accessing the loop expansion of
the Kontsevich integral of a knot. The setting is an application of the LMO invariant
to certain surgery presentations of knots by framed links in the solid torus. 
A consequence of this technology is a certain recent conjecture of
Rozansky's. Rozansky conjectured that the Kontsevich integral could be
organised into a series of ``lines'' which could be represented by
finite $\Qset$-linear combinations of diagrams whose edges were labelled,
in an appropriate sense, with rational functions. Furthermore, the 
conjecture requires that the denominator of the rational functions
be at most the Alexander polynomial of the knot. This conjecture is
obtained from an Aarhus-style surgery formula for this setting which we
expect will have other applications.
\end{abstract}
\maketitle

\pagestyle{myheadings}
\markboth{ANDREW KRICKER}{ON ROZANSKY'S RATIONALITY CONJECTURE}

\tableofcontents
\def\typeexamp{
\Picture{
\qbezier(0,4)(2,3)(4,2)
\qbezier(0,4)(-2,3)(-4,2)
\qbezier(0,-4)(2,-3)(4,-2)
\qbezier(0,-4)(-2,-3)(-4,-2)
\qbezier(4,2)(4,0)(4,-2)
\qbezier(-4,2)(-4,0)(-4,-2)
\put(-0.25,4.5){$x_i$}
\put(-0.25,-5.25){$x_l$}
\put(4.25,2){$x_n$}
\put(4.25,-2){$x_m$}
\put(-5,2){$x_j$}
\put(-5,-2){$x_k$}
\put(-2,2){$r$}
\put(-2,-2.6){$r$}
\put(3,0){$r$}
\put(-3.5,0){$g$}
\put(1,2){$g$}
\put(1,-2.75){$g$}
}}

\def\seriesexamp{
\Picture{
\qbezier[18](-2,0)(-2,2)(0,2)
\qbezier[18](2,0)(2,2)(0,2)
\qbezier[18](-2,0)(-2,-2)(0,-2)
\qbezier[18](2,0)(2,-2)(0,-2)
\put(0,-2){\vector(1,0){0.01}}
\put(0,-3){$F$}
}}

\def\redleft{
\Picture{
\qbezier(-3,0)(0,0)(3,0)
\put(-3.6,-1){$x_i$}
\put(3,-1){$x_j$}
\put(-0.25,0.5){$r$}
}}

\def\redright{
\Picture{
\qbezier[24](-3,0)(0,0)(3,0)
\put(-2,0.5){$W_{ij}(T,e^{u})$}
\put(2,0){\vector(1,0){0.01}}
}}

\def\greenright{
\Picture{
\qbezier[24](-3,0)(0,0)(3,0)
\put(-0.8,0.5){$L^{-1}_{ij}$}
\put(2,0){\vector(1,0){0.01}}
}}

\def\greenleft{
\Picture{
\qbezier(-3,0)(0,0)(3,0)
\put(-3.6,-1){$x_i$}
\put(3,-1){$x_j$}
\put(-0.25,0.5){$g$}
}}

\def\crossingchordself{
\Picture{
\put(-4.75,-1){$1$}
\put(4.25,-1){$1$}
\put(-4,-3){\vector(0,1){6}}
\put(4,-3){\vector(0,1){6}}
\put(-0.5,-3){\vector(0,1){6}}
\put(-4.5,-4){$x_1$}
\put(-1,-4){$x_i$}
\qbezier[7](-3.5,-2)(-2.25,-2)(-1,-2)
\qbezier[12](0,-2)(1.75,-2)(3.5,-2)
\qbezier[16](-0.5,1.5)(1.5,1.5)(1.5,0.5)
\qbezier[16](-0.5,-0.5)(1.5,-0.5)(1.5,0.5)
\put(-1.25,2){1}
\put(-1.25,0.25){1}
\put(-1.25,-1.75){1}
\put(-0.5,1.5){\vector(-1,0){0.01}}
\put(0.75,1.5){$t^{\epsilon(i,i,c)}$}
\put(3.5,-4){$x_\mu$}
}}

\def\crossingchordselfagain{
\Picture{
\put(-4.75,-1){$1$}
\put(4.25,-1){$1$}
\put(-4,-3){\vector(0,1){6}}
\put(4,-3){\vector(0,1){6}}
\put(-0.5,-3){\vector(0,1){6}}
\put(-4.5,-4){$x_1$}
\put(-1,-4){$x_i$}
\qbezier[7](-3.5,-2)(-2.25,-2)(-1,-2)
\qbezier[12](0,-2)(1.75,-2)(3.5,-2)
\qbezier[16](-0.5,1.5)(1.5,1.5)(1.5,0.5)
\qbezier[16](-0.5,-0.5)(1.5,-0.5)(1.5,0.5)
\put(-1.25,2){1}
\put(-1.25,0.25){1}
\put(-1.25,-1.75){1}
\put(-0.5,1.5){\vector(-1,0){0.01}}
\put(0.75,1.5){$t^{-\epsilon(i,i,c)}$}
\put(3.5,-4){$x_\mu$}
}}

\def\crossingchordselfloop{
\Picture{
\put(-4.75,-1){$1$}
\put(4.25,-1){$1$}
\put(-4,-3){\vector(0,1){6}}
\put(4,-3){\vector(0,1){6}}
\put(-0.5,-3){\vector(0,1){6}}
\put(-4.5,-4){$x_1$}
\put(-1,-4){$x_i$}
\qbezier[7](-3.5,-2)(-2.25,-2)(-1,-2)
\qbezier[12](0,-2)(1.75,-2)(3.5,-2)
\qbezier[8](-0.5,0.5)(0,0.5)(0.5,0.5)
\qbezier[12](0.5,0.5)(0.5,-0.5)(1.5,-0.5)
\qbezier[12](0.5,0.5)(0.5,1.5)(1.5,1.5)
\qbezier[12](1.5,1.5)(2.5,1.5)(2.5,0.5)
\qbezier[12](1.5,-0.5)(2.5,-0.5)(2.5,0.5)
\put(1.5,-0.5){\vector(1,0){0.01}}
\put(-1.25,2){1}
\put(-0.2,-0.4){1}
\put(-1.25,-1.75){1}
\put(0.75,1.65){$t^{\epsilon(i,i,c)}$}
\put(3.5,-4){$x_\mu$}
}}

\def\crossingchordselfb{
\Picture{
\put(-4.75,-1){$1$}
\put(4.25,-1){$1$}
\put(-4,-3){\vector(0,1){6}}
\put(4,-3){\vector(0,1){6}}
\put(-0.5,-3){\vector(0,1){6}}
\put(-4.5,-4){$x_1$}
\put(-1,-4){$x_i$}
\qbezier[7](-3.5,-2)(-2.25,-2)(-1,-2)
\qbezier[12](0,-2)(1.75,-2)(3.5,-2)
\qbezier[16](-0.5,1.5)(1.5,1.5)(1.5,0.5)
\qbezier[16](-0.5,-0.5)(1.5,-0.5)(1.5,0.5)
\put(-1.25,2){1}
\put(-1.25,0.25){1}
\put(-1.25,-1.75){1}
\put(-0.5,-0.5){\vector(-1,0){0.01}}
\put(0.75,1.5){$t^{\epsilon(i,i,c)}$}
\put(3.5,-4){$x_\mu$}
}}

\def\crossingchord{
\Picture{
\put(-4.75,-1){$1$}
\put(-2.5,-1){$1$}
\put(-2.5,2){$1$}
\put(2,2){$1$}
\put(3.25,-1){$1$}
\put(1,-1){$1$}
\put(-4,-3){\vector(0,1){6}}
\put(4,-3){\vector(0,1){6}}
\put(1.75,-3){\vector(0,1){6}}
\put(-1.75,-3){\vector(0,1){6}}
\put(-4.5,-4){$x_1$}
\put(-2.25,-4){$x_i$}
\put(1.25,-4){$x_j$}
\put(3.5,-4){$x_\mu$}
\qbezier[3](-3.5,-2)(-2.91,-2)(-2.25,-2)
\qbezier[3](3.5,-2)(2.91,-2)(2.25,-2)
\qbezier[7](-1.25,-2)(0,-2)(1.25,-2)
\qbezier[24](-1.75,0.5)(0,0.5)(1.75,0.5)
\put(-1.35,0.85){$t^{\epsilon(i,j,c)}$}
\put(1,0.5){\vector(1,0){0.01}}
}}

\def\crossingchordb{
\Picture{
\put(-4.75,-1){$1$}
\put(-2.5,-1){$1$}
\put(-2.5,2){$1$}
\put(2,2){$1$}
\put(3.25,-1){$1$}
\put(1,-1){$1$}
\put(-4,-3){\vector(0,1){6}}
\put(4,-3){\vector(0,1){6}}
\put(1.75,-3){\vector(0,1){6}}
\put(-1.75,-3){\vector(0,1){6}}
\put(-4.5,-4){$x_1$}
\put(-2.25,-4){$x_i$}
\put(1.25,-4){$x_j$}
\put(3.5,-4){$x_\mu$}
\qbezier[3](-3.5,-2)(-2.91,-2)(-2.25,-2)
\qbezier[3](3.5,-2)(2.91,-2)(2.25,-2)
\qbezier[7](-1.25,-2)(0,-2)(1.25,-2)
\qbezier[24](-1.75,0.5)(0,0.5)(1.75,0.5)
\put(-1.35,0.85){$t^{\epsilon(i,j,c)}$}
\put(1,0.5){\vector(1,0){0.01}}
}}

\def\crossingchordc{
\Picture{
\put(-4.75,-1){$1$}
\put(-2.5,-1){$1$}
\put(-2.5,2){$1$}
\put(2,2){$1$}
\put(3.25,-1){$1$}
\put(1,-1){$1$}
\put(-4,-3){\vector(0,1){6}}
\put(4,-3){\vector(0,1){6}}
\put(1.75,-3){\vector(0,1){6}}
\put(-1.75,-3){\vector(0,1){6}}
\put(-4.5,-4){$x_1$}
\put(-2.25,-4){$x_i$}
\put(1.25,-4){$x_j$}
\put(3.5,-4){$x_\mu$}
\qbezier[3](-3.5,-2)(-2.91,-2)(-2.25,-2)
\qbezier[3](3.5,-2)(2.91,-2)(2.25,-2)
\qbezier[7](-1.25,-2)(0,-2)(1.25,-2)
\qbezier[24](-1.75,0.5)(0,0.5)(1.75,0.5)
\put(-1.35,0.85){$t^{\epsilon(j,i,c)}$}
\put(1,0.5){\vector(-1,0){0.01}}
}}

\def\selfa{
\Picture{
\qbezier[24](-1,-1)(-1,1.5)(0,1.5)
\qbezier[24](1,-1)(1,1.5)(0,1.5)
\put(1,-1){\vector(0,-1){0.01}}
\put(-1.4,-2){$x_i$}
\put(0.6,-2){$x_i$}
\put(-0.75,2){$t^{\epsilon(i,i,c)}$}
}}

\def\selfap{
\Picture{
\qbezier[24](-1,-1)(-1,1.5)(0,1.5)
\qbezier[24](1,-1)(1,1.5)(0,1.5)
\put(1,-1){\vector(0,-1){0.01}}
\put(-1.4,-2){$x_i$}
\put(0.6,-2){$x_j$}
\put(-0.75,2){$t^{\epsilon(i,j,c)}$}
}}

\def\selfbp{
\Picture{
\qbezier[24](-1,-1)(-1,1.5)(0,1.5)
\qbezier[24](1,-1)(1,1.5)(0,1.5)
\put(1,-1){\vector(0,-1){0.01}}
\put(-1.4,-2){$x_j$}
\put(0.6,-2){$x_i$}
\put(-0.75,2){$t^{\epsilon(j,i,c)}$}
}}

\def\selfc{
\Picture{
\qbezier[24](-1,-1)(-1,1.5)(0,1.5)
\qbezier[24](1,-1)(1,1.5)(0,1.5)
\put(1,-1){\vector(0,-1){0.01}}
\put(-1.4,-2){$x_i$}
\put(0.6,-2){$x_i$}
\put(-0.75,2){$t^{-\epsilon(i,i,c)}$}
}}

\def\selfb{
\Picture{
\qbezier[12](0,-1)(0,0)(0,1)
\qbezier[12](0,1)(1,1)(1,2)
\qbezier[12](0,1)(-1,1)(-1,2)
\qbezier[12](-1,2)(-1,3)(0,3)
\qbezier[12](1,2)(1,3)(0,3)
\put(-0.4,-2){$x_i$}
\put(1.5,2){$t^{\epsilon(i,i,c)}$}
\put(-1,2){\vector(0,1){0.01}}
}}

\def\gluechorda{
\Picture{
\qbezier[24](-2,2)(-2,0)(-2,-2)
\put(-2.25,-3){$x_a$}
\qbezier[24](2,2)(2,0)(2,-2)
\put(1.75,-3){$x_b$}
}}

\def\gluechordb{
\Picture{
\qbezier[24](-2,2)(-2,0)(-2,-2)
\put(-2,-2){\circle*{0.25}}
\qbezier[24](2,2)(2,0)(2,-2)
\put(2,-2){\circle*{0.25}}
\qbezier[16](-2,-2)(-2,-4)(0,-4)
\qbezier[16](2,-2)(2,-4)(0,-4)
\put(-4.5,-5.25){$-\frac{1}{2} \sum_{i,j} \delta_{a,i} L^{-1}_{ij} \delta_{b,j}$}
}}

\def\gluechordc{
\Picture{
\qbezier[24](-2,2)(-2,0)(-2,-2)
\put(-2,-2){\circle*{0.25}}
\qbezier[24](2,2)(2,0)(2,-2)
\put(2,-2){\circle*{0.25}}
\qbezier[16](-2,-2)(-2,-4)(0,-4)
\qbezier[16](2,-2)(2,-4)(0,-4)
\put(-1.2,-5.25){$-L^{-1}_{ab}$}
}}

\def\Ddiag{
\Picture{
\put(0,-0.5){$
\Picture{
\qbezier[16](2,2)(2,1)(2,0)
\qbezier[16](6,2)(6,1)(6,0)
\qbezier[5](2.5,0.5)(4,0.5)(5.5,0.5)
\qbezier(1,2)(4,2)(7,2)
\qbezier(1,2)(0,2)(0,3)
\qbezier(1,4)(0,4)(0,3)
\qbezier(1,4)(4,4)(7,4)
\qbezier(7,2)(8,2)(8,3)
\qbezier(7,4)(8,4)(8,3)
\put(3.5,2.75){$D$}
\put(1.75,-1.1){${y_1}$}
\put(5.75,-1.1){${y_L}$}
}$}}
}

\def\subdoublea{
\Picture{
\qbezier[16](-2,-3)(-2,1)(0,1)
\qbezier[16](2,-3)(2,1)(0,1)
\put(2,-3){\vector(0,-1){0.01}}
\put(-2.5,-4.5){$x_i'$}
\put(1.5,-4.5){$x_l'$}
\put(-2.5,2){$\frac{1}{2}M_{ij}L_{jk}M_{kl}$}
}
}

\def\subdoubleb{
\Picture{
\qbezier[16](-2,-3)(-2,-7)(0,-7)
\qbezier[16](2,-3)(2,-7)(0,-7)
\qbezier[16](-6,-3)(-6,1)(-4,1)
\qbezier[16](-2,-3)(-2,1)(-4,1)
\qbezier[16](6,-3)(6,1)(4,1)
\qbezier[16](2,-3)(2,1)(4,1)
\put(-2,-8.5){$-\frac{1}{2} L^{-1}_{kl}$}
\put(-2,-3){\vector(0,-1){0.01}}
\put(-6.5,-4.5){$x_i'$}
\put(-5.5,2){$M^{-1}_{ij}L_{jk}
$}
}
}
\def\subadiag{
\Picture{
\qbezier[24](0,3)(0,0)(0,-3)
\put(-0.4,-4.5){$x_i$}
}}

\def\subadiagstrange{
\Picture{
\qbezier[24](0,3)(0,0)(0,-3)
\put(0,-3){\circle*{0.4}}
\put(4,-3){\circle*{0.4}}
\put(5,2){$L_{jk} M_{kl}$}
\put(0.2,-8.5){$-L^{-1}_{ij} $}
\qbezier[16](0,-3)(0,-7)(2,-7)
\qbezier[16](4,-3)(4,-7)(2,-7)
\qbezier[16](4,-3)(4,1)(6,1)
\qbezier[16](8,-3)(8,1)(6,1)
\put(8,-3){\vector(0,-1){0.01}}
\put(7.5,-4.5){$x_l'$}
}}

\def\subbdiag{
\Picture{
\put(0,-3){\vector(0,-1){0.01}}
\qbezier[24](0,3)(0,0)(0,-3)
\put(0.5,-0.5){$M_{ij}$}
\put(-0.5,-4.5){$x'_j$}
}}

\def\Mijkchord{
\ \ \Picture{
\qbezier[24](0.5,-0.25)(0.5,1.5)(1,1.5)
\qbezier[24](1.5,-0.25)(1.5,1.5)(1,1.5)
\put(-0.1,-1){$x_i$}
\put(1.27,-1){$x_k'$}
\put(1.5,-0.25){\vector(0,-1){0.01}}
\put(-0.5,2){$L_{ij}M_{jk}$}
}\ \ }

\def\Mchord{
\ \ \Picture{
\qbezier[24](0.5,-0.25)(0.5,1.5)(1,1.5)
\qbezier[24](1.5,-0.25)(1.5,1.5)(1,1.5)
\put(-0.1,-1){$x_i$}
\put(1.27,-1){$x_j'$}
\put(1.5,-0.25){\vector(0,-1){0.01}}
\put(0.5,2){$M_{ij}$}
}\ \ }

\def\Mijkchordy{
\ \ \Picture{
\qbezier[24](0.5,-0.25)(0.5,1.5)(1,1.5)
\qbezier[24](1.5,-0.25)(1.5,1.5)(1,1.5)
\put(-0.1,-1){$y_i$}
\put(1.27,-1){$x_k'$}
\put(1.5,-0.25){\vector(0,-1){0.01}}
\put(-0.5,2){$L_{ij}M_{jk}$}
}\ \ }

\def\Mijklchord{
\ \ \ \ \Picture{
\qbezier[24](0.5,-0.25)(0.5,1.5)(1,1.5)
\qbezier[24](1.5,-0.25)(1.5,1.5)(1,1.5)
\put(-0.1,-1){$x_i'$}
\put(1.27,-1){$x_l'$}
\put(1.5,-0.25){\vector(0,-1){0.01}}
\put(-1.45,2.25){$M_{ij}L_{jk}M_{kl}$}
}\ \ \ \ }

\def\legf{
\Picture{
\qbezier[24](0,-3)(0,0)(0,3)
\put(-0.25,3.5){$x_1$}
}}

\def\legfboth{
\Picture{
\qbezier[20](0,-3)(0,0)(0,3)
\put(-0.25,-4.5){$x_1$}
\put(0,3.25){$x_1$}
}}

\def\legfbothb{
\Picture{
\put(0,-1){\vector(0,-1){0.01}}
\put(0.5,1){$f\overline{f}$}
\qbezier[20](0,-3)(0,0)(0,3)
\put(-0.25,-4.5){$x_1$}
\put(0,3.25){$x_1$}
}}

\def\legfb{
\Picture{
\put(-1.5,-0.75){\line(1,0){3}}
\put(-1.5,0.75){\line(1,0){3}}
\put(-1.5,-0.75){\line(0,1){1.5}}
\put(1.5,-0.75){\line(0,1){1.5}}
\qbezier[12](0,0.75)(0,1.9)(0,3)
\qbezier[12](0,-0.75)(0,-1.9)(0,-3)
\put(-0.25,3.5){$x_1$}
\put(-0.8,-0.2){$f(k)$}
}}

\def\bijchord{
\ \ 
\Picture{
\qbezier[16](0,-1)(1,0.25)(2,1.5)
\put(-0.5,-1){$i$}
\put(2.25,1.5){$j$}
}\ \ }

\def\ijchord{
\ \ \Picture{
\qbezier[24](0.5,-0.25)(0.5,1.5)(1,1.5)
\qbezier[24](1.5,-0.25)(1.5,1.5)(1,1.5)
\put(-0.1,-1){$x_i$}
\put(1.27,-1){$x_j$}
}\ \ 
}

\def\izchord{
\ \ \Picture{
\qbezier[24](0.5,-0.25)(0.5,1.5)(1,1.5)
\qbezier[24](1.5,-0.25)(1.5,1.5)(1,1.5)
\put(-0.1,-1){$x_i$}
\put(1.27,-1){$z$}
\put(0.7,2){$L_i$}
}\ \ 
}

\def\ijchordxy{
\ \ \Picture{
\qbezier[24](0.5,-0.25)(0.5,1.5)(1,1.5)
\qbezier[24](1.5,-0.25)(1.5,1.5)(1,1.5)
\put(-0.1,-1){$x_i$}
\put(1.27,-1){$y_j$}
}\ \ 
}

\def\ijchordyy{
\ \ \Picture{
\qbezier[24](0.5,-0.25)(0.5,1.5)(1,1.5)
\qbezier[24](1.5,-0.25)(1.5,1.5)(1,1.5)
\put(-0.1,-1){$y_i$}
\put(1.27,-1){$y_j$}
}\ \ 
}

\def\iichord{
\ \ \Picture{
\qbezier[24](0.5,-0.25)(0.5,1.5)(1,1.5)
\qbezier[24](1.5,-0.25)(1.5,1.5)(1,1.5)
\put(-0.1,-1){$x_i$}
\put(1.27,-1){$x_i'$}
}\ \ 
}

\def\iiichord{
\ \ \Picture{
\qbezier[24](0.5,-0.25)(0.5,1.5)(1,1.5)
\qbezier[24](1.5,-0.25)(1.5,1.5)(1,1.5)
\put(-0.1,-1){$x_i$}
\put(1.27,-1){$x_i$}
}\ \ 
}

\def\ijchordup{
\ \ \Picture{
\qbezier[24](0.5,0.75)(0.5,-0.25)(1,-0.25)
\qbezier[24](1.5,0.75)(1.5,-0.25)(1,-0.25)
\put(-0.1,-1){$x_i$}
\put(1.27,-1){$x_j$}
}\ \ }

\def\corolpict{
\Picture{
\put(-4,-3){\line(0,1){2}}
\put(-4,1){\vector(0,1){2}}
\put(-5,-1){\line(0,1){2}}
\put(-3,-1){\line(0,1){2}}
\put(-4.3,-0.25){$e^k$}
\qbezier(-5,-1)(-4,-1)(-3,-1)
\qbezier(-5,1)(-4,1)(-3,1)
\put(0,-3){\line(0,1){2}}
\put(0,1){\vector(0,1){2}}
\put(-1,-1){\line(0,1){2}}
\put(1,-1){\line(0,1){2}}
\put(-0.3,-0.25){$e^k$}
\qbezier(-1,-1)(0,-1)(1,-1)
\qbezier(-1,1)(0,1)(1,1)
\put(3,-1){\vector(0,-1){2}}
\put(3,1){\line(0,1){2}}
\put(2,-1){\line(0,1){2}}
\put(4,-1){\line(0,1){2}}
\put(2.7,-0.25){$e^k$}
\qbezier(2,-1)(3,-1)(4,-1)
\qbezier(2,1)(3,1)(4,1)
\put(7,-1){\vector(0,-1){2}}
\put(7,1){\line(0,1){2}}
\put(6,-1){\line(0,1){2}}
\put(8,-1){\line(0,1){2}}
\put(6.7,-0.25){$e^k$}
\qbezier(6,-1)(7,-1)(8,-1)
\qbezier(6,1)(7,1)(8,1)
\qbezier[6](-3,2)(-2,2)(-1,2)
\qbezier[6](4,2)(5,2)(6,2)
}}

\def\corolpicts{
\Picture{
\put(-4,-1){\vector(0,1){2}}
\put(-1,-1){\vector(0,1){2}}
\put(1,1){\vector(0,-1){2}}
\put(4,1){\vector(0,-1){2}}
\qbezier[5](-3.5,-0.5)(-2.5,-0.5)(-1.5,-0.5)
\qbezier[5](3.5,-0.5)(2.5,-0.5)(1.5,-0.5)
\put(-4,-1.5){$\underbrace{\ \ \ \ \ \ \ \ \ }_r$}
\put(1,-1.5){$\underbrace{\ \ \ \ \ \ \  \ \ }_s$}
\put(-3.8,0.5){$e^{u}$}
\put(-0.8,0.5){$e^{u}$}
\put(1.2,0.5){$e^{-u}$}
\put(4.2,0.5){$e^{-u}$}
}}

\def\corolpic{
\Picture{
\qbezier(-4,-0.5)(-4.5,-0.5)(-4.5,0)
\qbezier(-4,0.5)(-4.5,0.5)(-4.5,0)
\qbezier(-4,0.5)(-4,0.5)(0,0.5)
\put(0,0.5){\vector(1,0){0.01}}
\qbezier(4,-0.5)(4.5,-0.5)(4.5,0)
\qbezier(4,0.5)(4.5,0.5)(4.5,0)
\qbezier(4,0.5)(4,0.5)(0,0.5)
\qbezier(-3.5,0)(-3.5,-2)(-3.5,-4)
\qbezier(-3.5,1)(-3.5,2)(-3.5,4)
\qbezier(3.5,0)(3.5,-2)(3.5,-4)
\qbezier(3.5,1)(3.5,2)(3.5,4)
\qbezier(0.5,0)(0.5,-2)(0.5,-4)
\qbezier(0.5,1)(0.5,2)(0.5,4)
\qbezier(-0.5,0)(-0.5,-2)(-0.5,-4)
\qbezier(-0.5,1)(-0.5,2)(-0.5,4)
\qbezier[12](-4,-0.5)(0,-0.5)(4,-0.5)
\put(-3.5,4){\vector(0,1){0.01}}
\put(-0.5,4){\vector(0,1){0.01}}
\put(3.5,-4){\vector(0,-1){0.01}}
\put(0.5,-4){\vector(0,-1){0.01}}
\qbezier[4](-3,-2)(-2,-2)(-1,-2)
\qbezier[4](3,-2)(2,-2)(1,-2)
\put(-3.5,-5){$\underbrace{ \ \ \ \ \ \ \ \ \ }_{r}$}
\put(0.5,-5){$\underbrace{ \ \ \ \ \ \ \ \ \ }_{s}$}
}}

\def\ijchordWh{
\Picture{
\qbezier[16](0.8,1.5)(0.8,0)(0.8,-1.5)
\put(1.1,1.5){$x$}
\put(1.1,-1.75){$k$}
}}

\def\ijchordWha{
\Picture{
\qbezier[24](0.5,-0.25)(0.5,1.5)(1,1.5)
\qbezier[24](1.5,-0.25)(1.5,1.5)(1,1.5)
\put(0,-1.4){$x_1$}
\put(1.25,-1.4){$k$}
}}

\def\ijchordWhrs{
\Picture{
\qbezier[24](0.5,-0.25)(0.5,1.5)(1,1.5)
\qbezier[24](1.5,-0.25)(1.5,1.5)(1,1.5)
\put(-0.75,-1.4){$x_{r+s}$}
\put(1.5,-1.4){$k$}
}}

\def\ijchordWhr{
\Picture{
\qbezier[24](0.5,-0.25)(0.5,1.5)(1,1.5)
\qbezier[24](1.5,-0.25)(1.5,1.5)(1,1.5)
\put(0,-1.4){$x_r$}
\put(1.25,-1.4){$k$}
}}

\def\ijchordWhrp{
\Picture{
\qbezier[24](0.5,-0.25)(0.5,1.5)(1,1.5)
\qbezier[24](1.5,-0.25)(1.5,1.5)(1,1.5)
\put(-0.75,-1.4){$x_{r+1}$}
\put(1.5,-1.4){$k$}
}}

\def\hopfundone{
\Picture{
\put(0.5,2.25){$x$}
\qbezier(0,0.25)(0,-1)(0,-3)
\put(0,2){\vector(0,1){0.01}}
\put(1.25,0.75){$k$}
\put(2,-0.2){\vector(0,-1){0.4}}
\qbezier(0.25,-0.5)(2,-0.5)(2,0)
\qbezier(0,0.5)(2,0.5)(2,0)
\qbezier(-0.25,-0.5)(-2,-0.5)(-2,0)
\qbezier(0,0.5)(-2,0.5)(-2,0)
\qbezier(0,0.75)(0,1)(0,3)
}}

\def\ijchordW{
\ \ \Picture{
\qbezier[24](0.5,-0.25)(0.5,1.5)(1,1.5)
\qbezier[24](1.5,-0.25)(1.5,1.5)(1,1.5)
\put(-0.1,-1){$x_i$}
\put(1.27,-1){$x_j$}
\put(1.5,-0.25){\vector(0,-1){0.01}}
\put(-0.9,2){$W_{ij}(T,t)$}
}\ \ }

\def\ijchordWhh{
\ \ \ \ \Picture{
\qbezier[24](0.5,-0.25)(0.5,1.5)(1,1.5)
\qbezier[24](1.5,-0.25)(1.5,1.5)(1,1.5)
\put(-0.1,-1){$x_i$}
\put(1.27,-1){$x_j$}
\put(1.5,-0.25){\vector(0,-1){0.01}}
\put(-1.4,2){$W_{ij}(T,e^{k})$}
}\ \ \ \ }

\def\ijchordM{
\ \ \ \ \Picture{
\qbezier[24](0.5,-0.25)(0.5,1.5)(1,1.5)
\qbezier[24](1.5,-0.25)(1.5,1.5)(1,1.5)
\put(-0.1,-1){$x_i$}
\put(1.27,-1){$x_j$}
\put(1.5,-0.25){\vector(0,-1){0.01}}
\put(-0.75,2){$M_{ij}(e^{u})$}
}\ \ \ \ }

\def\cijchordM{
\ \ \ \ \Picture{
\qbezier[24](0.5,-0.25)(0.5,1.5)(1,1.5)
\qbezier[24](1.5,-0.25)(1.5,1.5)(1,1.5)
\put(-0.1,-1){$x_i$}
\put(1.27,-1){$x_j$}
\put(1.5,-0.25){\vector(0,-1){0.01}}
\put(-0.75,2){$\widetilde{M}_{ij}(e^{u})$}
}\ \ \ \ }

\def\ijchordWo{
\ \ \ \ \Picture{
\qbezier[24](0.5,-0.25)(0.5,1.5)(1,1.5)
\qbezier[24](1.5,-0.25)(1.5,1.5)(1,1.5)
\put(-0.1,-1){$x_i$}
\put(1.27,-1){$x_j$}
\put(1.5,-0.25){\vector(0,-1){0.01}}
\put(-1.4,2){$W_{ij}(T,1)$}
}\ \ \ \ }

\def\ijchordWinv{
\ \ \ \ \Picture{
\qbezier[24](0.5,-0.25)(0.5,1.5)(1,1.5)
\qbezier[24](1.5,-0.25)(1.5,1.5)(1,1.5)
\put(-0.1,-1){$x_i$}
\put(1.27,-1){$x_j$}
\put(1.5,-0.25){\vector(0,-1){0.01}}
\put(-1.6,2){$W^{-1}_{ij}(T,e^k)$}
}\ \ \ \ }

\def\ijchordWinvp{
\ \ \ \ \Picture{
\qbezier[24](0.5,-0.25)(0.5,1.5)(1,1.5)
\qbezier[24](1.5,-0.25)(1.5,1.5)(1,1.5)
\put(-0.1,-1){$x_i'$}
\put(1.27,-1){$x_j'$}
\put(1.5,-0.25){\vector(0,-1){0.01}}
\put(-1.7,2){$W^{-1}_{ij}(T,e^k)$}
}\ \ \ \ }

\def\ijchordWt{
\ \ \Picture{
\qbezier[24](0.5,-0.25)(0.5,1.5)(1,1.5)
\qbezier[24](1.5,-0.25)(1.5,1.5)(1,1.5)
\put(-0.1,-1){$x_i$}
\put(1.27,-1){$x_j$}
\put(1.5,-0.25){\vector(0,-1){0.01}}
\put(-0.5,2){$\widetilde{W}_{ij}(T,t)$}
}\ \ }

\def\ASR{
\Picture{
\qbezier[24](0,-0.5)(1,1)(2,2)
\qbezier[24](0,-0.5)(-1,1)(-2,2)
\qbezier[12](0,-0.5)(0,-1)(0,-2)}}

\def\ASP{
\Picture{
\qbezier[24](0,-0.5)(1,1)(2,2)
\qbezier[24](0,-0.5)(-1,1)(-2,2)
\qbezier[12](0,-0.5)(0,-1)(0,-2)}}

\def\ASBP{
\Picture{
\qbezier[10](0,-0.5)(1,-0.0)(1,0.5)
\qbezier[10](0,-0.5)(-1,-0.0)(-1,0.5)
\qbezier[18](1,0.5)(1,1)(-2,2)
\qbezier[18](-1,0.5)(-1,1)(2,2)
\qbezier[12](0,-0.5)(0,-1)(0,-2)}}

\def\ASRB{
\Picture{
\qbezier[10](0,-0.5)(1,-0.0)(1,0.5)
\qbezier[10](0,-0.5)(-1,-0.0)(-1,0.5)
\qbezier[18](1,0.5)(1,1)(-2,2)
\qbezier[18](-1,0.5)(-1,1)(2,2)
\qbezier[12](0,-0.5)(0,-1)(0,-2)}}

\def\IHX{
\Picture{
\qbezier[30](-2,2)(0,0)(3,-3)
\qbezier[16](-1,1)(0,1)(3,1)
\qbezier[16](0,0)(0,-1)(0,-3)
}}

\def\IHXB{
\Picture{
\qbezier[30](-2,2)(0,0)(3,-3)
\qbezier[16](0,-1.5)(0,-1.5)(4,-1.5)
\qbezier[16](0,0)(0,-1)(0,-3)
}}

\def\IHXC{
\Picture{
\qbezier[30](-2,2)(0,0)(3,-3)
\qbezier[16](1.5,-1.5)(1.5,-1.5)(5.5,-1.5)
\qbezier[16](0,0)(0,-1)(0,-3)
}}

\def\OR{
\Picture{
\qbezier[30](-3,-1)(0,0)(3,1)
\put(-1.5,-0.5){\vector(3,1){0.01}}
\put(-0.5,0.5){$a$}
}
}

\def\ORB{
\Picture{
\qbezier[30](-3,-1)(0,0)(3,1)
\put(-1.5,-0.5){\vector(-3,-1){0.01}}
\put(-0.5,0.5){$\bar{a}$}
}
}

\def\PSA{
\Picture{
\thicklines
\qbezier(-3,0)(0,0)(0,0)
\qbezier(0,0)(1.5,0)(3,0)
\thinlines
\qbezier[20](0,0)(0,1.5)(0,3)
\put(-1.5,-0){\vector(1,0){0.01}}
\put(-2.5,-1){$at$}
\put(0,1.5){\vector(0,1){0.01}}
\put(1.5,0){\vector(1,0){0.01}}
\put(0.5,1.5){$b$}
\put(1.5,-1){$c$}
}
}

\def\PSX{
\Picture{
\qbezier[20](0,0)(0,1.5)(0,3)
\put(0.5,1.5){$at^{\pm 1}$}
\put(2.5,-2){$b$}
\put(-2,-2){\vector(0,1){0.01}}
\qbezier(0,0)(2,0)(2,-2)
\qbezier(0,-4)(2,-4)(2,-2)
\qbezier(0,0)(-2,0)(-2,-2)
\qbezier(0,-4)(-2,-4)(-2,-2)
}
}

\def\PSY{
\Picture{
\qbezier[20](0,0)(0,1.5)(0,3)
\put(0.5,1.5){$a$}
\put(2.5,-2){$b$}
\put(-2,-2){\vector(0,1){0.01}}
\qbezier(0,0)(2,0)(2,-2)
\qbezier(0,-4)(2,-4)(2,-2)
\qbezier(0,0)(-2,0)(-2,-2)
\qbezier(0,-4)(-2,-4)(-2,-2)
}
}

\def\PSB{
\Picture{
\thicklines
\qbezier(-3,0)(0,0)(0,0)
\qbezier(0,0)(1.5,0)(3,0)
\thinlines
\qbezier[20](0,0)(0,1.5)(0,3)
\put(-1.5,-0){\vector(1,0){0.01}}
\put(-2.5,-1){$a$}
\put(0,1.5){\vector(0,1){0.01}}
\put(1.5,0){\vector(1,0){0.01}}
\put(0.5,1.5){$bt$}
\put(1.5,-1){$ct$}
}
}

\def\STUA{
\Picture{
\put(0,-2){$
\Picture{
\put(-2,0){\vector(1,0){0.01}}
\qbezier(-3,0)(0,0)(0,0)
\qbezier(0,0)(1.5,0)(3,0)
\qbezier[16](0,0)(0,1)(0,2)
\qbezier[18](0,2)(1.5,3.5)(3,5)
\qbezier[18](0,2)(-1.5,3.5)(-3,5)
}$
}}
}

\def\STUB{
\Picture{
\put(0,-2){$
\Picture{
\put(-2,0){\vector(1,0){0.01}}
\qbezier(-3,0)(0,0)(0,0)
\qbezier(0,0)(1.5,0)(3,0)
\qbezier[24](1,0)(1,4.5)(3,5)
\qbezier[24](-1,0)(-1,4.5)(-3,5)
}$
}}
}

\def\STUC{
\Picture{
\put(0,-2){$
\Picture{
\put(-2,0){\vector(1,0){0.01}}
\qbezier(-3,0)(0,0)(0,0)
\qbezier(0,0)(1.5,0)(3,0)
\qbezier[24](-1,0)(1,4.5)(3,5)
\qbezier[24](1,0)(-1,4.5)(-3,5)
}$
}}
}

\def\PA{
\Picture{
\qbezier[24](-3,-1)(0,0)(0,0)
\qbezier[20](0,0)(1.5,0)(3,0)
\qbezier[20](0,0)(0.5,1.5)(1,3)
\put(-1.5,-0.5){\vector(3,1){0.01}}
\put(-2.5,-2){$at$}
\put(0.5,1.5){\vector(1,3){0.01}}
\put(1.5,0){\vector(1,0){0.01}}
\put(1,1.5){$b$}
\put(1.5,-1){$c$}
}
}

\def\PAX{
\Picture{
\qbezier[24](-3,-1)(0,0)(0,0)
\qbezier[20](0,0)(3,0)(3,1)
\put(-1.5,-0.5){\vector(3,1){0.01}}
\put(-2,-1.5){$at^{\pm 1}$}
\qbezier[20](0,0)(0,3)(1,3)
\qbezier[20](1,3)(3,3)(3,1)
\put(3,1){\vector(0,1){0.01}}
\put(3,3){$b$}
}
}

\def\PAY{
\Picture{
\qbezier[24](-3,-1)(0,0)(0,0)
\qbezier[20](0,0)(3,0)(3,1)
\put(-1.5,-0.5){\vector(3,1){0.01}}
\put(-2,-1.5){$a$}
\qbezier[20](0,0)(0,3)(1,3)
\qbezier[20](1,3)(3,3)(3,1)
\put(3,1){\vector(0,1){0.01}}
\put(3,3){$b$}
}
}

\def\PAP{
\Picture{
\qbezier[24](-3,-1)(0,0)(0,0)
\qbezier[20](0,0)(1.5,0)(3,0)
\qbezier[20](0,0)(0.5,1.5)(1,3)
\put(-1.5,-0.5){\vector(3,1){0.01}}
\put(-2.5,-2){$at$}
\put(0.5,1.5){\vector(-1,-3){0.01}}
\put(1.5,0){\vector(1,0){0.01}}
\put(1,1.5){$b$}
\put(1.5,-1){$c$}
}
}

\def\PAD{
\Picture{
\qbezier[24](-3,-1)(0,0)(0,0)
\qbezier[20](0,0)(1.5,0)(3,0)
\qbezier[20](0,0)(0.5,1.5)(1,3)
\put(-1.5,-0.5){\vector(3,1){0.01}}
\put(-2.5,-2){$t$}
\put(0.5,1.5){\vector(1,3){0.01}}
\put(1.5,0){\vector(1,0){0.01}}
\put(1,1.5){$1$}
\put(1.5,-1){$1$}
}
}

\def\relchecka{
\Picture{
\qbezier[24](-3,-1)(0,0)(0,0)
\qbezier[20](0,0)(1.5,0)(3,0)
\qbezier[20](0,0)(0.5,1.5)(1,3)
\qbezier[12](-2.5,-0.84)(-2.18,-1.84)(-1.86,-2.84)
\qbezier[12](-0.5,-0.16)(-0.18,-1.16)(0.14,-2.16)
\put(-1.9,-3.5){$u$}
\put(0.1,-2.7){$u$}
\qbezier[4](-1.5,-2.1)(-1,-1.95)(-0.5,-1.8)
\put(-1.5,-4){$\underbrace{\ \ \ \ }_{r}$}
}
}

\def\relcheckb{
\Picture{
\qbezier[24](-3,-1)(0,0)(0,0)
\qbezier[20](0,0)(1.5,0)(3,0)
\qbezier[20](0,0)(0.5,1.5)(1,3)
\qbezier[12](0.16,0.5)(1.16,0.5)(2.16,0.5)
\qbezier[12](0.84,2.5)(1.84,2.5)(2.84,2.5)
\put(2.95,2.25){$u$}
\put(2.25,0.25){$u$}
\qbezier[4](2,1)(2.16,1.5)(2.33,2)
\put(3,1.4){$
\left.
\begin{array}{l}
\\
\\
\end{array}
\right\}$}
\put(5,1.5){$s$}
\qbezier[12](0.5,0)(0.5,-1)(0.5,-2)
\qbezier[12](2.5,0)(2.5,-1)(2.5,-2)
\put(0.25,-2.75){$u$}
\put(2.25,-2.75){$u$}
\qbezier[4](0.5,-1.5)(1.25,-1.5)(2.25,-1.5)
\put(0.25,-3.25){$\underbrace{\ \ \ \ \ \ \ \ }_{r-s}$}
}
}

\def\PB{
\Picture{
\qbezier[24](-3,-1)(0,0)(0,0)
\qbezier[20](0,0)(0.5,1.5)(1,3)
\qbezier[20](0,0)(1.5,0)(3,0)
\put(-1.5,-0.5){\vector(3,1){0.01}}
\put(-2.5,-2){$a$}
\put(0.5,1.5){\vector(1,3){0.01}}
\put(1.5,0){\vector(1,0){0.01}}
\put(1,1.5){$bt$}
\put(1.5,-1){$ct$}
}
}

\def\PBP{
\Picture{
\qbezier[24](-3,-1)(0,0)(0,0)
\qbezier[20](0,0)(0.5,1.5)(1,3)
\qbezier[20](0,0)(1.5,0)(3,0)
\put(-1.5,-0.5){\vector(3,1){0.01}}
\put(-2.5,-2){$a$}
\put(0.5,1.5){\vector(-1,-3){0.01}}
\put(1.5,0){\vector(1,0){0.01}}
\put(1,1.5){$bt^{-1}$}
\put(1.5,-1){$ct$}
}
}

\def\PBD{
\Picture{
\qbezier[24](-3,-1)(0,0)(0,0)
\qbezier[20](0,0)(0.5,1.5)(1,3)
\qbezier[20](0,0)(1.5,0)(3,0)
\put(-1.5,-0.5){\vector(3,1){0.01}}
\put(-2.5,-2){$1$}
\put(0.5,1.5){\vector(1,3){0.01}}
\put(1.5,0){\vector(1,0){0.01}}
\put(1,1.5){$t$}
\put(1.5,-1){$t$}
}
}

\def\AA{
\Picture{
\qbezier[30](-3,-1)(0,0)(3,1)
\put(-1.5,-0.5){\vector(3,1){0.01}}
\put(-3.5,0.5){$q_1a+q_2b$}
}
}

\def\AB{
\Picture{
\qbezier[30](-3,-1)(0,0)(3,1)
\put(-1.5,-0.5){\vector(3,1){0.01}}
\put(-0.5,0.5){$a$}
}
}

\def\AC{
\Picture{
\qbezier[30](-3,-1)(0,0)(3,1)
\put(-1.5,-0.5){\vector(3,1){0.01}}
\put(-0.5,0.5){$b$}
}
}

\def\boxedline{
\Picture{
\put(0.75,-0.55){\line(0,1){1.1}}
\put(-0.75,-0.55){\line(1,0){1.5}}
\put(-0.75,-0.55){\line(0,1){1.1}}
\put(-0.75,0.55){\line(1,0){1.5}}
\put(-0.3,-0.2){$f$}
\qbezier[16](-2.75,0)(-1.75,0)(-0.75,0)
\qbezier[16](0.75,0)(1.75,0)(2.75,0)
\put(2.75,0){\vector(1,0){0.01}}
}
}

\def\boxedlinea{
\Picture{
\put(-0.25,0.5){$f$}
\qbezier[28](-2.75,0)(0,0)(2.75,0)
\put(2,0){\vector(1,0){0.01}}
}
}
\def\boxedlineb{
\Picture{
\put(-0.25,0.5){$t^n$}
\qbezier[28](-2.75,0)(0,0)(2.75,0)
\put(2,0){\vector(1,0){0.01}}
}
}
\def\boxedlinec{
\Picture{
\put(-0.25,0.5){$e^{nu}$}
\qbezier[28](-2.75,0)(0,0)(2.75,0)
\put(2,0){\vector(1,0){0.01}}
}
}

\def\terma{
\Picture{
\qbezier[20](-1.5,0)(0.25,0)(2.5,0)}}

\def\termb{
\Picture{
\qbezier[20](-1.5,0)(0.25,0)(2.5,0)
\qbezier[12](0.25,0)(0.25,-1)(0.25,-2)
\put(0,-3){$k$}
}
}
\def\termc{
\Picture{
\qbezier[20](-1.5,0)(0.25,0)(2.5,0)
\qbezier[12](0,0)(0,-1)(0,-2)
\qbezier[12](0.5,0)(0.5,-1)(0.5,-2)
\put(-0.25,-3){$k$}
\put(0.25,-3){$k$}
}}
\def\termd{
\Picture{
\qbezier[20](-1.5,0)(0.25,0)(2.5,0)
\qbezier[12](0.75,0)(0.75,-1)(0.75,-2)
\qbezier[12](0.25,0)(0.25,-1)(0.25,-2)
\qbezier[12](-0.25,0)(-0.25,-1)(-0.25,-2)
\put(-0.5,-3){$k$}
\put(0,-3){$k$}
\put(0.5,-3){$k$}
}}

\def\clasp{
\Picture{
\qbezier(0,1)(0,0.6)(0,0.6)
\qbezier(0,-1)(0,-0.6)(0,0.2)
\qbezier(0,0.4)(0.4,0.4)(0.4,0)
\qbezier(0,0.4)(-0.4,0.4)(-0.4,0)
\qbezier(0.4,0)(0.4,-0.2)(0.2,-0.36)
\qbezier(-0.4,0)(-0.4,-0.2)(-0.2,-0.36)
\qbezier(0.4,0)(0.6,0)(1,0)
}
}

\def\claspf{
\Picture{
\qbezier(0,1)(0,0)(0,-1)
\put(1,0){$\claspc$}
}
}

\def\claspff{
\Picture{
\qbezier[5](0,1)(0,0.5)(0,0)
\qbezier(0,0)(0,-0.5)(0,-1)
}
}

\def\claspffp{
\Picture{
\qbezier(0,0)(0,-0.5)(0,-1)
}
}

\def\claspfft{
\Picture{
\qbezier(0,1)(0,0.5)(0,0)
\qbezier[5](0,0)(0,-0.5)(0,-1)
}
}

\def\claspfftp{
\Picture{
\qbezier(0,1)(0,0.5)(0,0)
}
}

\def\claspc{
\Picture{
\qbezier(0,0.4)(0.4,0.4)(0.4,0)
\qbezier(0,0.4)(-0.4,0.4)(-0.4,0)
\qbezier(0.4,0)(0.4,-0.4)(0,-0.4)
\qbezier(-0.4,0)(-0.4,-0.4)(-0,-0.4)
\qbezier(0.4,0)(0.6,0)(1,0)
}
}

\def\claspr{
\Picture{
\qbezier(-0,-1)(-0,-0.6)(-0,-0.6)
\qbezier(-0,--1)(-0,--0.6)(-0,-0.2)
\qbezier(-0,-0.4)(-0.4,-0.4)(-0.4,-0)
\qbezier(-0,-0.4)(--0.4,-0.4)(--0.4,-0)
\qbezier(-0.4,-0)(-0.4,--0.2)(-0.2,--0.36)
\qbezier(--0.4,-0)(--0.4,--0.2)(--0.2,--0.36)
\qbezier(-0.4,-0)(-0.6,-0)(-1,-0)
}
}

\def\capl{
\Picture{
\qbezier(0,-1)(0,0)(-1,0)
\qbezier(-2,-1)(-2,0)(-1,0)
}}

\def\cupl{
\Picture{
\qbezier(0,1)(0,0)(-1,0)
\qbezier(-2,1)(-2,0)(-1,0)
}}

\def\caplf{
\Picture{
\qbezier[5](-2,-1)(-2,-1.5)(-2,-2)
\qbezier[10](0,-1)(0,0)(-1,0)
\qbezier[10](-2,-1)(-2,0)(-1,0)
}}

\def\cuplf{
\Picture{
\qbezier[5](-2,1)(-2,1.5)(-2,2)
\qbezier[10](0,1)(0,0)(-1,0)
\qbezier[10](-2,1)(-2,0)(-1,0)
}}

\def\leftbox{
\Picture{
\put(1,0){$\clasp$}
\put(1,-2){$\clasp$}
\put(1,2){$\capl$}
\put(1,-4){$\cupl$}
\put(-1,-3){\line(0,1){4}}
}}

\def\leftboxn{
\Picture{
\put(1,0){$\claspff$}
\put(1,-2){$\claspfft$}
\put(1,2){$\caplf$}
\put(1,-4){$\cuplf$}
\put(-1,-2){\line(0,1){2}}
}}

\def\leftboxnb{
\Picture{
\put(1,0){$\claspff$}
\put(1,-2){$\claspfftp$}
\put(1,2){$\caplf$}
\put(-1,-2){\line(0,1){2}}
\qbezier(-1,-2)(-1,-3)(0,-3)
\qbezier(1,-2)(1,-3)(0,-3)
}}

\def\leftboxnc{
\Picture{
\put(1,0){$\claspffp$}
\put(1,-2){$\claspfft$}
\put(1,-4){$\cuplf$}
\put(-1,-2){\line(0,1){2}}
\qbezier(-1,0)(-1,1)(0,1)
\qbezier(1,0)(1,1)(0,1)
}}

\def\reltripa{
\Picture{
\put(-9,5){\circle*{0.7}}
\put(-4,4){$\leftboxn$}
\qbezier(-5,3)(-7,3)(-9,5)
\put(-9,5){\line(0,1){2}}
\put(-9,5){\line(-1,0){2}}
\put(-6,-2){$(M,L_A)$}}
}

\def\reltripb{
\Picture{
\put(-9,5){\circle*{0.7}}
\put(-4,4){$\leftboxnb$}
\qbezier(-5,3)(-7,3)(-9,5)
\put(-9,5){\line(0,1){2}}
\put(-9,5){\line(-1,0){2}}
\put(-6,-2){$(M,L_B)$}
}
}

\def\reltripc{
\Picture{
\put(-9,5){\circle*{0.7}}
\put(-4,4){$\leftboxnc$}
\qbezier(-5,3)(-7,3)(-9,5)
\put(-9,5){\line(0,1){2}}
\put(-9,5){\line(-1,0){2}}
\put(-6,-2){$(M,L_C)$}
}
}

\def\leftboxfa{
\Picture{
\put(1,0){$\claspf$}
\put(1,-2){$\clasp$}
\put(1,2){$\capl$}
\put(1,-4){$\cupl$}
\put(-1,-3){\line(0,1){4}}
}}

\def\leftboxfb{
\Picture{
\put(1,0){$\clasp$}
\put(1,-2){$\claspf$}
\put(1,2){$\capl$}
\put(1,-4){$\cupl$}
\put(-1,-3){\line(0,1){4}}
}}

\def\leftboxf{
\Picture{
\put(1,0){$\claspc$}
\put(1,-2){$\claspc$}
}}

\def\leftboxp{
\Picture{
\put(1,0){$\clasp$}
\put(1,-2){$\clasp$}
\put(1,2){$\capl$}
\put(1,-4){$\cupl$}
\qbezier(-1,-3)(-3,-3)(-5,-1)
}}


\def\capr{
\Picture{
\qbezier(-0,--1)(-0,-0)(--1,-0)
\qbezier(--2,--1)(--2,-0)(--1,-0)
}}

\def\cupr{
\Picture{
\qbezier(-0,-1)(-0,-0)(--1,-0)
\qbezier(--2,-1)(--2,-0)(--1,-0)
}}

\def\rightbox{
\Picture{
\put(-1,-0){$\claspr$}
\put(-1,--2){$\claspr$}
\put(-1,-2){$\capr$}
\put(-1,--4){$\cupr$}
\put(1,3){\line(-0,-1){4}}
}}

\def\vertclasp{
\Picture{
\qbezier(-1,0)(-0.8,0)(-0.6,0)
\qbezier(--1,0)(--0.6,-0)(-0.2,-0)
\qbezier(-0.4,0)(-0.4,0.4)(-0,0.4)
\qbezier(-0.4,0)(-0.4,-0.4)(--0,-0.4)
\qbezier(-0,0.4)(--0.2,0.4)(--0.36,0.2)
\qbezier(-0,-0.4)(--0.2,-0.4)(--0.36,-0.2)
\qbezier(-0,0.4)(-0,0.6)(-0,1)
}}

\def\vertclaspud{
\Picture{
\qbezier(--1,0)(--0.8,0)(--0.6,0)
\qbezier(---1,0)(---0.6,-0)(--0.2,-0)
\qbezier(--0.4,0)(--0.4,-0.4)(--0,-0.4)
\qbezier(--0.4,0)(--0.4,--0.4)(---0,--0.4)
\qbezier(--0,-0.4)(---0.2,-0.4)(---0.36,-0.2)
\qbezier(--0,--0.4)(---0.2,--0.4)(---0.36,--0.2)
\qbezier(--0,-0.4)(--0,-0.6)(--0,-1)
}}

\def\ygraph{
\Picture{
\put(0,-3.75){$\vertclaspf$}
\put(-3.45,6){$\vertclaspfu$}
\put(3.45,6){$\vertclaspfu$}
\put(0,0){\circle{1.5}}
\qbezier(0.75,0)(3.75,1)(3.75,3)
\qbezier(-0.75,0)(-3.75,1)(-3.75,3)
\qbezier(-0.45,0.5)(-3.25,1.5)(-3.15,3)
\qbezier(0.45,0.5)(3.25,1.5)(3.15,3)
}}

\def\vertclaspf{
\Picture{
\qbezier(-0.3,1)(-0.3,2)(-0.3,3)
\qbezier(0.3,1)(0.3,2)(0.3,3)
\qbezier(-2,-0.2)(-3,-0.2)(-3,-4)
\qbezier(-2,0.2)(-3.6,0.2)(-3.6,-4)
\qbezier(2,-0.2)(3,-0.2)(3,-4)
\qbezier(2,0.2)(3.6,0.2)(3.6,-4)
\qbezier(1,-0.2)(1,-1)(0,-1)
\qbezier(-1,0)(-1,1)(0,1)
\qbezier(1,0.2)(1,1)(0,1)
\qbezier(-1,0)(-1,-1)(0,-1)
\qbezier(-2,-0.2)(-1,-0.2)(-1,-0.2)
\qbezier(--2,-0.2)(--0.6,-0.2)(-0.36,-0.2)
\qbezier(-2,0.2)(-1,0.2)(-1,0.2)
\qbezier(--2,0.2)(--0.6,0.2)(-0.36,0.2)
\qbezier(-0.4,0)(-0.4,0.4)(-0,0.4)
\qbezier(-0.4,0)(-0.4,-0.4)(--0,-0.4)
\qbezier(-0,0.4)(--0.2,0.4)(--0.36,0.2)
\qbezier(-0,-0.4)(--0.2,-0.4)(--0.36,-0.2)
}}

\def\vertclaspfu{
\Picture{
\qbezier(-0.3,-1)(-0.3,-2)(-0.3,-3)
\qbezier(0.3,-1)(0.3,-2)(0.3,-3)
\qbezier(-2,--0.2)(-3,--0.2)(-3,--4)
\qbezier(-2,-0.2)(-3.6,-0.2)(-3.6,--4)
\qbezier(2,--0.2)(3,--0.2)(3,--4)
\qbezier(2,-0.2)(3.6,-0.2)(3.6,--4)
\qbezier(1,--0.2)(1,--1)(0,--1)
\qbezier(-1,-0)(-1,-1)(0,-1)
\qbezier(1,-0.2)(1,-1)(0,-1)
\qbezier(-1,-0)(-1,--1)(0,--1)
\qbezier(-2,--0.2)(-1,--0.2)(-1,--0.2)
\qbezier(--2,--0.2)(--0.6,--0.2)(-0.36,--0.2)
\qbezier(-2,-0.2)(-1,-0.2)(-1,-0.2)
\qbezier(--2,-0.2)(--0.6,-0.2)(-0.36,-0.2)
\qbezier(-0.4,-0)(-0.4,-0.4)(-0,-0.4)
\qbezier(-0.4,-0)(-0.4,--0.4)(--0,--0.4)
\qbezier(-0,-0.4)(--0.2,-0.4)(--0.36,-0.2)
\qbezier(-0,--0.4)(--0.2,--0.4)(--0.36,--0.2)
}}

\def\extvertclasp{
\Picture{
\put(0,0){$\vertclasp$}
\qbezier(-4,0)(-1,0)(-1,0)
\qbezier(4,0)(1,0)(1,0)
}}

\def\verthook{
\Picture{
\qbezier(-1,0)(-1,-2)(0,-2)
\qbezier(0,-2)(1,-2)(1,-0.2)
\qbezier(-4,0)(-1.25,0)(-1.25,0)
\qbezier(4,0)(1.25,0)(-0.75,0)
\qbezier(-1,0)(-1,1)(-1,2)
\qbezier(1,0.2)(1,1)(1,2)
}}

\def\verthookc{
\Picture{
\qbezier(-1,0)(-1,-2)(0,-2)
\qbezier(0,-2)(1,-2)(1,-0.7)
\qbezier(-4,0)(-1.25,0)(-1.25,0)
\qbezier(-4,-0.5)(-1.25,-0.5)(-1.25,-0.5)
\qbezier(4,0)(1.25,0)(-0.75,0)
\qbezier(4,-0.5)(1.25,-0.5)(-0.75,-0.5)
\qbezier(-1,0)(-1,1)(-1,2)
\qbezier(1,0.2)(1,1)(1,2)
}}

\def\verthookp{
\Picture{
\qbezier(-1,0)(-1,-2)(0,-2)
\qbezier(-4,0)(-1.25,0)(-1.25,0)
\qbezier(5,0)(1.25,0)(-0.75,0)
\qbezier(-1,0)(-1,1)(-1,2)
\qbezier(1,0.2)(1,1)(1,2)
}}

\def\halftwist{
\Picture{
\put(0,0){\circle{1.4}}
\put(0,0.7){\line(0,1){0.3}}
\put(0,-0.7){\line(0,-1){0.3}}
\qbezier(0.49,0.49)(0,0)(-0.49,-0.49)
}
}

\def\halftwistp{
\Picture{
\put(0,0){\circle{1.4}}
\put(0.7,0){\line(1,0){0.3}}
\put(-0.7,0){\line(-1,0){0.3}}
\qbezier(0.49,-0.49)(0,0)(-0.49,0.49)
}
}

\def\shade{
\Picture{
\qbezier[20](0,-3)(-2,-5)(-4,-7)
\qbezier[20](1,-3)(-1,-5)(-3,-7)
\qbezier[20](2,-3)(0,-5)(-2,-7)
\qbezier[20](-2,-3)(-4,-5)(-6,-7)
\qbezier[20](-1,-3)(-3,-5)(-5,-7)
\qbezier[20](-0.2,-3)(-2.2,-5)(-4.2,-7)
\qbezier[20](-0.4,-3)(-2.4,-5)(-4.4,-7)
\qbezier[20](-0.6,-3)(-2.6,-5)(-4.6,-7)
\qbezier[20](-0.8,-3)(-2.8,-5)(-4.8,-7)
\qbezier[20](0.2,-3)(-1.8,-5)(-3.8,-7)
\qbezier[20](0.4,-3)(-1.6,-5)(-3.6,-7)
\qbezier[20](0.6,-3)(-1.4,-5)(-3.4,-7)
\qbezier[20](0.8,-3)(-1.2,-5)(-3.2,-7)
\qbezier[20](-2.2,-3)(-4.2,-5)(-6.2,-7)
\qbezier[20](-2.4,-3)(-4.4,-5)(-6.4,-7)
\qbezier[20](-2.6,-3)(-4.6,-5)(-6.6,-7)
\qbezier[20](-2.8,-3)(-4.8,-5)(-6.8,-7)
\qbezier[20](2.2,-3)(0.2,-5)(-1.8,-7)
\qbezier[20](2.4,-3)(0.4,-5)(-1.6,-7)
\qbezier[20](2.6,-3)(0.6,-5)(-1.4,-7)
\qbezier[20](2.8,-3)(0.8,-5)(-1.2,-7)
\qbezier[20](-1.2,-3)(-3.2,-5)(-5.2,-7)
\qbezier[20](-1.4,-3)(-3.4,-5)(-5.4,-7)
\qbezier[20](-1.6,-3)(-3.6,-5)(-5.6,-7)
\qbezier[20](-1.8,-3)(-3.8,-5)(-5.8,-7)
\qbezier[20](1.2,-3)(-0.8,-5)(-2.8,-7)
\qbezier[20](1.4,-3)(-0.6,-5)(-2.6,-7)
\qbezier[20](1.6,-3)(-0.4,-5)(-2.4,-7)
\qbezier[20](1.8,-3)(-0.2,-5)(-2.2,-7)
}}

\def\shadet{
\Picture{
\qbezier[10](1,-3)(-1,-5)(-3,-7)
\qbezier[10](-1,-3)(-3,-5)(-5,-7)
\qbezier[10](-0.2,-3)(-2.2,-5)(-4.2,-7)
\qbezier[10](-0.6,-3)(-2.6,-5)(-4.6,-7)
\qbezier[10](0.2,-3)(-1.8,-5)(-3.8,-7)
\qbezier[10](0.6,-3)(-1.4,-5)(-3.4,-7)
\qbezier[10](-2.2,-3)(-4.2,-5)(-6.2,-7)
\qbezier[10](-2.6,-3)(-4.6,-5)(-6.6,-7)
\qbezier[10](2.2,-3)(0.2,-5)(-1.8,-7)
\qbezier[10](2.6,-3)(0.6,-5)(-1.4,-7)
\qbezier[10](-1.4,-3)(-3.4,-5)(-5.4,-7)
\qbezier[10](-1.8,-3)(-3.8,-5)(-5.8,-7)
\qbezier[10](1.4,-3)(-0.6,-5)(-2.6,-7)
\qbezier[10](1.8,-3)(-0.2,-5)(-2.2,-7)
}}

\def\leg{
\Picture{
\put(0,-3){$\extvertclasp$}
\qbezier(0,-2)(0,-1)(0,-1)
\put(0,-1){\circle*{0.7}}
\qbezier(0,-1)(2.5,3)(4,3)
\qbezier(0,-1)(-2.5,3)(-4,3)
\put(0,0){$\shade$}
}}

\def\legc{
\Picture{
\qbezier(0,-2)(0,-1)(0,-1)
\put(0,-1){\circle*{0.7}}
\qbezier(0,-1)(2,1)(3,2)
\qbezier(0,-1)(-2,1)(-3,2)
\put(0,-3){$\halftwist$}
\qbezier(0,-4)(0,-6)(0,-6)
\put(-2,-8){$(M,L_A)$}
}}

\def\legcc{
\Picture{
\qbezier(0,-2)(0,-1)(0,-1)
\put(0,-1){\circle*{0.7}}
\qbezier(0,-1)(2,1)(3,2)
\qbezier(0,-1)(-2,1)(-3,2)
\qbezier(0,-2)(0,-6)(0,-6)
\put(-2,-8){$(M,L_B)$}
}}

\def\legccp{
\Picture{
\qbezier(0,-2)(0,-1)(0,-1)
\put(0,-1){\circle*{0.7}}
\qbezier(0,-1)(2,1)(4,3)
\qbezier(0,-1)(-2,1)(-4,3)
\qbezier(0,-2)(0,-6)(0,-9)
}}

\def\legp{
\Picture{
\put(0,-3){$\extvertclasp$}
\qbezier(0,-2)(0,-1)(0,-1)
\put(0,-1){\circle*{0.7}}
\qbezier(0,-1)(2,1)(3,2)
\qbezier(0,-1)(-2,1)(-3,2)
\put(-3,2){\circle*{0.7}}
\qbezier(-3,2)(-3,3)(-3,5)
\qbezier(-3,2)(-4,2)(-6,2)
\put(-1.5,-6.5){$(K,L)$}
}}

\def\legppb{
\Picture{
\put(0,-3){$\extvertclasp$}
\put(4,-4){$K$}
\qbezier(0,-2)(0,-1)(0,-1)
\put(0,-1){\circle*{0.7}}
\qbezier(0,-1)(2,1)(3,2)
\qbezier(0,-1)(-2,1)(-3,2)
}}

\def\leglink{
\Picture{
\put(0,-3){$\extvertclasp$}
\put(4,-4){$K$}
\qbezier(0,-2)(0,-1)(0,2)
}}

\def\polya{
\Picture{
\put(-2.8,2.2){$x_1$}
\put(-4.3,-5){\vector(0,1){4}}
\qbezier(-7,-3)(-5,-3)(-5,-3)
\put(0,0){$\shade$}
\put(0,-3){$\extvertclasp$}
\put(5,-2){$K$}
\put(-2,1){\vector(1,-1){0.01}}
\put(2,1){\vector(1,1){0.01}}
\qbezier(0,-2)(0,-1)(0,-1)
\put(0,-1){\circle*{0.7}}
\qbezier(0,-1)(2,1)(3,2)
\qbezier(0,-1)(-2,1)(-3,2)
\put(-2.8,-10){$1-t^{-x^i}$}
}}

\def\polyb{
\Picture{
\put(-2.8,2.2){$x_1$}
\put(-4.3,-5){\vector(0,1){4}}
\qbezier(-7,-3)(-5,-3)(-5,-3)
\put(0,0){$\shade$}
\put(0,-3){$\extvertclasp$}
\put(5,-2){$K$}
\put(-2,1){\vector(-1,1){0.01}}
\put(2,1){\vector(-1,-1){0.01}}
\qbezier(0,-2)(0,-1)(0,-1)
\put(0,-1){\circle*{0.7}}
\qbezier(0,-1)(2,1)(3,2)
\qbezier(0,-1)(-2,1)(-3,2)
\put(-2.8,-10){$1-t^{x^i}$}
}}

\def\addlegs{
\Picture{
\put(-2.8,2.2){$e$}
\put(-4.3,-5){\vector(0,1){4}}
\qbezier(-7,-3)(-5,-3)(-5,-3)
\put(0,0){$\shade$}
\put(0,-3){$\extvertclasp$}
\put(5,-2){$K$}
\put(-2,1){\vector(-1,1){0.01}}
\put(2,1){\vector(-1,-1){0.01}}
\qbezier(0,-2)(0,-1)(0,-1)
\put(0,-1){\circle*{0.7}}
\qbezier(0,-1)(2,1)(3,2)
\qbezier(0,-1)(-2,1)(-3,2)
\put(-2.8,-10){$1-t^{\hat{e}}$}
}}

\def\legpp{
\Picture{
\put(-3,2){\circle*{0.7}}
\qbezier(-3,2)(-3,3)(-3,5)
\qbezier(-3,2)(-4,2)(-6,2)
\qbezier(-3,2)(0,-1)(3,2)
\qbezier(-4,-3)(0,-3)(4,-3)
\put(-1.5,-6.5){$(K,L')$}
}}

\def\legppc{
\Picture{
\put(-3,2){\circle*{0.7}}
\qbezier(-3,2)(-3,3)(-3,5)
\qbezier(-3,2)(-4,2)(-6,2)
\put(3,2){\circle*{0.7}}
\qbezier(3,2)(3,3)(3,5)
\qbezier(3,2)(4,2)(6,2)
\qbezier(-3,2)(-2,1)(0,1)
\qbezier(3,2)(2,1)(0,1)
\put(-1.5,-6.5){$(M,L_A)$}
}}

\def\yot{
\Picture{
\put(0,0){$\cross$}
\put(-1.6,-0.75){$\crossr$}
\put(-1,1.4){\line(0,1){3}}
\put(-0.6,-1.15){\line(0,-1){3}}
}}

\def\crossh{
\Picture{
\qbezier(1,-1.4)(0.5,-1.4)(0,-1.4)
\qbezier(1,-1.4)(1.4,-1.4)(1.4,-1)
\qbezier(1,-0.6)(1.4,-0.6)(1.4,-1)
\qbezier(0,-1.4)(-0.4,-1.4)(-0.4,-1)
\qbezier(-0.4,-1)(-0.4,-0.6)(0,-0.6)
}}

\def\crosshb{
\Picture{
\qbezier(-1,--1.4)(-0.5,--1.4)(-0,--1.4)
\qbezier(-1,--1.4)(-1.4,--1.4)(-1.4,--1)
\qbezier(-1,--0.6)(-1.4,--0.6)(-1.4,--1)
\qbezier(-0,--1.4)(--0.4,--1.4)(--0.4,--1)
\qbezier(--0.4,--1)(--0.4,--0.6)(-0,--0.6)
}}

\def\brok{
\Picture{
\qbezier(-1.4,0)(-3,0)(-3,0)
\qbezier(3,0.15)(1.2,0.15)(1.2,0.15)
\put(-0.85,0.95){$\crossh$}
\put(0.8,-0.75){$\crosshb$}
\put(-6,2){\circle*{0.7}}
\qbezier(-6,2)(-6,3)(-6,5)
\qbezier(-6,2)(-7,2)(-9,2)
\put(6,2){\circle*{0.7}}
\qbezier(6,2)(6,3)(6,5)
\qbezier(6,2)(7,2)(9,2)
\qbezier(-6,2)(-5,0)(-3,0)
\qbezier(6,2)(5,0.15)(3,0.15)
}}

\def\brokb{
\Picture{
\qbezier(-3,0)(-3,0)(3,0)
\put(-6,2){\circle*{0.7}}
\qbezier(-6,2)(-6,3)(-6,5)
\qbezier(-6,2)(-7,2)(-9,2)
\put(6,2){\circle*{0.7}}
\qbezier(6,2)(6,3)(6,5)
\qbezier(6,2)(7,2)(9,2)
\qbezier(-6,2)(-5,0)(-3,0)
\qbezier(6,2)(5,0)(3,0)
}}

\def\brek{
\Picture{
\qbezier(-1.4,0)(-3,0)(-3,0)
\qbezier(3,0.15)(1.2,0.15)(1.2,0.15)
\put(-0.85,0.95){$\crossh$}
\put(0.8,-0.75){$\crosshb$}
}}

\def\brekb{
\Picture{
\put(-3,0){$\brek$}
\put(3,0.15){$\brek$}
}}

\def\legpppb{
\Picture{
\put(4,-4){$K$}
\qbezier(-3,2)(0,-1)(3,2)
\qbezier(-4,-3)(0,-3)(4,-3)
}}

\def\legppppb{
\Picture{
\qbezier(-3,2)(0,-1)(3,2)
}}

\def\legppp{
\Picture{
\put(-3,2){\circle*{0.7}}
\qbezier(-3,2)(-3,3)(-3,5)
\qbezier(-3,2)(-4,2)(-6,2)
\qbezier(-3,2)(-1,2)(-1,-1)
\put(0,-3){$\verthook$}
\qbezier(1,-1)(1,1)(3,2)
\put(-1.5,-6.5){$(K,L'')$}
}}

\def\legpppc{
\Picture{
\put(-3,2){\circle*{0.7}}
\qbezier(-3,2)(-3,3)(-3,5)
\qbezier(-3,2)(-4,2)(-6,2)
\put(3,2){\circle*{0.7}}
\qbezier(3,2)(3,3)(3,5)
\qbezier(3,2)(4,2)(6,2)
\qbezier(-3,2)(-1,2)(-1,-1)
\put(0,-3){$\verthookc$}
\qbezier(1,-1)(1,1)(3,2)
\put(-1.5,-6.5){$(M,L_B)$}
}}

\def\legged{
\Picture{
\put(-4,4){$\leftbox$}
\put(4,2){$\rightbox$}
\qbezier(-5,3)(-7,3)(-9,3)
\qbezier(5,3)(7,3)(9,3)
\qbezier(-2,4)(0,4)(2,4)
\put(0,-4){$\verthookp$}
\put(1,-1){$\shade$}
\qbezier(-2,2)(-1,2)(-1,0)
\qbezier(2,2)(1,2)(1,0)
\put(-1,-1){$\halftwist$}
\qbezier(1,0)(1,-1)(1,-2)
}}

\def\leggedy{
\Picture{
\put(-4,4){$\leftbox$}
\put(4,2){$\rightbox$}
\put(1,2){$\vertclaspud$}
\put(1,1){\line(0,-1){6}}
\qbezier(-5,3)(-7,3)(-9,7)
\qbezier(5,3)(7,3)(9,7)
\qbezier(-2,4)(0,4)(2,4)
\put(-1,2){$\halftwistp$}
}}

\def\leggedp{
\Picture{
\put(-4,4){$\leftbox$}
\put(4,2){$\rightbox$}
\qbezier(-5,3)(-7,3)(-9,5)
\qbezier(5,3)(7,3)(9,3)
\qbezier(-2,4)(0,4)(2,4)
\put(6,-8.75){$\check{F}$}
\qbezier(-2,2)(-1,2)(-1,0)
\qbezier(2,2)(1,2)(1,0)
\put(-1,-1){$\halftwist$}
\put(-0.75,1.75){B}
\qbezier(1,0)(1,-1)(1,-2)
\put(0,-5){$\modhookp$}
\qbezier(3,-7)(3,-5.5)(2,-4.5)
\qbezier(1,-2)(1,-3.5)(2,-4.5)
\put(3,-4.5){C}
\put(-9,5){\circle*{0.7}}
\put(-8,2){Y}
\put(-0.4,5){A}
\put(-9,5){\line(0,1){2}}
\put(-9,5){\line(-1,0){2}}
}}

\def\tripa{
\Picture{
\put(-9,2){Y$^{a_1}$}
\put(2,3.75){A$^{a_1}$}
\put(2,1.75){B$^{a_1}$}
\put(-9,5){\circle*{0.7}}
\put(-4,4){$\leftbox$}
\qbezier(-5,3)(-7,3)(-9,5)
\qbezier(-2,4)(0,4)(1,4)
\qbezier(-2,2)(-1,2)(1,2)
\put(-9,5){\line(0,1){2}}
\put(-9,5){\line(-1,0){2}}
\put(-6,-2){$\widetilde{L}^{(a_1,\ldots,a_s)}$}}
}

\def\tripb{
\Picture{
\put(-9,2){Y$^{a_1}$}
\put(2,3.75){A$^{a_1}$}
\put(2,1.75){B$^{a_1}$}
\put(-9,5){\circle*{0.7}}
\put(-4,4){$\leftboxfa$}
\qbezier(-5,3)(-7,3)(-9,5)
\qbezier(-1,4)(0,4)(1,4)
\qbezier(-2,2)(-1,2)(1,2)
\put(-9,5){\line(0,1){2}}
\put(-9,5){\line(-1,0){2}}
\put(-6,-2){${R}^{(a_1,\ldots,a_s)}$}
}}

\def\tripc{
\Picture{
\put(-9,2){Y$^{a_1}$}
\put(2,3.75){A$^{a_1}$}
\put(2,1.75){B$^{a_1}$}
\put(-9,5){\circle*{0.7}}
\put(-4,4){$\leftboxfb$}
\qbezier(-5,3)(-7,3)(-9,5)
\qbezier(-2,4)(0,4)(1,4)
\qbezier(-1,2)(-1,2)(1,2)
\put(-9,5){\line(0,1){2}}
\put(-9,5){\line(-1,0){2}}
\put(-6,-2){${S}^{(a_1,\ldots,a_s)}$}
}}

\def\leggedpp{
\Picture{
\put(-4,4){$\leftboxf$}
\put(4,2){$\rightbox$}
\qbezier(5,3)(7,3)(9,3)
\qbezier(-2,4)(0,4)(2,4)
\qbezier(-2,2)(-1,2)(-1,0)
\qbezier(2,2)(1,2)(1,0)
\put(-1,-1){$\halftwist$}
\qbezier(1,0)(1,-1)(1,-2)
\put(0,-5){$\modhookp$}
\qbezier(3,-7)(3,-5.5)(2,-4.5)
\qbezier(1,-2)(1,-3.5)(2,-4.5)
}}

\def\pierce{
\Picture{
\qbezier(-1,5)(-1,-2)(-1,-2)
\qbezier[28](-1,-2)(-1,-4)(-1,-6)
\put(0.5,3){$\shadet$}
\qbezier(-5,0)(-1.2,0)(-1.2,0)
\qbezier(-0.8,0)(5,0)(5,0)
}}

\def\modhook{
\Picture{
\qbezier(-1.4,1)(-1.4,0.5)(-1.4,0)
\qbezier(-1.4,1)(-1.4,1.4)(-1,1.4)
\qbezier(-0.6,1)(-0.6,1.4)(-1,1.4)
\qbezier(-1.4,0)(-1.4,-0.4)(-1,-0.4)
\qbezier(-1,-0.4)(-0.6,-0.4)(-0.6,0)
\qbezier(-1.2,1)(-0.6,1)(-0.4,1)
\qbezier(-1.2,0)(-0.6,0)(-0.4,0)
\qbezier(-1.6,1)(-1.8,1)(-1.8,1)
\qbezier(-1.6,0)(-1.8,0)(-1.8,0)
\qbezier(-0.4,0)(0.2,0)(0.2,-2)
\qbezier(-0.4,1)(1.2,1)(1.2,-2)
\qbezier(-1.8,0)(-2.4,0)(-2.4,-2)
\qbezier(-1.8,1)(-3.4,1)(-3.4,-2)
\qbezier(0.2,-2)(0.2,-3)(-0.4,-3)
\qbezier(-3.4,-2)(-3.4,-3)(-3.8,-3)
\qbezier(-2.4,-2)(-2.4,-3)(-2,-3)
\qbezier(1.2,-2)(1.2,-3)(1.6,-3)
\qbezier(-5,-2)(-5,-2)(-3.4,-2)
\qbezier(-2.4,-2)(0,-2)(0.2,-2)
\qbezier(1.2,-2)(2.8,-2)(2.8,-2)
\qbezier(-1,1.4)(-1,3)(-1,3)
\qbezier[20](-1.6,0.5)(-2.9,0.5)(-2.9,-2)
\qbezier[20](-0.4,0.5)(0.7,0.5)(0.7,-2)
\qbezier[6](-1.2,0.5)(-0.8,0.5)(-0.4,0.5)
\qbezier[20](0.7,-2)(0.7,-4)(-1,-4)
\qbezier[20](0.7,-2)(0.7,-4)(-1,-4)
\qbezier[20](-2.9,-2)(-2.9,-4)(-1,-4)
\qbezier[28](-1,-4)(-1,-6.25)(-1,-8.5)
}}

\def\modhookp{
\Picture{
\qbezier(-1.4,1)(-1.4,0.5)(-1.4,0)
\qbezier(-1.4,1)(-1.4,1.4)(-1,1.4)
\qbezier(-0.6,1)(-0.6,1.4)(-1,1.4)
\qbezier(-1.4,0)(-1.4,-0.4)(-1,-0.4)
\qbezier(-1,-0.4)(-0.6,-0.4)(-0.6,0)
\qbezier(-1.2,1)(-0.6,1)(-0.4,1)
\qbezier(-1.2,0)(-0.6,0)(-0.4,0)
\qbezier(-1.6,1)(-1.8,1)(-1.8,1)
\qbezier(-1.6,0)(-1.8,0)(-1.8,0)
\qbezier(-0.4,0)(0.2,0)(0.2,-2)
\qbezier(-0.4,1)(1.2,1)(1.2,-2)
\qbezier(-1.8,0)(-2.4,0)(-2.4,-2)
\qbezier(-1.8,1)(-3.4,1)(-3.4,-2)
\qbezier(0.2,-2)(0.2,-3)(-0.4,-3)
\qbezier(-3.4,-2)(-3.4,-3)(-3.8,-3)
\qbezier(-2.4,-2)(-2.4,-3)(-2,-3)
\qbezier(1.2,-2)(1.2,-3)(1.6,-3)
\qbezier(-9,-2)(-5,-2)(-3.4,-2)
\qbezier(-2.4,-2)(0,-2)(0.2,-2)
\qbezier(1.2,-2)(2.8,-2)(7,-2)
\qbezier(-1,1.4)(-1,3)(-1,3)
\qbezier[20](-1.6,0.5)(-2.9,0.5)(-2.9,-2)
\qbezier[20](-0.4,0.5)(0.7,0.5)(0.7,-2)
\qbezier[6](-1.2,0.5)(-0.8,0.5)(-0.4,0.5)
\qbezier[20](0.7,-2)(0.7,-4)(-1,-4)
\qbezier[20](0.7,-2)(0.7,-4)(-1,-4)
\qbezier[20](-2.9,-2)(-2.9,-4)(-1,-4)
\qbezier[28](-1,-4)(-1,-6)(-1,-6)
\qbezier[20](-1,-6)(-1,-8)(1,-8)
\qbezier[20](3,-6)(3,-8)(1,-8)
\qbezier[32](3,-6)(3,-4)(3,-2)
}}

\def\cross{
\Picture{
\qbezier(-1.4,1)(-1.4,0.5)(-1.4,0)
\qbezier(-1.4,1)(-1.4,1.4)(-1,1.4)
\qbezier(-0.6,1)(-0.6,1.4)(-1,1.4)
\qbezier(-1.4,0)(-1.4,-0.4)(-1,-0.4)
\qbezier(-1,-0.4)(-0.6,-0.4)(-0.6,0)
}}

\def\crossr{
\Picture{
\qbezier(--1.4,1)(--1.4,0.5)(--1.4,0)
\qbezier(--1.4,1)(--1.4,1.4)(--1,1.4)
\qbezier(--0.6,1)(--0.6,1.4)(--1,1.4)
\qbezier(--1.4,0)(--1.4,-0.4)(--1,-0.4)
\qbezier(--1,-0.4)(--0.6,-0.4)(--0.6,0)
}}

\def\seifcross{
\Picture{
\put(0,0){$\cross$}
\put(-1.6,-0.75){$\crossr$}
\put(-1,1.4){\line(0,1){3}}
\put(-0.6,-1.15){\line(0,-1){3}}
\put(-0.5,2.5){$B^{i+1}$}
\put(-0.2,-3){$C^{i}$}
}}

\def\seifcrossb{
\Picture{
\put(0,0){$\cross$}
\put(-1.6,-0.75){$\crossr$}
\put(-1,1.4){\line(0,1){3}}
\put(-0.6,-1.15){\line(0,-1){3}}
\put(-0.5,2.5){Y$_a$}
\put(-0.2,-3){Y$_b$}
}}

\def\ycomp{
\Picture{
\put(5,4.8){\circle{2}}
\put(-5,4.8){\circle{2}}
\put(5,4.8){\circle{1}}
\put(-5,4.8){\circle{1}}
\put(0,0){\circle{1.2}}
\qbezier(-0.3,-0.55)(-0.3,-1)(-0.3,-4)
\qbezier(0.3,-0.55)(0.3,-1)(0.3,-4)
\qbezier(0.58,-0.1)(2.6,1.9)(4.6,3.9)
\qbezier(-0.59,-0.1)(-2.6,1.9)(-4.6,3.9)
\qbezier(0.4,0.4)(2.4,2.4)(4.2,4.2)
\qbezier(-0.4,0.4)(-2.4,2.4)(-4.2,4.2)
\put(0,-5){\circle{2}}
\put(0,-5){\circle{1}}
\put(0,-5){\circle{1}}
}}

\def\actycompa{
\Picture{
\put(5,4.8){\circle{2}}
\put(-5,4.8){\circle{2}}
\put(5,4.8){\circle{1}}
\put(-6,4.6){\vector(0,-1){0.01}}
\put(-5.5,5){\vector(0,1){0.01}}
\put(-1,-5.2){\vector(0,-1){0.01}}
\put(-0.5,-4.8){\vector(0,1){0.01}}
\put(4,4.6){\vector(0,-1){0.01}}
\put(4.5,5){\vector(0,1){0.01}}
\put(-5,4.8){\circle{1}}
\put(0,0){\circle{1.2}}
\put(-0.3,-2){\vector(0,-1){0.01}}
\put(0.3,-2){\vector(0,1){0.01}}
\qbezier(-0.3,-0.55)(-0.3,-1)(-0.3,-4)
\qbezier(0.3,-0.55)(0.3,-1)(0.3,-4)
\qbezier(0.58,-0.1)(2.6,1.9)(4.6,3.9)
\put(2.6,1.9){\vector(1,1){0.01}}
\put(-2.6,1.9){\vector(1,-1){0.01}}
\qbezier(-0.59,-0.1)(-2.6,1.9)(-4.6,3.9)
\qbezier(0.4,0.4)(2.4,2.4)(4.2,4.2)
\put(2.4,2.4){\vector(-1,-1){0.01}}
\put(-2.4,2.4){\vector(-1,1){0.01}}
\qbezier(-0.4,0.4)(-2.4,2.4)(-4.2,4.2)
\put(0,-5){\circle{2}}
\put(0,-5){\circle{1}}
\put(0,-5){\circle{1}}
}}

\def\actycompb{
\Picture{
\put(5,4.8){\circle{2}}
\put(-5,4.8){\circle{2}}
\put(5,4.8){\circle{1}}
\put(-5,4.8){\circle{1}}
\put(0,0){\circle{1.2}}
\qbezier(-0.3,-0.55)(-0.3,-1)(-0.3,-4)
\qbezier(0.3,-0.55)(0.3,-1)(0.3,-4)
\qbezier(0.58,-0.1)(2.6,1.9)(4.6,3.9)
\qbezier(-0.59,-0.1)(-2.6,1.9)(-4.6,3.9)
\qbezier(0.4,0.4)(2.4,2.4)(4.2,4.2)
\qbezier(-0.4,0.4)(-2.4,2.4)(-4.2,4.2)
\put(0,-5){\circle{2}}
\put(0,-5){\circle{1}}
\put(0,-5){\circle{1}}
}}

\def\actycomp{
\Picture{
\put(5,4.8){\circle{2}}
\put(-5,4.8){\circle{2}}
\put(0,0){\circle*{1.2}}
\qbezier(0,-0.6)(0,-4)(0,-4)
\qbezier(0,0)(2.25,2)(4.5,4)
\qbezier(0,0)(-2.25,2)(-4.5,4)
\put(0,-5){\circle{2}}
}}

\def\chord{
\Picture{
\put(-1.5,2){$K$}
\put(4.5,2){$K$}
\qbezier(0,1)(0,2)(0,3)
\qbezier(0,-1)(0,-2)(0,-3)
\qbezier(4,1)(4,2)(4,3)
\qbezier(4,-1)(4,-2)(4,-3)
\put(0,0){$\clasp$}
\qbezier(1,0)(2,0)(3,0)
\put(4,0){$\claspr$}}
}

\def\strut{
\Picture{
\put(-3,0){$\vertclasp$}
\put(3,0){$\vertclasp$}
\qbezier[32](-6,0)(0,0)(6,0)
\put(4,-2){$K$}
\put(0,6){\circle*{1}}
\qbezier(-3,1)(-3,4)(0,6)
\qbezier(3,1)(3,4)(0,6)
\qbezier(0,6)(0,7)(0,9)
}}

\def\thetaA{
\Picture{
\put(-3,0){$\vertclasp$}
\qbezier(-3,1)(-3,2)(-3,3)
\put(3,0){$\vertclasp$}
\qbezier(3,1)(3,2)(3,3)
\qbezier[20](-3,0)(0,0)(3,0)
\put(-3,3){\circle*{1}}
\put(0,5){\circle*{1}}
\put(0,7){\circle*{1}}
\qbezier(0,5)(0,6)(0,7)
\put(3,3){\circle*{1}}
\qbezier[5](-2,1.5)(0,1.5)(2,1.5)
\qbezier(-5,3)(-3,3)(-2,3)
\qbezier(5,3)(3,3)(2,3)
\qbezier[16](-2,3)(0,3)(2,3)
\qbezier(-5,3)(-7,3)(-7,5)
\qbezier(5,3)(7,3)(7,5)
\put(-7,5){\circle*{1}}
\put(7,5){\circle*{1}}
\qbezier(-7,5)(0,5)(7,5)
\qbezier(-5,7)(-7,7)(-7,5)
\qbezier(5,7)(7,7)(7,5)
\qbezier(-5,7)(0,7)(5,7)
\qbezier(-4,0)(-6,0)(-6,-2)
\qbezier(4,0)(6,0)(6,-2)
\qbezier(-6,-2)(-6,-4)(0,-4)
\qbezier(6,-2)(6,-4)(0,-4)
\put(-3.5,-1){$\underbrace{\ \ \ \ \ \ \ \ \ \ \ \ \ \ \ \ \ \ \ }_n$}
\put(7.5,6){$\kappa_n$}
\put(6.75,-2){$U$}
}}

\def\thetab{
\Picture{
\put(0,5){\circle*{1}}
\put(0,7){\circle*{1}}
\qbezier(0,5)(0,6)(0,7)
\put(-3,5){\circle*{1}}
\put(3,5){\circle*{1}}
\qbezier(-3,5)(0,5)(3,5)
\qbezier(-3,5)(-3,7)(0,7)
\qbezier(3,5)(3,7)(0,7)
\qbezier(-3,5)(-3,3)(0,3)
\qbezier(3,5)(3,3)(0,3)
\put(3.8,5){$\lambda_{\kappa_n}$}
}}

\def\thetac{
\Picture{
\put(0,5){\circle*{0.3}}
\put(0,7){\circle*{0.3}}
\qbezier[12](0,5)(0,6)(0,7)
\put(-3,5){\circle*{0.3}}
\put(3,5){\circle*{0.3}}
\qbezier[20](-3,5)(0,5)(3,5)
\qbezier[12](-3,5)(-3,7)(0,7)
\qbezier[12](3,5)(3,7)(0,7)
\qbezier[18](-3,5)(-3,3)(0,3)
\qbezier[18](3,5)(3,3)(0,3)
\put(3.8,5){$\kappa$}
}
}

\def\boxdefn{
\Picture{
\put(2,-0.5){\vector(-1,0){4}}
\put(-2,0.5){\vector(1,0){4}}
\put(-1,-1){\line(0,1){2}}
\put(1,-1){\line(0,1){2}}
\put(-1,-1){\line(1,0){2}}
\put(-1,1){\line(1,0){2}}
\qbezier[20](-0.5,-1)(-0.5,-2.5)(-0.5,-4)
\qbezier[20](0.5,-1)(0.5,-2.5)(0.5,-4)
}}
\def\boxdefna{
\Picture{
\put(2,-0.5){\vector(-1,0){4}}
\put(-2,0.5){\vector(1,0){4}}
\qbezier[24](-0.5,-0.5)(-0.5,-2.25)(-0.5,-4)
\qbezier[24](0.5,-0.5)(0.5,-2.25)(0.5,-4)
}}
\def\boxdefnb{
\Picture{
\put(2,-0.5){\vector(-1,0){4}}
\put(-2,0.5){\vector(1,0){4}}
\qbezier[24](-0.5,-0.5)(-0.5,-2.25)(-0.5,-4)
\qbezier[24](0.5,0.5)(0.5,-1.75)(0.5,-4)
}}
\def\boxdefnc{
\Picture{
\put(2,-0.5){\vector(-1,0){4}}
\put(-2,0.5){\vector(1,0){4}}
\qbezier[24](-0.5,0.5)(-0.5,-1.75)(-0.5,-4)
\qbezier[24](0.5,-0.5)(0.5,-2.25)(0.5,-4)
}}
\def\boxdefnd{
\Picture{
\put(2,-0.5){\vector(-1,0){4}}
\put(-2,0.5){\vector(1,0){4}}
\qbezier[24](0.5,0.5)(0.5,-1.75)(0.5,-4)
\qbezier[24](-0.5,0.5)(-0.5,-1.75)(-0.5,-4)
}}
\def\shiftl{
\Picture{
\qbezier(-2,0.5)(-0,0.5)(2,0.5)
\qbezier(-2,-0.5)(-0,-0.5)(2,-0.5)
\qbezier(2,0.5)(3,0.5)(4,1.5)
\qbezier(4,1.5)(5,2.5)(6,2.5)
\qbezier(2,-0.5)(3,-0.5)(4,-1.5)
\qbezier(4,-1.5)(5,-2.5)(6,-2.5)
\qbezier(6,2.5)(7.5,2.5)(7.5,1)
\qbezier(6,-2.5)(7.5,-2.5)(7.5,-1)
\qbezier(7.5,1)(7.5,0)(7.5,-1)
\put(4,1.5){\vector(1,1){0.05}}
\qbezier[18](7.5,0.8)(9.5,0.8)(11.5,0.8)
\qbezier[18](7.5,-0.8)(9.5,-0.8)(11.5,-0.8)
\qbezier[3](11,0.6)(11,0)(11,-0.6)
\put(-1,-1){\line(1,0){2}}
\put(-1,1){\line(1,0){2}}
\put(-1,-1){\line(0,1){2}}
\put(1,-1){\line(0,1){2}}
\qbezier[24](-0.6,-1)(-0.6,-3)(-0.6,-5)
\qbezier[24](0.6,-1)(0.6,-3)(0.6,-5)
\qbezier[2](-0.4,-3)(0,-3)(0.4,-3)
\put(-0.3,-6.5){$n$}
\put(12,0){$m$}
}}
\def\shiftr{
\Picture{
\qbezier(-2,0.5)(-0,0.5)(2,0.5)
\qbezier(-2,-0.5)(-0,-0.5)(2,-0.5)
\qbezier(2,0.5)(3,0.5)(4,1.5)
\qbezier(4,1.5)(5,2.5)(6,2.5)
\qbezier(2,-0.5)(3,-0.5)(4,-1.5)
\qbezier(4,-1.5)(5,-2.5)(6,-2.5)
\qbezier(6,2.5)(7.5,2.5)(7.5,1)
\qbezier(6,-2.5)(7.5,-2.5)(7.5,-1)
\qbezier(7.5,1)(7.5,0)(7.5,-1)
\put(4,1.5){\vector(1,1){0.05}}
\qbezier[18](7.5,0.8)(9.5,0.8)(11.5,0.8)
\qbezier[18](7.5,-0.8)(9.5,-0.8)(11.5,-0.8)
\qbezier[3](11,0.6)(11,0)(11,-0.6)
\put(8.2,-1.2){\line(1,0){2}}
\put(8.2,1.2){\line(1,0){2}}
\put(8.2,-1.2){\line(0,1){2.4}}
\put(10.2,-1.2){\line(0,1){2.4}}
\qbezier[24](8.6,-1.2)(8.6,-3.2)(8.6,-5.2)
\qbezier[24](9.8,-1.2)(9.8,-3.2)(9.8,-5.2)
\qbezier[2](8.8,-3)(9.2,-3)(9.6,-3)
\put(9,-6.5){$n$}
\put(12,0){$m$}
}}

\def\glueexamp{
\Picture{
\put(-1,-0.75){\line(1,0){2}}
\put(-1,0.75){\line(1,0){2}}
\put(-1,-0.75){\line(0,1){1.5}}
\put(1,-0.75){\line(0,1){1.5}}
\put(-4,-0.75){\line(1,0){2}}
\put(-4,0.75){\line(1,0){2}}
\put(-4,-0.75){\line(0,1){1.5}}
\put(-2,-0.75){\line(0,1){1.5}}
\put(2,-0.75){\line(1,0){2}}
\put(2,0.75){\line(1,0){2}}
\put(2,-0.75){\line(0,1){1.5}}
\put(4,-0.75){\line(0,1){1.5}}
\put(-3,0.75){\vector(0,1){2}}
\put(0,0.75){\vector(0,1){2}}
\put(3,-0.75){\vector(0,-1){2}}
\put(-3,-0.75){\line(0,-1){2}}
\put(0,-0.75){\line(0,-1){2}}
\put(3,0.75){\line(0,1){2}}
\put(-3.3,-0.3){$t$}
\put(-0.3,-0.3){$t$}
\put(2.7,-0.3){$t$}
}}

\def\sweepa{
\Picture{
\put(-0.05,2.25){$\circ$}
\put(-2,4){\vector(0,-1){1}}
\put(2,3){\vector(0,1){1}}
\put(-2,2){\vector(0,-1){4}}
\put(2,-2){\vector(0,1){4}}
\put(-1.75,3){$t^{-1}$}
\put(2.25,3){$t$}
\qbezier[24](-2,0)(0,0)(2,0)
\put(-2.75,1){$f$}
\put(2.25,1){$g$}
\put(0,0){\vector(-1,0){0.01}}
\put(0,-1){$h$}
\put(-2.75,-1.5){$j$}
\put(2.25,-1.5){$k$}
}}

\def\sweepb{
\Picture{
\qbezier[22](-6,1)(0,1)(6,1)
\put(-2,2){\vector(0,-1){4}}
\put(2,-2){\vector(0,1){4}}
\qbezier[24](-2,-0.5)(0,0)(2,0.5)
\put(-4.5,1){$ft^{-1}$}
\put(2.25,1){$gt$}
\put(0,0){\vector(-4,-1){0.01}}
\put(0,-1){$h$}
\put(-2.75,-1.5){$j$}
\put(2.25,-1.5){$k$}
}}

\def\sweepc{
\Picture{
\qbezier[22](-6,0)(0,0)(6,0)
\put(-2,2){\vector(0,-1){4}}
\put(2,-2){\vector(0,1){4}}
\qbezier[24](-2,-0.5)(0,0)(2,0.5)
\put(-4.5,1){$ft^{-1}$}
\put(2.25,1){$g$}
\put(0,0){\vector(-4,-1){0.01}}
\put(0,-1){$ht^{-1}$}
\put(-2.75,-1.5){$j$}
\put(2.25,-1.5){$kt$}
}}

\def\sweepd{
\Picture{
\qbezier[22](-6,-1.25)(0,-1.25)(6,-1.25)
\put(-2,2){\vector(0,-1){4}}
\put(2,-2){\vector(0,1){4}}
\qbezier[24](-2,-0.5)(0,0)(2,0.5)
\put(-2.75,1){$f$}
\put(2.25,1){$g$}
\put(0,0){\vector(-4,-1){0.01}}
\put(0,-1){$h$}
\put(-4.5,-1.5){$jt^{-1}$}
\put(2.25,-1.5){$kt$}
}}

\def\sweepe{
\Picture{
\put(-0.05,-2.25){$\circ$}
\put(-2,-3){\vector(0,-1){1}}
\put(2,-4){\vector(0,1){1}}
\put(-2,2){\vector(0,-1){4}}
\put(2,-2){\vector(0,1){4}}
\put(-1.75,-3.25){$t^{-1}$}
\put(2.25,-3.25){$t$}
\qbezier[24](-2,0)(0,0)(2,0)
\put(-2.75,1){$f$}
\put(2.25,1){$g$}
\put(0,0){\vector(-1,0){0.01}}
\put(0,-1){$h$}
\put(-2.75,-1.5){$j$}
\put(2.25,-1.5){$k$}
}}

\def\figate{
\Picture{
\qbezier(6,6)(6,2)(5,2)
\qbezier(5,2)(4,2)(4,3)
\qbezier(5,2.5)(5,3.5)(4,3.5)
\qbezier(5,1.5)(5,0)(5,-1)
\qbezier(4,3.5)(3,3.5)(3,2)
\qbezier(3,2)(3,0)(3,-1)
\qbezier(4,4)(4,5)(2.5,5)
\qbezier(2.5,5)(1,5)(1,4)
\qbezier(1,4)(1,2)(1,-1)
}}

\def\figates{
\Picture{
\qbezier(-4,2)(0,2)(0,0)
\qbezier(0,0)(0,-1)(-0.1,-1)
\qbezier(-1,0)(-1,-2)(3,-2)
\qbezier(-1,0)(-1,1)(-0.9,1)
\qbezier(-0.5,1.4)(-0.25,2)(3,2)
\qbezier(-0.5,-1.5)(-0.75,-2)(-4,-2)
\put(3,2){$\twistunits$}
\put(5,2){$\twistunits$}
\put(7,2){$\twistunits$}
\put(3,-2){\line(1,0){6}}
}}

\def\figatesass{
\Picture{
\put(-16,0){$\figates$}
\put(-3,0){$\figates$}
\put(10,0){$\figates$}
\put(-11,-6){$t^{-1}\widetilde{K}$}
\put(2,-6){$\widetilde{K}$}
\put(15,-6){$t \widetilde{K}$}
\put(-20,3.5){\line(1,0){39}}
\put(-20,-3.5){\line(1,0){39}}
}}

\def\figateb{
\Picture{
\qbezier(1,-1)(1,-4)(2,-4)
\qbezier(3,-1)(3,-4)(2,-4)
\put(5,-1){$\twistunit$}
\put(5,-3){$\twistunit$}
\put(5,-5){$\twistunit$}
\qbezier(5,-7)(5,-7.5)(5,-8)
}}

\def\figatediag{
\Picture{
\put(5,6){\vector(0,1){0.01}}
\qbezier(5,-8)(5,-2)(4,-2)
\qbezier(4,-2)(3,-2)(3,-3.5)
\qbezier(3,-3.5)(3,-5)(2,-5)
\qbezier(2,-5)(1,-5)(1,-1)
\qbezier(5,6)(5,0)(4,0)
\qbezier(4,0)(3,0)(3,1.5)
\qbezier(3,1.5)(3,3)(2,3)
\qbezier(2,3)(1,3)(1,-1)
\qbezier[16](4,0)(4,-1)(4,-2)
\put(5.5,-4){$t$}
\put(5.5,2){$1$}
\put(-0.6,-2.5){$t^{-1}$}
\put(4.5,-1.3){$1$}
}
}

\def\figatediagb{
\Picture{
\put(5,6){\vector(0,1){0.01}}
\qbezier(5,-8)(5,-3)(5,-3)
\qbezier(5,-3)(5,-2)(4,-2)
\qbezier(4,-2)(3,-2)(3,-3.5)
\qbezier(3,-3.5)(3,-5)(2,-5)
\qbezier(2,-5)(1,-5)(1,-1)
\qbezier(5,6)(5,0)(4,0)
\qbezier(4,0)(3,0)(3,1.5)
\qbezier(3,1.5)(3,3)(2,3)
\qbezier(2,3)(1,3)(1,-1)
\qbezier[16](5,-4)(7,-4)(7,-5)
\qbezier[16](5,-6)(7,-6)(7,-5)
\put(7.3,-5.25){1}
\put(5.3,-5.25){1}
\put(5.25,-7.6){1}
\put(-1.25,-2.5){$t.t^{-1}$}
}
}

\def\figatediagc{
\Picture{
\put(0,-4){\vector(0,1){8}}
\qbezier[16](0,-1.5)(3,-1.5)(3,0)
\qbezier[16](0,1.5)(3,1.5)(3,0)
\put(0,1.5){\vector(-1,0){0.01}}
\put(3.5,1.5){$\frac{1}{2}({t-3+t^{-1}})$}
}}

\def\figateass{
\Picture{
\put(6,6){\vector(0,1){0.01}}
\put(0,0){$\figate$}
\put(0,-0.5){$\figateb$}
\qbezier[10](3,-1.25)(2,-1.25)(2,-0.75)
\qbezier[14](3,-1.25)(4,-1.25)(5,-1.25)
\qbezier[10](5,-1.25)(6,-1.25)(6,-0.75)
\qbezier[8](2,-0.75)(2,-0.25)(2.75,-0.25)
\qbezier[14](3.15,-0.25)(4,-0.25)(4.85,-0.25)
\qbezier[8](5.25,-0.25)(6,-0.25)(6,-0.75)
\qbezier(1,-1.5)(1,-0)(1,0)
}}

\def\figateassc{
\Picture{
\put(6,6){\vector(0,1){0.01}}
\put(0,0){$\figate$}
\put(0,-0.5){$\figateb$}
\qbezier[24](2,-1.25)(5,-1.25)(8,-1.25)
\qbezier(3,-1)(3,-1.25)(3,-1.5)
\qbezier(5,-1)(5,-1.25)(5,-1.5)
\qbezier(1,-1.5)(1,-0)(1,0)
}}

\def\figateassb{
\Picture{
\qbezier(-1,-8.5)(-1,-11)(2,-11)
\qbezier(2,-11)(5,-11)(5,-8.5)
\qbezier(6,6)(6,8)(2.5,8)
\qbezier(-1,6)(-1,8)(2.5,8)
\qbezier(-1,6)(-1,0)(-1,-8.5)
\put(-2.5,-4.5){$*$}
\put(0,0){$\figate$}
\put(0,-0.5){$\figateb$}
\qbezier(3,-1.25)(2,-1.25)(2,-0.75)
\qbezier(3,-1.25)(4,-1.25)(5,-1.25)
\qbezier(5,-1.25)(6,-1.25)(6,-0.75)
\qbezier(2,-0.75)(2,-0.25)(2.75,-0.25)
\qbezier(3.15,-0.25)(4,-0.25)(4.85,-0.25)
\qbezier(5.25,-0.25)(6,-0.25)(6,-0.75)
\qbezier(1,-1.5)(1,-0)(1,0)
}}

\def\figatep{
\Picture{
\qbezier(6,6)(6,2)(5,2)
\qbezier(5,2)(4,2)(4,3)
\qbezier(5,2.5)(5,3.5)(4,3.5)
\qbezier(5,1.5)(5,0)(5,-1)
\qbezier(4,3.5)(3,3.5)(3,2)
\qbezier(3,2)(3,0)(3,-1)
\qbezier(4,4)(4,5)(2.5,5)
\qbezier(2.5,5)(1,5)(1,4)
\qbezier(1,4)(1,2)(1,-1)
\put(6,7){\vector(0,1){1}}
\put(1,-3){\vector(0,1){1}}
\put(3,-2){\vector(0,-1){1}}
\put(5,-3){\vector(0,1){1}}
\put(2.1,-2.75){(}
\put(5.6,-2.75){)}
}}

\def\figatebp{
\Picture{
\put(2.1,0.25){(}
\put(5.6,0.25){)}
\put(1,0){\vector(0,1){1}}
\put(3,1){\vector(0,-1){1}}
\put(5,0){\vector(0,1){1}}
\put(5,-10){\vector(0,1){1}}
\qbezier(1,-1)(1,-4)(2,-4)
\qbezier(3,-1)(3,-4)(2,-4)
\put(5,-1){$\twistunit$}
\put(5,-3){$\twistunit$}
\put(5,-5){$\twistunit$}
\qbezier(5,-7)(5,-7.5)(5,-8)
}}

\def\twistunit{
\Picture{
\qbezier(0,0)(0,-1.5)(0.5,-1.5)
\qbezier(0.5,-1.5)(1,-1.5)(1,-1)
\qbezier(1,-1)(1,-0.5)(0.5,-0.5)
\qbezier(0.5,-0.5)(0.3,-0.5)(0.2,-1)
\qbezier(0,-2)(0,-1.5)(0.1,-1.3)
}}

\def\twistunits{
\Picture{
\qbezier(0,0)(1.5,0)(1.5,0.5)
\qbezier(1.5,0.5)(1.5,1)(1,1)
\qbezier(1,1)(0.5,1)(0.5,0.5)
\qbezier(0.5,0.5)(0.5,0.3)(1,0.2)
\qbezier(2,0)(1.5,0)(1.3,0.1)
}}

\def\qunit{
\Picture{
\put(-2,-0.75){\line(1,0){4}}
\put(-2,0.75){\line(1,0){4}}
\put(-2,-0.75){\line(0,1){1.5}}
\put(2,-0.75){\line(0,1){1.5}}
\put(-0.65,-0.25){$q(t)$}
\qbezier[8](-1.75,0.75)(-1.75,1.25)(-1.75,1.75)
\qbezier[8](1.75,0.75)(1.75,1.25)(1.75,1.75)
\qbezier[8](-1.75,-0.75)(-1.75,-1.25)(-1.75,-1.75)
\qbezier[8](1.75,-0.75)(1.75,-1.25)(1.75,-1.75)
\qbezier[4](-1.25,1.25)(0,1.25)(1.25,1.25)
\qbezier[4](-1.25,-1.25)(0,-1.25)(1.25,-1.25)
}
}

\def\qrev{
\Picture{
\put(-2,-0.75){\line(1,0){4}}
\put(-2,0.75){\line(1,0){4}}
\put(-2,-0.75){\line(0,1){1.5}}
\put(2,-0.75){\line(0,1){1.5}}
\put(-0.65,-0.25){$q(t)$}
\qbezier[8](-1.75,0.75)(-1.75,1.5)(-2,1.5)
\qbezier[8](-1.75,0.75)(-1.75,1.5)(-2,1.5)
\qbezier[8](-2,1.5)(-2.25,1.5)(-2.25,0.75)
\qbezier[24](-2.25,0.75)(-2.25,-1)(-2.25,-3)
\qbezier[24](-4,0.75)(-4,-1)(-4,-3)
\qbezier[4](-3.5,-2)(-3.15,-2)(-2.75,-2)
\qbezier[16](1.75,0.75)(1.75,3)(-0.5,3)
\qbezier[16](-4,0.75)(-4,3)(-0.5,3)
\qbezier[8](1.75,-0.75)(1.75,-1.5)(2,-1.5)
\qbezier[8](--1.75,-0.75)(--1.75,-1.5)(--2,-1.5)
\qbezier[8](--2,-1.5)(--2.25,-1.5)(--2.25,-0.75)
\qbezier[24](--2.25,-0.75)(--2.25,--1)(--2.25,--3)
\qbezier[24](--4,-0.75)(--4,--1)(--4,--3)
\qbezier[4](--3.5,--2)(3.15,--2)(--2.75,--2)
\qbezier[16](-1.75,-0.75)(-1.75,-3)(--0.5,-3)
\qbezier[16](--4,-0.75)(--4,-3)(--0.5,-3)
\qbezier[4](--1.25,-1.25)(-0,-1.25)(-1.25,-1.25)
\qbezier[4](--1.25,--1.25)(-0,--1.25)(-1.25,--1.25)
}
}

\def\throughunit{
\Picture{
\qbezier[4](-1.25,0)(0,0)(1.25,0)
\qbezier[20](-1.75,-1.75)(-1.75,0)(-1.75,1.75)
\qbezier[20](1.75,-1.75)(1.75,0)(1.75,1.75)
}
}

\def\runit{
\Picture{
\put(-2,-0.75){\line(1,0){4}}
\put(-2,0.75){\line(1,0){4}}
\put(-2,-0.75){\line(0,1){1.5}}
\put(2,-0.75){\line(0,1){1.5}}
\put(-0.65,-0.25){$r(t)$}
\qbezier[8](-1.75,0.75)(-1.75,1.25)(-1.75,1.75)
\qbezier[8](1.75,0.75)(1.75,1.25)(1.75,1.75)
\qbezier[8](-1.75,-0.75)(-1.75,-1.25)(-1.75,-1.75)
\qbezier[8](1.75,-0.75)(1.75,-1.25)(1.75,-1.75)
\qbezier[4](-1.25,1.25)(0,1.25)(1.25,1.25)
\qbezier[4](-1.25,-1.25)(0,-1.25)(1.25,-1.25)
}
}

\def\qcomm{
\Picture{
\put(-1.75,1.25){$\qunit$}
\put(1.75,-1.25){$\runit$}
\qbezier[20](3.5,0.5)(3.5,1.75)(3.5,3)
\qbezier[20](-3.5,-0.5)(-3.5,-1.75)(-3.5,-3)
\qbezier[4](0.5,2.5)(1.75,2.5)(3,2.5)
\qbezier[4](-0.5,-2.5)(-1.75,-2.5)(-3,-2.5)
}}

\def\qcommb{
\Picture{
\put(-1.75,-1.25){$\qunit$}
\put(1.75,1.25){$\runit$}
\qbezier[20](3.5,-0.5)(3.5,-1.75)(3.5,-3)
\qbezier[20](-3.5,0.5)(-3.5,1.75)(-3.5,3)
\qbezier[4](0.5,-2.5)(1.75,-2.5)(3,-2.5)
\qbezier[4](-0.5,2.5)(-1.75,2.5)(-3,2.5)
}}

\def\qrunit{
\Picture{
\put(-2,-0.75){\line(1,0){4}}
\put(-2,0.75){\line(1,0){4}}
\put(-2,-0.75){\line(0,1){1.5}}
\put(2,-0.75){\line(0,1){1.5}}
\put(-1.5,-0.25){$q(t)r(t)$}
\qbezier[18](-1.75,0.75)(-1.75,2)(-1.75,3)
\qbezier[18](1.75,0.75)(1.75,2)(1.75,3)
\qbezier[18](-1.75,-0.75)(-1.75,-2)(-1.75,-3)
\qbezier[18](1.75,-0.75)(1.75,-2)(1.75,-3)
\qbezier[4](-1.25,2)(0,2)(1.25,2)
\qbezier[4](-1.25,-2)(0,-2)(1.25,-2)
}
}

\def\tunit{
\Picture{
\put(-2,-0.75){\line(1,0){4}}
\put(-2,0.75){\line(1,0){4}}
\put(-2,-0.75){\line(0,1){1.5}}
\put(2,-0.75){\line(0,1){1.5}}
\put(-0.25,-0.25){$t$}
\qbezier[18](-1.75,0.75)(-1.75,2)(-1.75,3)
\qbezier[18](1.75,0.75)(1.75,2)(1.75,3)
\qbezier[18](-1.75,-0.75)(-1.75,-2)(-1.75,-3)
\qbezier[18](1.75,-0.75)(1.75,-2)(1.75,-3)
\qbezier[4](-1.25,2)(0,2)(1.25,2)
\qbezier[4](-1.25,-2)(0,-2)(1.25,-2)
}
}

\def\funit{
\Picture{
\put(-2,-0.75){\line(1,0){4}}
\put(-2,0.75){\line(1,0){4}}
\put(-2,-0.75){\line(0,1){1.5}}
\put(2,-0.75){\line(0,1){1.5}}
\put(-1,-0.25){$f(k)$}
\qbezier[18](-1.75,0.75)(-1.75,2)(-1.75,3)
\qbezier[18](1.75,0.75)(1.75,2)(1.75,3)
\qbezier[18](-1.75,-0.75)(-1.75,-2)(-1.75,-3)
\qbezier[18](1.75,-0.75)(1.75,-2)(1.75,-3)
\qbezier[4](-1.25,2)(0,2)(1.25,2)
\qbezier[4](-1.25,-2)(0,-2)(1.25,-2)
}
}

\def\addunit{
\Picture{
\put(-3,-0.75){\line(1,0){6}}
\put(-3,0.75){\line(1,0){6}}
\put(-3,-0.75){\line(0,1){1.5}}
\put(3,-0.75){\line(0,1){1.5}}
\put(-2.5,-0.25){$aq(t) + br(t)$}
\qbezier[18](-1.75,0.75)(-1.75,2)(-1.75,3)
\qbezier[18](1.75,0.75)(1.75,2)(1.75,3)
\qbezier[18](-1.75,-0.75)(-1.75,-2)(-1.75,-3)
\qbezier[18](1.75,-0.75)(1.75,-2)(1.75,-3)
\qbezier[4](-1.25,2)(0,2)(1.25,2)
\qbezier[4](-1.25,-2)(0,-2)(1.25,-2)
}
}

\def\invunit{
\Picture{
\put(-2,-0.75){\line(1,0){4}}
\put(-2,0.75){\line(1,0){4}}
\put(-2,-0.75){\line(0,1){1.5}}
\put(2,-0.75){\line(0,1){1.5}}
\put(-1.2,-0.25){$q(t^{-1})$}
\qbezier[18](-1.75,0.75)(-1.75,2)(-1.75,3)
\qbezier[18](1.75,0.75)(1.75,2)(1.75,3)
\qbezier[18](-1.75,-0.75)(-1.75,-2)(-1.75,-3)
\qbezier[18](1.75,-0.75)(1.75,-2)(1.75,-3)
\qbezier[4](-1.25,2)(0,2)(1.25,2)
\qbezier[4](-1.25,-2)(0,-2)(1.25,-2)
}
}

\def\qlong{
\Picture{
\put(-2,-0.75){\line(1,0){4}}
\put(-2,0.75){\line(1,0){4}}
\put(-2,-0.75){\line(0,1){1.5}}
\put(2,-0.75){\line(0,1){1.5}}
\put(-0.6,-0.25){$q(t)$}
\qbezier[18](-1.75,0.75)(-1.75,2)(-1.75,3)
\qbezier[18](1.75,0.75)(1.75,2)(1.75,3)
\qbezier[18](-1.75,-0.75)(-1.75,-2)(-1.75,-3)
\qbezier[18](1.75,-0.75)(1.75,-2)(1.75,-3)
\qbezier[4](-1.25,2)(0,2)(1.25,2)
\qbezier[4](-1.25,-2)(0,-2)(1.25,-2)
}
}

\def\rlong{
\Picture{
\put(-2,-0.75){\line(1,0){4}}
\put(-2,0.75){\line(1,0){4}}
\put(-2,-0.75){\line(0,1){1.5}}
\put(2,-0.75){\line(0,1){1.5}}
\put(-0.6,-0.25){$r(t)$}
\qbezier[18](-1.75,0.75)(-1.75,2)(-1.75,3)
\qbezier[18](1.75,0.75)(1.75,2)(1.75,3)
\qbezier[18](-1.75,-0.75)(-1.75,-2)(-1.75,-3)
\qbezier[18](1.75,-0.75)(1.75,-2)(1.75,-3)
\qbezier[4](-1.25,2)(0,2)(1.25,2)
\qbezier[4](-1.25,-2)(0,-2)(1.25,-2)
}
}

\def\manyq{
\Picture{
\put(-2,-0.25){$t$}
\put(1.5,-0.25){$t$}
\put(-2.75,-0.75){\line(1,0){2}}
\put(-2.75,0.75){\line(1,0){2}}
\put(-2.75,-0.75){\line(0,1){1.5}}
\put(-0.75,-0.75){\line(0,1){1.5}}
\put(2.75,-0.75){\line(-1,0){2}}
\put(2.75,0.75){\line(-1,0){2}}
\put(2.75,-0.75){\line(0,1){1.5}}
\put(0.75,-0.75){\line(0,1){1.5}}
\qbezier[18](-1.75,0.75)(-1.75,2)(-1.75,3)
\qbezier[18](1.75,0.75)(1.75,2)(1.75,3)
\qbezier[18](-1.75,-0.75)(-1.75,-2)(-1.75,-3)
\qbezier[18](1.75,-0.75)(1.75,-2)(1.75,-3)
\qbezier[4](-1.25,2)(0,2)(1.25,2)
\qbezier[4](-1.25,-2)(0,-2)(1.25,-2)
}
}

\def\qforka{
\Picture{
\put(-2,-0.75){\line(1,0){4}}
\put(-2,0.75){\line(1,0){4}}
\put(-2,-0.75){\line(0,1){1.5}}
\put(2,-0.75){\line(0,1){1.5}}
\put(-0.65,-0.25){$q(t)$}
\qbezier[18](-1.75,0.75)(-1.75,2)(-1.75,3)
\qbezier[18](1.75,0.75)(1.75,2)(1.75,3)
\qbezier[18](-1.75,-0.75)(-1.75,-2)(-1.75,-3)
\qbezier[18](1.75,-0.75)(1.75,-2)(1.75,-3)
\qbezier[18](0,0.75)(0,2)(0,3)
\qbezier[4](0,-0.75)(0,-1)(0,-1.25)
\put(0,-1.25){\circle{0.01}}
\qbezier[18](0,-1.25)(-0.5,-1.25)(-0.5,-3)
\qbezier[18](0,-1.25)(0.5,-1.25)(0.5,-3)
\qbezier[4](-1.25,2)(0,2)(1.25,2)
\qbezier[2](-1.5,-2)(-1.1,-2)(-0.75,-2)
\qbezier[2](1.5,-2)(1.1,-2)(0.75,-2)
}
}

\def\qforkb{
\Picture{
\put(-2,-0.75){\line(1,0){4}}
\put(-2,0.75){\line(1,0){4}}
\put(-2,-0.75){\line(0,1){1.5}}
\put(2,-0.75){\line(0,1){1.5}}
\put(-0.65,-0.25){$q(t)$}
\qbezier[18](-1.75,0.75)(-1.75,2)(-1.75,3)
\qbezier[18](1.75,0.75)(1.75,2)(1.75,3)
\qbezier[18](-1.75,-0.75)(-1.75,-2)(-1.75,-3)
\qbezier[18](1.75,-0.75)(1.75,-2)(1.75,-3)
\qbezier[18](0,1.75)(0,2.25)(0,3)
\put(0,1.75){\circle{0.01}}
\qbezier[8](0,1.75)(0.5,1.75)(0.5,0.75)
\qbezier[8](0,1.75)(-0.5,1.75)(-0.5,0.75)
\qbezier[18](-0.5,-0.75)(-0.5,-1.9)(-0.5,-3)
\qbezier[18](0.5,-0.75)(0.5,-1.9)(0.5,-3)
\qbezier[4](-1.25,2.5)(0,2.5)(1.25,2.5)
\qbezier[2](-1.5,-2)(-1.1,-2)(-0.75,-2)
\qbezier[2](1.5,-2)(1.1,-2)(0.75,-2)
}
}

\def\qloopa{
\Picture{
\put(-2,-0.75){\line(1,0){4}}
\put(-2,0.75){\line(1,0){4}}
\put(-2,-0.75){\line(0,1){1.5}}
\put(2,-0.75){\line(0,1){1.5}}
\put(-0.65,-0.25){$q(t)$}
\qbezier[18](-1.75,0.75)(-1.75,2)(-1.75,3)
\qbezier[18](1.75,0.75)(1.75,2)(1.75,3)
\qbezier[18](-1.75,-0.75)(-1.75,-2)(-1.75,-3)
\qbezier[18](1.75,-0.75)(1.75,-2)(1.75,-3)
\qbezier[18](0,-1.25)(-0.5,-1.25)(-0.5,-3)
\qbezier[18](0,-1.25)(0.5,-1.25)(0.5,-3)
\qbezier[2](-1.5,-2)(-1.1,-2)(-0.75,-2)
\qbezier[2](1.5,-2)(1.1,-2)(0.75,-2)
\qbezier[4](-1.25,2)(0,2)(1.25,2)
}
}

\def\qloopb{
\Picture{
\qbezier[4](-1.25,2.5)(0,2.5)(1.25,2.5)
\put(-2,-0.75){\line(1,0){4}}
\put(-2,0.75){\line(1,0){4}}
\put(-2,-0.75){\line(0,1){1.5}}
\put(2,-0.75){\line(0,1){1.5}}
\put(-0.65,-0.25){$q(t)$}
\qbezier[18](-1.75,0.75)(-1.75,2)(-1.75,3)
\qbezier[18](1.75,0.75)(1.75,2)(1.75,3)
\qbezier[18](-1.75,-0.75)(-1.75,-2)(-1.75,-3)
\qbezier[18](1.75,-0.75)(1.75,-2)(1.75,-3)
\qbezier[8](0,1.75)(0.5,1.75)(0.5,0.75)
\qbezier[8](0,1.75)(-0.5,1.75)(-0.5,0.75)
\qbezier[18](-0.5,-0.75)(-0.5,-1.9)(-0.5,-3)
\qbezier[18](0.5,-0.75)(0.5,-1.9)(0.5,-3)
\qbezier[2](-1.5,-2)(-1.1,-2)(-0.75,-2)
\qbezier[2](1.5,-2)(1.1,-2)(0.75,-2)
}
}

\def\qtimesr{
\Picture{
\put(0,1.25){$
\Picture{
\put(-2,-0.75){\line(1,0){4}}
\put(-2,0.75){\line(1,0){4}}
\put(-2,-0.75){\line(0,1){1.5}}
\put(2,-0.75){\line(0,1){1.5}}
\put(-0.65,-0.25){$q(t)$}
\qbezier[8](-1.75,0.75)(-1.75,1.25)(-1.75,1.75)
\qbezier[8](1.75,0.75)(1.75,1.25)(1.75,1.75)
\qbezier[8](-1.75,-0.75)(-1.75,-1.25)(-1.75,-1.75)
\qbezier[8](1.75,-0.75)(1.75,-1.25)(1.75,-1.75)
\qbezier[4](-1.25,1.25)(0,1.25)(1.25,1.25)
\qbezier[4](-1.25,-1.25)(0,-1.25)(1.25,-1.25)
}$}
\put(0,-1.25){$
\Picture{
\put(-2,-0.75){\line(1,0){4}}
\put(-2,0.75){\line(1,0){4}}
\put(-2,-0.75){\line(0,1){1.5}}
\put(2,-0.75){\line(0,1){1.5}}
\put(-0.65,-0.25){$r(t)$}
\qbezier[8](-1.75,0.75)(-1.75,1.25)(-1.75,1.75)
\qbezier[8](1.75,0.75)(1.75,1.25)(1.75,1.75)
\qbezier[8](-1.75,-0.75)(-1.75,-1.25)(-1.75,-1.75)
\qbezier[8](1.75,-0.75)(1.75,-1.25)(1.75,-1.75)
\qbezier[4](-1.25,1.25)(0,1.25)(1.25,1.25)
\qbezier[4](-1.25,-1.25)(0,-1.25)(1.25,-1.25)
}$}}
}

\def\coupon{
\Picture{
\put(-1,-0.75){\line(1,0){2}}
\put(-1,0.75){\line(1,0){2}}
\put(-1,-0.75){\line(0,1){1.5}}
\put(1,-0.75){\line(0,1){1.5}}
\put(-0.1,-0.2){$t$}
\qbezier[12](0,0.75)(0,1.75)(0,2.75)
\qbezier[12](0,-0.75)(0,-1.75)(0,-2.75)
}}

\def\fcoupon{
\Picture{
\put(-1,-0.75){\line(1,0){2}}
\put(-1,0.75){\line(1,0){2}}
\put(-1,-0.75){\line(0,1){1.5}}
\put(1,-0.75){\line(0,1){1.5}}
\put(-0.95,-0.2){$f(k)$}
\qbezier[12](0,0.75)(0,1.75)(0,2.75)
\qbezier[12](0,-0.75)(0,-1.75)(0,-2.75)
}}

\def\fcoupa{
\Picture{
\qbezier[24](0,-2.75)(0,0)(0,2.75)
}}

\def\fcoupb{
\Picture{
\qbezier[24](0,-2.75)(0,0)(0,2.75)
\qbezier[12](0,0)(1,0)(2,0)
\put(3,-0.25){$k$}
}}

\def\fskx{
\Picture{
\qbezier(0,-2.75)(0,0)(0,2.75)
\put(0,2.75){\vector(0,1){0.01}}
\put(-0.2,-0.2){x}
\qbezier[12](0,0)(1,0)(2,0)
\put(3,-0.25){$k$}
}}

\def\fsk{
\Picture{
\qbezier(0,-2.75)(0,0)(0,2.75)
\put(0,2.75){\vector(0,1){0.01}}
\qbezier[12](0,0)(1,0)(2,0)
\put(3,-0.25){$k$}
}}

\def\fskxr{
\Picture{
\qbezier(0,-2.75)(0,0)(0,2.75)
\put(0,-2.75){\vector(0,-1){0.01}}
\put(-0.2,-0.2){x}
\qbezier[12](0,0)(1,0)(2,0)
\put(3,-0.25){$k$}
}}

\def\fskr{
\Picture{
\qbezier(0,-2.75)(0,0)(0,2.75)
\put(0,-2.75){\vector(0,-1){0.01}}
\qbezier[12](0,0)(1,0)(2,0)
\put(3,-0.25){$k$}
}}

\def\fcoupbx{
\Picture{
\qbezier[24](0,-2.75)(0,0)(0,2.75)
\qbezier[12](0,0)(1,0)(2,0)
\put(-0.2,-0.2){x}
\put(3,-0.25){$k$}
}}

\def\fcoupc{
\Picture{
\qbezier[24](0,-2.75)(0,0)(0,2.75)
\qbezier[12](0,1)(1,1)(2,1)
\qbezier[12](0,-1)(1,-1)(2,-1)
\put(3,0.75){$k$}
\put(3,-1.25){$k$}
}}

\def\fcoupcx{
\Picture{
\qbezier[24](0,-2.75)(0,0)(0,2.75)
\qbezier[12](0,1)(1,1)(2,1)
\put(-0.2,0.8){x}
\put(-0.2,-1.2){x}
\qbezier[12](0,-1)(1,-1)(2,-1)
\put(3,0.75){$k$}
\put(3,-1.25){$k$}
}}

\def\fcoupd{
\Picture{
\qbezier[24](0,-2.75)(0,0)(0,2.75)
\qbezier[12](0,1.5)(1,1.5)(2,1.5)
\qbezier[12](0,-1.5)(1,-1.5)(2,-1.5)
\qbezier[12](0,0)(1,0)(2,0)
\put(3,1.25){$k$}
\put(3,-0.25){$k$}
\put(3,-1.75){$k$}
}}

\def\fcoupdx{
\Picture{
\qbezier[24](0,-2.75)(0,0)(0,2.75)
\put(-0.2,1.3){x}
\put(-0.2,-0.2){x}
\put(-0.2,-1.7){x}
\qbezier[12](0,1.5)(1,1.5)(2,1.5)
\qbezier[12](0,-1.5)(1,-1.5)(2,-1.5)
\qbezier[12](0,0)(1,0)(2,0)
\put(3,1.25){$k$}
\put(3,-0.25){$k$}
\put(3,-1.75){$k$}
}}

\def\coupup{
\Picture{
\put(-1,-0.75){\line(1,0){2}}
\put(-1,0.75){\line(1,0){2}}
\put(-1,-0.75){\line(0,1){1.5}}
\put(1,-0.75){\line(0,1){1.5}}
\put(-0.6,-0.2){$t^{-1}$}
\qbezier[12](0,0.75)(0,2.5)(-1,2.5)
\qbezier[32](-1,2.5)(-2,2.5)(-2,-2.75)
\qbezier[12](0,-0.75)(0,-2.5)(1,-2.5)
\qbezier[32](1,-2.5)(2,-2.5)(2,2.75)
}}

\def\couponinv{
\Picture{
\put(-1,-0.75){\line(1,0){2}}
\put(-1,0.75){\line(1,0){2}}
\put(-1,-0.75){\line(0,1){1.5}}
\put(1,-0.75){\line(0,1){1.5}}
\put(-0.6,-0.2){$t^{-1}$}
\qbezier[12](0,0.75)(0,1.75)(0,2.75)
\qbezier[12](0,-0.75)(0,-1.75)(0,-2.75)
}}

\def\unit{
\Picture{
\qbezier[24](0,4.5)(0,0)(0,-4.5)
}}

\def\couponmult{
\Picture{
\put(0,1.75){$\coupon$}
\put(0,-1.75){$\couponinv$}
}}

\def\pusha{
\Picture{
\put(0,1.75){$\coupon$}
\put(0,-1){\circle*{0.02}}
\qbezier[16](0,-1)(1,-2)(1,-4)
\qbezier[16](0,-1)(-1,-2)(-1,-4)}}

\def\pushb{
\Picture{
\qbezier[24](0,4.5)(0,1)(0,1)
\put(0,1){\circle*{0.02}}
\qbezier[12](0,1)(-1.5,0)(-1.5,-1)
\qbezier[12](0,1)(1.5,0)(1.5,-1)
\put(-2.5,-1){\line(1,0){2}}
\put(-2.5,-2.5){\line(1,0){2}}
\put(-2.5,-2.5){\line(0,1){1.5}}
\put(-0.5,-2.5){\line(0,1){1.5}}
\put(2.5,-1){\line(-1,0){2}}
\put(2.5,-2.5){\line(-1,0){2}}
\put(2.5,-2.5){\line(0,1){1.5}}
\put(0.5,-2.5){\line(0,1){1.5}}
\put(-1.6,-1.9){$t$}
\put(1.4,-1.9){$t$}
\qbezier[8](-1.5,-2.5)(-1.5,-3.25)(-1.5,-4)
\qbezier[8](1.5,-2.5)(1.5,-3.25)(1.5,-4)
}}

\def\wcoup{
\ \ \ \ \ \ 
\Picture{
\put(0,0.3){$\Picture{
\put(-0.5,-3.5){$x_i$}
\put(-0.5,3.2){$x_j$}
\put(-1.5,-0.75){\line(1,0){3}}
\put(-1.5,0.75){\line(1,0){3}}
\put(-1.5,-0.75){\line(0,1){1.5}}
\put(1.5,-0.75){\line(0,1){1.5}}
\put(-1.4,-0.2){$W_{ij}(t)$}
\qbezier[12](0,0.75)(0,1.75)(0,2.75)
\qbezier[12](0,-0.75)(0,-1.75)(0,-2.75)}$}
}\ 
}

\def\wcoupk{
\ \ \ \ \ \ 
\Picture{
\put(0,0.3){$\Picture{
\put(-0.5,-3.5){$x_i$}
\put(-0.5,3.2){$x_j$}
\put(-1.5,-0.75){\line(1,0){3}}
\put(-1.5,0.75){\line(1,0){3}}
\put(-1.5,-0.75){\line(0,1){1.5}}
\put(1.5,-0.75){\line(0,1){1.5}}
\put(-1.4,-0.2){$W_{ij}(e^k)$}
\qbezier[12](0,0.75)(0,1.75)(0,2.75)
\qbezier[12](0,-0.75)(0,-1.75)(0,-2.75)}$}
}\ 
}

\def\wcoupii{
\ \ \ \ \ \ 
\Picture{
\put(0,0.3){$\Picture{
\put(-0.5,-3.5){$x_i$}
\put(-0.5,3.2){$x_i'$}
\qbezier[32](0,2.75)(0,0)(0,-2.75)}$}
}\ 
}

\def\wcoupinv{
\ \ \ \ \ \ 
\Picture{
\put(0,0.3){$\Picture{
\put(-0.5,-3.5){$x_i$}
\put(-0.5,3.2){$x_j$}
\put(-2,-0.75){\line(1,0){4}}
\put(-2,0.75){\line(1,0){4}}
\put(-2,-0.75){\line(0,1){1.5}}
\put(2,-0.75){\line(0,1){1.5}}
\put(-1.6,-0.2){$W_{ij}^{-1}(t)$}
\qbezier[12](0,0.75)(0,1.75)(0,2.75)
\qbezier[12](0,-0.75)(0,-1.75)(0,-2.75)}$}
}\ 
}

\def\wcoupinvk{
\ \ \ \ \ \ 
\Picture{
\put(0,0.3){$\Picture{
\put(-0.5,-3.5){$x_i$}
\put(-0.5,3.2){$x_j$}
\put(-2,-0.75){\line(1,0){4}}
\put(-2,0.75){\line(1,0){4}}
\put(-2,-0.75){\line(0,1){1.5}}
\put(2,-0.75){\line(0,1){1.5}}
\put(-1.9,-0.2){$W_{ij}^{-1}(e^k)$}
\qbezier[12](0,0.75)(0,1.75)(0,2.75)
\qbezier[12](0,-0.75)(0,-1.75)(0,-2.75)}$}
}\ 
}

\def\wcoupinvkp{
\ \ \ \ \ \ 
\Picture{
\put(0,0.3){$\Picture{
\put(-0.5,-3.5){$x_i'$}
\put(-0.5,3.2){$x_j'$}
\put(-2,-0.75){\line(1,0){4}}
\put(-2,0.75){\line(1,0){4}}
\put(-2,-0.75){\line(0,1){1.5}}
\put(2,-0.75){\line(0,1){1.5}}
\put(-1.9,-0.2){$W_{ij}^{-1}(e^k)$}
\qbezier[12](0,0.75)(0,1.75)(0,2.75)
\qbezier[12](0,-0.75)(0,-1.75)(0,-2.75)}$}
}\ 
}

\def\wcoupk{
\ \ \ \ \ \ 
\Picture{
\put(0,0.3){$\Picture{
\put(-0.5,-3.5){$x_i$}
\put(-0.5,3.2){$x_j$}
\put(-2,-0.75){\line(1,0){4}}
\put(-2,0.75){\line(1,0){4}}
\put(-2,-0.75){\line(0,1){1.5}}
\put(2,-0.75){\line(0,1){1.5}}
\put(-1.75,-0.2){$W_{ij}(e^k)$}
\qbezier[12](0,0.75)(0,1.75)(0,2.75)
\qbezier[12](0,-0.75)(0,-1.75)(0,-2.75)}$}
}\ 
}

\def\tbit{
\Picture{
\put(-0.75,-0.75){\line(1,0){1.5}}
\put(0.75,-0.75){\line(0,1){1.5}}
\put(-0.75,0.75){\line(1,0){1.5}}
\put(-0.75,-0.75){\line(0,1){1.5}}
\put(-0.2,-0.75){$^t$}
}}

\def\tinv{
\Picture{
\put(-0.75,-0.75){\line(1,0){1.5}}
\put(0.75,-0.75){\line(0,1){1.5}}
\put(-0.75,0.75){\line(1,0){1.5}}
\put(-0.75,-0.75){\line(0,1){1.5}}
\put(-0.5,-0.8){$^\frac{1}{t}$}
}}

\def\windexamp{
\Picture{
\put(0,-10){\vector(0,1){20}}
\put(7.5,0){$\tinv$}
\qbezier[16](5,0)(5,2.6)(6.25,3)
\qbezier[16](7.5,0.75)(7.5,2.6)(6.25,3)
\qbezier[16](5,0)(5,-2.6)(6.25,-3)
\qbezier[16](7.5,-0.75)(7.5,-2.6)(6.25,-3)
\put(12,2.5){$\tbit$}
\put(12,-2.5){$\tbit$}
\qbezier[16](12,-1.75)(12,0)(12,1.75)
\qbezier[24](12,3.25)(12,5)(9,5)
\qbezier[24](6.25,3)(6.25,5)(9,5)
\qbezier[24](12,-3.25)(12,-5)(9,-5)
\qbezier[24](6.25,-3)(6.25,-5)(9,-5)
\qbezier[36](9,5)(9,7.5)(0,7.5)
\qbezier[36](9,-5)(9,-7.5)(0,-7.5)
}}

\def\windexampb{
\Picture{
\put(0,-10){\vector(0,1){20}}
\qbezier(7.5,0)(10,-2)(12.5,0)
\qbezier(7.5,0)(10,2)(12.5,0)
\qbezier(6.25,1)(6.25,1)(7.5,0)
\qbezier(13.75,1)(13.75,1)(12.5,0)
\qbezier(-2,-10)(10,-11.5)(22,-10)
\qbezier(-2,-10)(10,-8.5)(22,-10)
\qbezier(-2,10)(10,11.5)(22,10)
\qbezier[80](-2,10)(10,8.5)(22,10)
\qbezier(22,10)(22,0)(22,-10)
\qbezier(-2,-10)(-2,0)(-2,10)
\qbezier[42](7.5,7.5)(20.5,7.5)(20.5,0)
\qbezier[42](10,-7.5)(20.5,-7.5)(20.5,0)
\qbezier[32](10,-7.5)(5,-7.5)(5,0)
\qbezier[32](5,0)(5,5)(12,5)
\qbezier[32](12,5)(18,5)(18,0)
\qbezier[50](18,0)(18,-7)(2,-7)
\qbezier[24](4,6)(6,6)(6,0)
\qbezier[24](6,0)(6,-5)(10,-5)
\qbezier[24](10,-5)(16,-5)(16,0)
\qbezier[24](16,0)(16,3)(10,3)
\qbezier[24](3,4)(4,3)(10,3)
\qbezier[12](3,4)(3,6)(4,6)
\qbezier[6](4,6)(4,6.75)(4,7.5)
\qbezier[16](7.5,7.5)(3.7,7.5)(0,7.5)
\qbezier[32](3,4)(2,3)(2,-7)
\qbezier[6](2,-7)(1,-7)(0,-7)
}}

\def\pcoupon{
\Picture{
\put(-1,-0.75){\line(1,0){2}}
\put(-1,0.75){\line(1,0){2}}
\put(-1,-0.75){\line(0,1){1.5}}
\put(1,-0.75){\line(0,1){1.5}}
\put(-0.75,-0.2){$p(t)$}
\qbezier[12](0,0.75)(0,3)(0,5)
\qbezier[12](0,-0.75)(0,-3)(0,-5)
}}

\def\pcouponb{
\Picture{
\qbezier[32](0,5)(0,0)(0,-5)
}}

\def\pcouponc{
\Picture{
\put(-1,-0.75){\line(1,0){2}}
\put(-1,0.75){\line(1,0){2}}
\put(-1,-0.75){\line(0,1){1.5}}
\put(1,-0.75){\line(0,1){1.5}}
\put(-0.1,-0.2){$t$}
\qbezier[12](0,0.75)(0,3)(0,5)
\qbezier[12](0,-0.75)(0,-3)(0,-5)
}}

\def\pcouponn{
\Picture{
\put(-1,2.25){\line(1,0){2}}
\put(-1,3.75){\line(1,0){2}}
\put(-1,2.25){\line(0,1){1.5}}
\put(1,2.25){\line(0,1){1.5}}
\put(-0.1,2.8){$t$}
\put(-0.1,-3.2){$t$}
\put(-1,-2.25){\line(1,0){2}}
\put(-1,-3.75){\line(1,0){2}}
\put(-1,-2.25){\line(0,-1){1.5}}
\put(1,-2.25){\line(0,-1){1.5}}
\qbezier[8](0,3.75)(0,4.4)(0,5)
\qbezier[8](0,-3.75)(0,-4.4)(0,-5)
\qbezier[4](0,2.25)(0,2)(0,1.75)
\qbezier[4](0,-2.25)(0,-2)(0,-1.75)
\qbezier[4](0,1.75)(0,0)(0,-1.75)
\put(1,-0.2){
$
\left.
\begin{array}{l}
\\ \\ \\ \\ \\ \\ 
\end{array} \right\}
$}
\put(3.75,-0.2){$n$}
}}

\def\pcoup{
\Picture{
\put(-1,-0.75){\line(1,0){2}}
\put(-1,0.75){\line(1,0){2}}
\put(-1,-0.75){\line(0,1){1.5}}
\put(1,-0.75){\line(0,1){1.5}}
\put(-0.75,-0.2){$p(t)$}
\qbezier[12](0,0.75)(0,1.75)(0,2.75)
\qbezier[12](0,-0.75)(0,-1.75)(0,-2.75)
}}

\def\nopcoup{
\Picture{
\qbezier[32](0,-2.75)(0,0)(0,2.75)
\put(0,1.5){\vector(0,1){0.5}}
\put(0.6,-0.2){$p(t)$}
}}

\section{Introduction}

In some lectures at the Joseph Fourier Institute in June of 1999,
Lev Rozansky formulated an important conecture concerning the structure of 
the Kontsevich integral \cite{Roz}, and mentioned something of a related
program for a ``finite-type theory of knots' complements''.

By means of some notation let us now describe what of this will be proved 
in this paper.

A {\bf generating diagram} is a diagram with oriented trivalent vertices
(which is to say that the incident edges at a trivalent vertex 
are cyclically
ordered)
and edges decorated with oriented 
bivalent vertices labelled by elements of $\Qset[[k]]$, the ring of formal
power series in a variable $k$. A generating
diagram represents a series of elements of $\CB(k)$ by expanding these 
power series into series of diagrams, as follows. 
If $f(k)\in \Qset[[k]]$ is
\[
f = f_0 + f_1 k + f_2 k^2 + f_3 k^3 \ldots,
\]
where $f_i\in \Qset$, then an edge labelled with $f(k)$
is to be expanded as follows.

\

\[
\fcoupon\hspace{0.05cm}= \hspace{0.05cm}f_0\hspace{0.25cm}
\fcoupa+\hspace{0.15cm}f_1\hspace{0.25cm} \fcoupb
\hspace{1cm}+\hspace{0.15cm}f_2\hspace{0.25cm} \fcoupc
\hspace{1cm}+\hspace{0.15cm}f_3\hspace{0.25cm} \fcoupd \hspace{1cm} \ldots
\]\vspace{0.75cm}

\begin{remark}
The incoming edges at the label (which Definition \ref{wcoupdefn} will
introduce as a ``winding coupon'') are ordered, which determines the
orientation of the introduced trivalent vertices, as shown. The opposite
ordering with the label $f(-k)$ gives the same series. 
\end{remark}

\begin{definition}
Define the {\bf gd-degree} of a generating diagram to be half the number of
trivalent vertices of the (original) diagram.
\end{definition}

\begin{definition}\label{alexdefn}\label{Alexdefn}
Take a knot $K$ in an integral homology three-sphere $M$. Let $A_{(M,K)}(t)$
denote the Alexander
polynomial of the pair $(M,K)$ 
fixed by the requirement that it
satisfy:
\begin{enumerate}
\item{$A_{(M,K)}(t) = A_{(M,K)}(\frac{1}{t}),$}
\item{$A_{(M,K)}(1) = 1.$}
\end{enumerate}
\end{definition}

Let $\Qset^1(t)$ be the ring of rational functions in $t$ that are non-singular
at 1. Denote the inclusion
\[
\iota  :  \Qset^1(t) \hookrightarrow \Qset[[k]],
\]
defined by subtituting $e^{k}$ into $t$.

\begin{definition}
Let $L_{(M,K)}$ be the $\Qset$-vector subspace of $\Qset[[k]]$ that is the
image under $\iota$ of rational functions of the form
\[
\frac{P(t)}{A_{(M,K)}(t)}
\]
where $P(t)\in \Qset[t,t^{-1}]$. This notation can be read as the space of labels.
\end{definition}

The next two definitions require the ``wheel with $2n$ spokes'':

\[
\wh_2 =\ \ \ \ \ \ \ \ 
\Picture{
\thicklines
\DottedCircle
\qbezier[6](1.95,0)(1.5,0)(1,0)
\qbezier[6](-2,0)(-1.5,0)(-1,0)}\ \ ,\ 
\wh_4 =\ \ \ \ \ \ \ \ 
\Picture{
\thicklines
\DottedCircle
\qbezier[6](1.95,0)(1.5,0)(1,0)
\qbezier[6](-2,0)(-1.5,0)(-1,0)
\qbezier[6](0,1)(0,1.5)(0,2)
\qbezier[6](0,-1)(0,-1.5)(0,-2)}\ \ ,\ 
\wh_6 =\ \ \ \ \ \ \ \ 
\Picture{
\thicklines
\DottedCircle
\qbezier[6](1.95,0)(1.5,0)(1,0)
\qbezier[6](1,1.73)(0.75,1.3)(0.5,0.87)
\qbezier[6](-1,1.73)(-0.75,1.3)(-0.5,0.87)
\qbezier[6](1,-1.73)(0.75,-1.3)(0.5,-0.87)
\qbezier[6](-1,-1.73)(-0.75,-1.3)(-0.5,-0.87)
\qbezier[6](-2,0)(-1.5,0)(-1,0)}\ \ ,\ \ldots
\lbb{wheels}
\]\vspace{0.2cm}

\begin{definition}
If the rational numbers $b_{2n}$ are determined by the equality:
\[
\sum b_{2n}x^{2n} = \frac{1}{2} \mbox{log}\left(
\frac{\mbox{sinh}( \frac{x}{2} )}{\frac{x}{2}}\right),
\]
then the series $\nu(k) \in \CB(k)$ is defined to be
\[
\mbox{exp}_{\sqcup}( \sum b_{2n} \omega_{2n} ).
\]
\end{definition}
\begin{remark}
This has recently been shown to be $\KIH(U)$, 
the Kontsevich integral of the unknot, 
by Bar-Natan, Le and Thurston \cite{TW}.
\end{remark}
\begin{definition}\label{wheeldefn}\label{wheelsdefn}
Let $Wh(M,K)$ be defined by
\[
Wh(M,K) = \mbox{exp}_{\sqcup}
\left(\left. \left[-\frac{1}{2} \mbox{log}\left( A_{(M,K)}(e^{h})\right) 
\right] \right|_{h^{2n}\rightarrow \omega_{2n}}\right) \sqcup \nu(k),
\]
where the operation indicated is to expand the term inside the square
brackets into a power series in $u$, and then to replace terms like
$ch^{2n}$ by $c \omega_{2n}$. 
\end{definition}

The LMO invariant was introduced by Thang Le, Jun Murakami and
Tomotada Ohtsuki \cite{LMO} (following an earlier investigation
also with Hitoshi Murakami \cite{LMMO}). 
Our (perhaps non-standard) normalisation of the non-surgered component
specialises to the three-sphere as follows:
\begin{equation}\label{LMOspec}
\KI^{LMO}(S^3,K) = \KIH(K).
\end{equation}

This brings us to the rationality conjecture.
In his lectures Rozansky conjectured the $M \simeq S^3$ case of the following
(see also the new paper by Garoufalidis and Rozansky \cite{GR}).
\begin{theorem}
Let $K$ be a zero-framed knot in an integral homology three-sphere $M$. The
LMO invariant of this pair may be represented
\[
\KI^{LMO}(M,K) =  Wh(M,K)\sqcup \mbox{exp}_{\sqcup}(r) \in \CB(k),
\]
where $r = \sum_{m=1}^{\infty} r^{(m)}$ with $r^{(m)}$ a finite 
$\Qset$-linear combination of connected generating diagrams
of $gd$-degree $m$ whose edges are labelled from $L_{(M,K)}$.
\end{theorem}

This conjecture is motivated by an analogous property of the 
coloured Jones function, as shown by Rozansky
(see, for example,
\cite{Roz2}).

The structures we describe which lead to the 
proof of this
depend on a delicate assembly of results from
the literature. A significant debt is to the
papers of the Aarhus group, who are Bar-Natan, Garoufalidis, Rozansky
and Thurston \cite{A1,A2,A3}. 

This theory was
described at the workshop ``Art of Low-Dimensional Topology VII'', January 8th, 2000, 
for which invitation the author thanks Toshitake Kohno. 

A sequel to this
paper \cite{KS} will develop some technical issues raised within this
work.

\begin{note}
Between versions of this paper, a new paper by Garoufalidis and
Rozansky appeared \cite{GR}. We (the author)
suggest reading these papers in tandem, as
their concerns are somewhat complementary.
\end{note}

\begin{acknowledgements}
The author is supported by a Japan Society for the Promotion of Science
Postdoctoral Fellowship. Thanks to Tomotada Ohtsuki and
the Department of Mathematical and Computing Sciences at 
the Tokyo Institute of Technology
for their support; 
to Kazuo Habiro and Dylan Thurston
for many helpful comments regarding this work; and also to Louis
Funar, Stavros Garoufalidis and Hitoshi Murakami.
\end{acknowledgements}

\newpage

\section{The outline}

\subsection{Special surgery presentations}\label{outlineintro}

\

The strategy of the calculation to be presently described is to apply
the LMO surgery formula to a special surgery presentation which
exists for any knot in a $\ZHS$. The following is well-known.
\begin{lemma}\label{techpreslem}
A zero-framed knot in a $\ZHS$ 
may be obtained from the zero-framed
unkot $U$ in $S^3$ by performing
surgery on some framed link which has the property that
every component of it has linking number 0 with $U$.
\end{lemma}
It may be useful to keep the following example in mind. The figure of 8
knot is obtained by performing surgery on the (blackboard-framed)
component marked with a $*$ below:
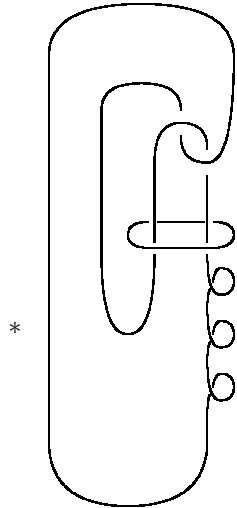
\vspace{2.5cm}

\[
\figateassb
\]
\vspace{4cm}

Thus a knot in a $\ZHS$ can be presented by a framed link in a solid torus
(fixed in $S^3$)
such that every component has linking zero with the core of the torus.

It will prove technically
advantageous to work with a slightly different object: a {\it framed
string link in
the solid torus}. For this definition, realise the
solid torus $ST$ as the complement in the cube
$\{(x,y,z)\in {\Rset}^3; 0\leq x\leq 1, 0\leq y\leq 1, 0\leq z\leq 1\}$
of the hole
$\{(x,y,z)\in \Rset^3; \frac{1}{4}<x<\frac{3}{4}, \frac{1}{4}<y<\frac{3}{4},
0\leq z\leq 1\}$. 

\begin{remark}
We thus use the definite article ``the'', as in ``the solid torus'',
to remind that we are referring to a particular solid torus embedded
in $S^3$.
\end{remark}

\begin{definition}
{\bf A $\mu$-string string link in the solid torus} is a proper
embedding $[0,1]\sqcup\ldots\sqcup [0,1] \hookrightarrow ST$ such that 
the $i$th $\{0\}$ is mapped to $(\frac{i}{\mu+1},0,\frac{1}{2})$, such that
the $i$th $\{1\}$ is mapped to $(\frac{i}{\mu+1},1,\frac{1}{2})$, and with
a framing in the familiar sense of a framed tangle. These are 
identified up to framed isotopies in the solid torus.
\end{definition}
\begin{remark}
We will refer to the $y=0$ plane as {\it the base}, and the $y=1$ plane
as {\it the top}. In this work string links in the solid torus
will always be oriented from the base to the top.
\end{remark}

\begin{definition}
We will call the {\bf meridional disc} the
disc $\{ x\in [\frac{3}{4},1], y=\frac{1}{2}, z\in [0,1] \}$ in the solid
torus.
\end{definition}

To draw a diagram of a string link in a solid torus, we will take
a projection in general position onto the $x-y$ plane, in the familiar
sense. It is convenient to represent the  ``hole'' as
as a fixed dashed loop, or to represent the ``meridional disc''
as a short dashed line segment. 
A diagram will be called
in {\bf general position with respect to the meridional disc} if
it intersects that dashed line transversally. 
For example:
\vspace{2cm}

\[
\figateass\hspace{3cm}\ \mbox{or}\hspace{1.25cm} \figateassc\hspace{2.5cm}
\]
\vspace{3.5cm}

\begin{definition}\label{thrdefn}
\

\begin{enumerate}
\item{Let a marked framed tangle (resp. link) be a framed tangle (resp. link),
possibly with some distinguished components.}
\item{For a string link in the solid torus $T$, let Thr$(T)$ denote the
marked framed tangle obtained
by marking all components,
threading the hole with an unmarked zero-framed unknot (to fix this
let us say we thread on the $x>1$ side),
and then forgetting the hole.}
\item{For a marked framed tangle $T'$, let Clos$(T')$ denote the
marked framed link obtained by closing the tangle.}
\item{For a marked framed link $L$, let KII$(L)$ denote the
class of $L$ modulo Kirby move IIs, where markings indicate to-be-surgered
components (slides of unmarked components over marked components are
also allowed).}
\end{enumerate}
\end{definition}

We now restrict to the presentations guaranteed by Lemma \ref{techpreslem}.
In the following, the adjective special will just mean one of examples
in question. In particular:

\begin{definition}\label{ssldefn}
Let a {\bf special string link in the solid torus} be a string link
in the solid torus $T$ such that:
\begin{enumerate}
\item{Every component has zero algebraic intersection number
with the meridional disc, for diagrams in general position with respect to
the meridional disc.}
\item{The determinant of the linking matrix of the marked components
of \\ Clos(Thr($T$)) is $\pm 1$.}
\end{enumerate}
\end{definition}

\begin{definition}\label{sftdefn}
Let a {\bf special tangle} be a marked, framed 
tangle with one closed, unmarked component;
forgetting that  component leaves the tangle a string link
whose linking matrix has determinant $\pm 1$.
\end{definition}
\begin{definition}
Let a {\bf special link} be a marked, framed link
with one unmarked component; the determinant of the linking
matrix of marked components is $\pm 1$.
\end{definition}


\subsection{The master diagram}

Here is the plan. The basic idea is that formulae for the Kontsevich integral
of a knot can be obtained by applying the LMO invariant to the surgery
presentations just considered. We will see that the resulting factorisation
of the calculation through a certain invariant of string links in the solid
torus has important implications for the result.

The following diagram records this factorisation, as we will now explain.

\begin{diagram}\label{mastercube}
\[
\Picture{
\put(-10,-2){$
\Picture{
\put(0,0){
$\left\{
\begin{array}{c}
Special\ string\ links \\
in\ the\ solid\ torus.
\end{array}
\right\}
$}
\put(20,0){
$\left\{
\begin{array}{c}
Special \\
tangles.
\end{array}
\right\}
$}
\put(17.75,-14){
$\left\{
\begin{array}{c}
Special\ links \\ \hline
Kirby\ move\ IIs.
\end{array}
\right\}
$}
\put(13,-7){$\CB(X,\underline{k})$}
\put(-4,-7){$\CB^{ST}(X)^{Int}$}
\put(13.5,0){\vector(1,0){5}}
\put(15,0.25){Thr}
\put(5,-1.5){\vector(-3,-2){6}}
\put(22,-1.5){\vector(-3,-2){6}}
\put(23.5,-2.5){\vector(0,-1){9}}
\put(24,-7){KII$\ \circ$ Clos}
\put(-2,-8.5){\vector(0,-1){9}}
\put(14.5,-8.5){\vector(0,-1){9}}
\put(-4,-20){$\CB^{QST}(\phi)\times \Zset \times \Zset^1[t,t^{-1}] $}
\put(-3.5,-22){$(S,\sigma,P(t)) \mapsto Thr^{D}(S)\sqcup (-1)^{n\sigma}
{Wh'}(P(t)) $}
\put(13,-20){$\CB_{\leq n}(\underline{k})$}
\put(23,-15.5){\vector(-3,-2){5.75}}
\put(21,-18){$\KI^{LMO,o}_n$}
\put(16,-3.5){$\sigma \circ \KIC $}
\put(-2,-3.5){$\sigma \circ \KIC^{ST} $}
\put(-1,-10.5){$\int^{FG\,in\,ST} dX \times \sigma_+ \times \mbox{det}$}
\put(15,-10.5){$\int^{(n)}dX$}
\put(3,-7){\vector(1,0){6}}
\put(2,-6.25){$S\mapsto Thr^D(S)\sqcup \nu(k)$}
\put(7.75,-19.75){\vector(1,0){4.5}}
\qbezier[24](17,-14)(12.5,-14)(8,-14)
\qbezier[18](8,-14)(5,-14)(2,-17)
\put(2,-17){\vector(-3,-2){0.01}}
}
$}
}
\]
\end{diagram}
\vspace{9cm}

The invariant
\[
\sigma \circ \KIC^{ST} :
\left\{
\begin{array}{c}
Special\ string\ links \\
in\ the\ solid\ torus.
\end{array}
\right\}\rightarrow \CB^{ST}(X)^{Int}
\]
is introduced in Section \ref{stinvsect}. This is an enhancement of the usual
Kontsevich integral, taking values in a space of {\it winding diagrams}.
A winding diagram is, in an appropriate sense, a uni-trivalent diagram
decorated by {\it winding coupons}.

Intuitively speaking, this decoration (modulo some relations) 
describes a homotopy class of proper mappings of that diagram into
the solid torus. In this intuitive picture winding coupons record
intersections of the edges of some representative in general position 
with respect to some fixed meridional disc, with that disc. This
is illustrated in the next figure.

This invariant is more or less pre-existent in the literature, though
our approach and the structures we describe depart from existent works
in certain ways. Our formal presentation, with labelled edges, is closest 
to that of Goryunov \cite{G}, and the intuitive picture is closest  
to that of Andersen-Mattes-Reshetikhin \cite{AMR} (see also Suetsugue 
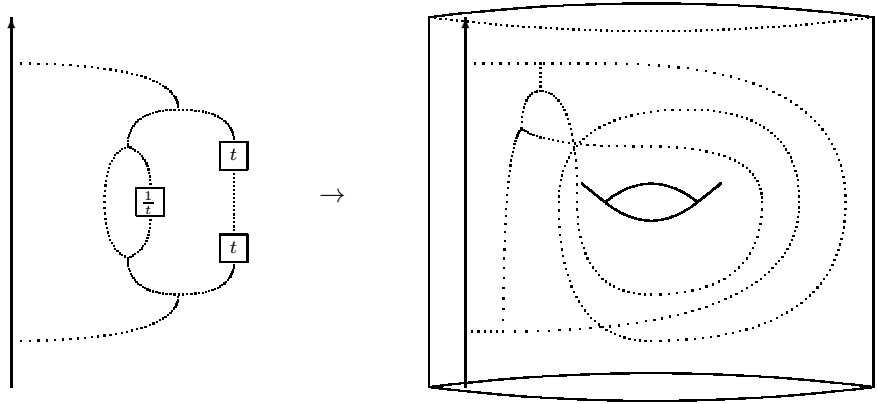
\cite{S}).


\

\vspace{2.5cm}

\setlength{\unitlength}{7pt}
\[
\Picture{
\put(-20,0){$
\windexamp\hspace{3.5cm} \rightarrow \hspace{1.5cm}\windexampb$}}
\]
\setlength{\unitlength}{10pt}
\vspace{3cm}

The destination indicated, $\CB^{ST}(X)$, is a space
of {\it symmetrised winding diagrams}. This plays the part of
Bar-Natan's algebra $\CB$: legs on skeletons are to be
symmetrised. This space is introduced in Section \ref{symsect}.
The map $\sigma$ is the appropriate version 
of Bar-Natan's ``formal Poincare-Birkhoff-Witt'' map.

The utility of this invariant is expressed by the top face of
the cube. Namely, the Kontsevich integral of one of the tangles
of interest, $Thr(T)$, factors through $\KIC^{ST}$. The map which
completes the square is expressed in terms of a map $Thr^D$, {\it threading
diagrams}. This is introduced in Section \ref{threadsect}. This map,
\[
Thr^D : \CB^{ST}(X) \rightarrow \CB(X,\underline{k} )
\]
is the operation of replacing winding coupons (``intersections
with the meridional disc'') with exponentials of legs. The commutativity
of this face depends crucially on a recent calculation due to 
Bar-Natan, Le and Thurston, as is indicated in that section.
\vspace{0.3cm}

\[
Thr^D\left(\ \ \ \ \ \coupon \right)\hspace{0.05cm}= \hspace{0.05cm}\hspace{0.25cm}
\fcoupa+\hspace{0.15cm}\hspace{0.25cm} \fcoupb
\hspace{1cm}+\hspace{0.15cm}\frac{1}{2!}\hspace{0.25cm} \fcoupc
\hspace{1cm}+\hspace{0.15cm}\frac{1}{3!}\hspace{0.25cm} \fcoupd \hspace{1cm} \ldots
\]\vspace{0.75cm}

The subspace $\CB^{ST}(X)^{Int} \subset \CB^{ST}(X)$ is the subspace of
``integrable'' elements (adapting a concept of the Aarhus papers to the
present context), as is defined in Section \ref{intsect}. 
In this case it refers to the subspace of elements of
the form \vspace{0.5cm}

\[
S = exp_{\sqcup}\left( \frac{1}{2} \sum_{i,j} \wcoup \right) \sqcup R,
\]\vspace{0.5cm}

\hspace{-0.45cm}where $W(t) \in M_\mu( \Zset[t,t^{-1}] )$ is a Hermitian matrix of
Laurent polynomials, such that $det(W(1))=\pm 1$, 
and $R$ is a series of diagrams without chords.
In the case at hand ($T$ a string link in the solid torus), 
$\sigma(\KIC^{ST}(T))$ is of this form with matrix
$W(T,t)$, the {\it winding matrix} of $T$, introduced in Section 
\ref{windsect}. 
This is a generalisation of
the notion of linking matrix which incorporates winding information
of the link around the hole of the solid torus.

The mapping $\int^{FG\, in\, ST}dX$ (again, an adaption of a concept
from the Aarhus papers) is defined in the following way. Take an
element $S\in \CB^{ST}(X)^{Int}$, with decomposition as above.
\vspace{0.25cm}

\[
\int^{FG\, in\, ST} dX S = 
\left<
exp_{\sqcup}\left(-\frac{1}{2}\sum_{i,j} \wcoupinv \right), R \right>.
\]
\vspace{0.5cm}

This takes values in the space $\CB^{QST}(\phi)$ of {\it rational winding
diagrams}. This space is in some sense an extension of $\CB^{ST}(\phi)$
which admits rational functions as labels on winding coupons. This mapping
and space are introduced in Section \ref{windsect}.

The front face of the master diagram is detailed in Section 
\ref{surgopsect}. 
The commutativity of the front face indicates that this formula calculates
the LMO invariant (in the event, an extension of a theorem due to the Aarhus
group \cite{A3}). Note that $\sigma_+$ and det are just some normalisation
factors that need to be carried along for the diagram to make sense.

The theorem to take home is the following.
\begin{theorem}[Surgery formula]\label{surgeryformula}
Let a pair of a zero-framed knot $K$ in an integral homology three-sphere
$M$ be presented by $T$, some special string link in the solid torus.
Then $\KI_n^{LMO}(M,K)$ is equal to
\[
\frac{
Wh(M,K) \sqcup Thr^D\left( \int^{FG\, in\, ST}dX \sigma(\KIC^{ST}(T)) \right)}
{ \left((-1)^n\int^{(n)}dU \sigma(\KIC(U_+))\right)^
{\sigma_+(W(T,1))}
\left(\int^{(n)}dU \sigma(\KIC(U_-))\right)^{\sigma_-(W(T,1))}}
\ \ \in\ \ \CB_{\leq n}(\underline{k}).
\]
\end{theorem}

Observe that in this setting ($\Zset HS^3$s), $\KI_n^{LMO}$ is
the degree less or equal to $n$ truncation of the full LMO invariant
$\KI^{LMO}$. We may alternatively present this formla as follows.

\begin{theorem}[Surgery formula, $\beta$-version.]
Let a pair of a zero-framed knot $K$ in an integral homology three-sphere
$M$ be presented by $T$, some special string link in the solid torus.
Then $\KI^{LMO}(M,K)$ is equal to
\[
\frac{
Wh(M,K) \sqcup Thr^D\left( \int^{FG\, in\, ST}dX \sigma(\KIC^{ST}(T)) \right)}
{ \left(\int^{FG}dU \sigma(\KIC(U_+))\right)^
{\sigma_+(W(T,1))}
\left(\int^{FG}dU \sigma(\KIC(U_-))\right)^{\sigma_-(W(T,1))}}
\ \ \in\ \ \CB(\underline{k}).
\]
\end{theorem}

\subsection{Conjecture - a winding diagram valued invariant of knots}

Clearly we are only seeing half of a cube in Diagram \ref{mastercube}.
It will be interesting to describe the other vertex and faces. 
More immediately, the
dashed line
\[
\left\{
\frac{
Special\ links}{
Kirby\ move\ IIs}\right\}
\rightarrow \CB^{QST}(\phi)\times \Zset \times \Zset^1[t,t^{-1}]
\]
would be a consequence of the following conjecture.
\begin{conjecture}\label{invcon}
$
Thr^D : \CB^{QST}(\phi) \rightarrow \CB(\underline{k})
$
is injective.
\end{conjecture}

Actually, this seems clear; we defer a careful explanation of this to
the sequel, which will also discuss the relation with some normalisation
and other technical
issues (see Section \ref{normsect}). Can this corollary be 
proved without reference to the Kontsevich integral?

\begin{corollary}[to the conjecture.]

Take a pair $(M,K)$. Choose $T$, a special string link in the solid
torus presenting $(M,K)$. Then
\[
\frac
{ \int^{FG\, in\, ST} dX \sigma( \KIC^{ST}(T)) }
{\left(\int^{FG} dU \sigma( \KIC( U_+))\right)^{\sigma_+(W(T,1))}
\left(\int^{FG} dU \sigma( \KIC( U_-))\right)^{\sigma_-(W(T,1))}} \in \CB^{QST}(\phi)
\]
is an invariant of the pair $(M,K)$.
\end{corollary}

Strictly speaking, this does not increase our pool of knot invariants. It
may, however, 
be a presentation more appropriate for topological applications.


\subsection{Notation and conventions}\label{notsect}

\

\underline{Spaces of diagrams}

We use standard definitions for spaces of uni-trivalent diagrams \cite{BN}.
There is one point that may be unfamiliar to some readers. To introduce this,
we note that in generality a space may be denoted:
\[
\CA_n(SK ; x_1,\ldots , x_p , \underline{w}_1,\ldots,\underline{w}_q)
\]
which indicates that univalent vertices may:
\begin{enumerate}
\item{be located, 
up to orientation-preserving diffeomorphisms, on a skeleton $SK$,}
\item{or be labelled from $\{x_1,\ldots,x_p,\underline{w}_1,\ldots,\underline{w}_q\}$.}
\end{enumerate}
The underlining of a variable indicates that {\it link relations for that
variable} are to be included in the definition. These were identified in Section 5.2 of 
\cite{A2}.
The point is that link relations are what 
must be included to obtain an
isomorphism:
\[
\CA_n(SK ; x_1,\,\ldots , x_p , \underline{w}_1,\,\ldots,\underline{w}_q)
\simeq 
\CA_n(SK \cup \underbrace{\uparrow\ldots\uparrow}_p
\Picture{
\put(1,0.25){\circle{0.75}}
\put(1.85,0){$\ldots$}
\put(4,0.25){\circle{0.75}}
\put(0.35,-0.2){$\underbrace{\ \ \ \ \ \ \ \ \ \ \ \ \ }_q$}}
\ \ \ \ \ \ \ \ ).
\]
On occasions when there is no skeleton (that is to say that all univalent
vertices are to be labelled), we may use a $\CB$ in place of the $\CA$, following the conventions of \cite{BN}. When there is only one label, then
$\CB(\underline{x})\simeq \CB(x)$.

In this work, we usually work with a set of labels $X=\{x_1,\,\ldots,x_\mu\}$
corresponding to string link components (ultimately providing surgery 
components), and/or a label $k$ corresponding to the (closed) knot component.
Such spaces will be denoted, for example, $\CB(X,\underline{k})$. 

An element of this space may be indicated $S(\overline{x},k)$, where
$\overline{x}$ is thought of as a vector of variables. The logic of
this notation should be clear after reading Section \ref{transsect}.

\

\underline{The invariants}

The notation $\KIH$ denotes precisely the functorial representation
of the category of framed $q$-tangles, according to the definition of 
\cite{LM} (see also \cite{BN2}). We use, unless otherwise stated, the {\it associator with
rational coefficients}. 

Our definition of $\KIC$ differs slightly from existent usage.
In general we will be considering tangles which have a closed component,
which, when forgotten, leaves the tangle a string link. Take such
a tangle $T$, which has $n$ such string link components. Our usage 
of $\KIC(T)$ is:
\[
\KIC(T) = (\underbrace{\nu \otimes \ldots \otimes \nu}_n)\circ 
\Delta^{n-1}(\nu) \circ \KIH(T). 
\]

\

\underline{The LMO invariant}

The LMO group are Thang Le, Jun Murakami and Tomotada Ohtsuki.
In Theorem 6.2 of \cite{LMO} there
is defined an invariant of a link in a three-manifold, which is denoted
there $\Omega(M,L)$. The invariant we use is related to the definition
given there as follows:
\[
\KI^{LMO}(M,K) = \Omega(M,K)\# \nu^{-1}.
\]

Thus, following Proposition 6.5 of \cite{LMO}, we have the restriction
\begin{equation}\label{speceqn}
\KI^{LMO}(S^3,K) = \KIH(K).
\end{equation}

\newpage

\section{A diagram-valued invariant of string links in the solid torus}
\label{stinvsect}

\begin{definition}
A {\bf skeleton} is an oriented one-manifold whose boundary points
are seperated into an ordered pair of ordered sets.
\end{definition}

Our theory employs a certain enhancement of the familiar notion of
uni-trivalent diagram, which will be called a winding diagram. These
will possibly contain a certain new type of vertex which we will call a
winding coupon. This may conveniently be thought of as some decoration
of an underlying uni-trivalent diagram.

\begin{definition}\label{wcoupdefn}
A {\bf winding coupon} is a bivalent vertex whose incoming edges are 
ordered.
\end{definition}

This is depicted as follows. The edges are ordered so that the edge
incoming at the base of $t$ is first, and the edge outgoing at the
top of $t$ is second. We make the convention that a coupon
labelled by $t^{-1}$ is given the reverse orientation.
\vspace{0.7cm}

\[
\coupon\hspace{0.5cm}\mbox{which is equivalent to}\hspace{1.25cm}\coupup
\]
\vspace{0.7cm}

\begin{definition}
A {\bf winding diagram} on a skeleton $SK$ is a graph with univalent vertices,
trivalent vertices and winding coupons, such that:
\begin{enumerate}
\item{(Boundary points.) The set of univalent vertices is seperated
into an ordered pair of ordered sets.}
\item{(Skeleton.) There are distinguished disjoint oriented cycles and
oriented paths between univalent vertices which are
labelled by components of the skeleton $SK$; forgetting
not distinguished edges leaves one with $SK$.}
\item{(Internal verices.) Trivalent vertices which are not met by
distinguished edges are vertex-oriented (have their incoming edges
cyclically ordered).}
\end{enumerate}
Two winding diagrams are identified if there is a graph isomporhism
between them respecting orientations at vertices (including winding
coupons) 
and respecting skeleton information (orientations and labels of distinguished
edges, ordering of boundary points).
\end{definition}

\begin{definition}
The grade of a winding diagram is half the number of trivalent vertices.
\end{definition}

The space in question will be a quotient of the space of finite $\Qset$-linear
combinations of winding coupons of a fixed grade. The quotient will
be by the span of the following classes of vectors. In the relations
MULT and PUSH below,  edges may be part of the skeleton.

\begin{eqnarray*}
AS: & \hspace{1.5cm} & \ASR\hspace{0.5cm}+\hspace{0.85cm}\ASRB \\
& & \\
& & \\
& & \\
& & \\
& & \\
IHX: & \hspace{1.5cm} & \IHX\hspace{1cm}-\hspace{1cm}\IHXB\hspace{1cm}-
\hspace{1cm}\IHXC \\
& & \\
& & \\
& & \\
& & \\
& & \\
STU: & \hspace{1.5cm} & \STUA\hspace{1cm}-\hspace{1.5cm}\STUB\hspace{1cm}+\hspace{1.5cm}\STUC \\
& & \\
& & \\
& & \\
& & \\
& & \\MULT: & \hspace{1.5cm} &
\couponmult\hspace{0.75cm} - \hspace{1.5cm}\unit \\
& & \\
& & \\
& & \\
& & \\
& & \\
& & \\
PUSH: & \hspace{1.5cm} & \pusha \hspace{0.75cm} - \hspace{1.5cm} \pushb \\
& & \\
& & \\
& &
\end{eqnarray*}

\begin{remark}
Looking ahead to the definition of the map $Thr^D$, Definition \ref{threaddefn}, may give some feeling for these orientation conventions and relations.
\end{remark}

\begin{definition}\label{winddefn}
Let $SK$ be a skeleton. Let 
\[
\CA^{ST}_m(SK) =
\frac{
\left\{
\begin{array}{l}
\mbox{Finite $\Qset$-linear combinations of} \\
\mbox{degree $m$  winding diagrams on $SK$}
\end{array}
\right\}
}
{\Qset-\mbox{span of above relations}}
\]
Let $\CA^{ST}(SK)$ denote the completion of $\oplus_{m=0}^{\infty} 
\CA^{ST}(SK)$
with respect to degree.
\end{definition}

\begin{remark}
If $SK$ has $\mu$ components, then $\CA^{ST}_0(SK)$ is isomorphic, as a
$\Qset$-vector space, to $\Qset[t_1^{\pm 1},\ldots,t_\mu^{\pm 1}]$.
\end{remark}
\begin{remark}
If the bottom boundary configuration (regarded as a word in the symbols
$\uparrow$ and $\downarrow$, in the familiar sense) of the skeleton
$K$ matches the top boundary configuration of a skeleton $L$, then
an operation $\circ: \CA^{ST}_n(K)\times \CA^{ST}_m(L) \rightarrow 
\CA^{ST}_{n+m}(K\circ L)$ is obviously defined, and is extended to the
completions. 
\end{remark}

\begin{definition}
For some skeleton $SK$ let
\[
\gamma : \CA(SK) \rightarrow \CA^{ST}(SK),
\]
be the mapping defined by linearly extending the operation of mapping
the element represented by some diagram in $\CA(SK)$ to the element
that that diagram represents in $\CA^{ST}(SK)$.
\end{definition}

\subsection{Some notation}\label{CoupNot}

A coupon labelled by a polynomial represents an element of 
$\CA^{ST}$ via the following expansion. 
At this stage this is best regarded
as a notation; later a space will be introduced ($\CB^{QST}(X)$,
Definition \ref{rationaldefn}) which will admit such labels in its 
definition.

If the polynomial is
\[
p(t) = p_0 + p_1 t \ldots + p_n t^n
\]
then the following expansion is to be understood. 
\vspace{1.25cm}

\[
\pcoupon = p_0\ \ \ \pcouponb +\ \ \ p_1 \ \ \ \ \pcouponc +\ \  \ldots\ \  + p_n\ \ \ \ \ \ \pcouponn
\]\vspace{1.4cm}

It will also prove helpful to have the following diagrammatic.
An oriented edge, labelled with a polynomial, indicates that a coupon
with that label is to be introduced, in the sense just introduced.
The orientation of the coupon is specified by the orientation of the edge.
\vspace{0.5cm}

\[
\nopcoup \hspace{0.5cm} \rightarrow\hspace{1cm} \pcoup
\]\vspace{0.5cm}

\subsection{The invariant}

\begin{definition}
A 4-tuple $(A,B,w_1,w_2)$, where
\begin{itemize}
\item{$A$ and $B$ are q-tangles,}
\item{the bottom boundary word of $A$ is equal to the top boundary 
word of $B$ is equal to $(w_1)(w_2)$,}
\item{the top boundary word of $A$ is equal to the bottom boundary
word of $B$ is equal to $(\ldots(\uparrow\uparrow)\uparrow)\ldots\uparrow)$,}
\end{itemize}
is called a {\bf presentation} 
for $T$, a $\mu$-string string link in the solid torus,
if the result of composing $A$ with $B$, while drilling the hole
at some point on the mutual bounding line between $w_1$ and $w_2$, gives 
$T$. It is clear that every string link in the solid torus has such a 
presentation.
\end{definition}

For example, the string link in the solid torus associated to the
previously considered
surgery presentation of the figure of 8 knot, has
the following presentation:

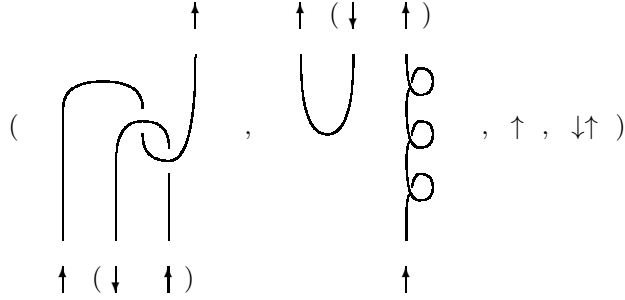
\[
(\put(0,-3){\figatep}\hspace{3cm},\put(0,4){\figatebp}\hspace{3cm}, \ \ \uparrow\ \ , \ \ \downarrow \uparrow\ \ )
\]

\vspace{1.7cm}

For a boundary word $w$, let $G_w$ be the winding
diagram obtained from the identity diagram $I_w$ by attaching
a winding coupon to each strand. This notation can be read as the gluing 
diagram. For example,
\vspace{0.5cm}

\[
G_{(\uparrow \uparrow)\downarrow} =\ \ \ \ \ \ \ \ \ \ \ \ \ \ \ \ \ \ \glueexamp\ \ \ \ \ \ \ \ .
\]

\

\

\begin{definition}
Let $\uparrow^\mu$ denote the skeleton underlying a $\mu$-string
string link.
\end{definition}

\begin{definition}\label{stinvdeaf}
If $T$ is a $\mu$-string string link in the solid torus, then let
\begin{equation}\label{invdeafeqn}
\KI^{ST}(T) = \gamma(\KIH(A_T)) \circ ( I_{w_1} \otimes G_{w_2} ) \circ \gamma(\KIH(B_T))
\in \CA^{ST}(\uparrow^\mu),
\end{equation}
where $(A_T,B_T,w_1,w_2)$ is a presentation for $T$.
\end{definition}

The most pressing issue is, of course, to show that this is well-defined.
For the time being, then, indicate the dependence on the presentation
$\KI^{ST}(A_T,B_T,w_1,w_2)$. The well-definedness
will follow from the following, clear,
observation.
\begin{lemma}\label{sweeplem}
Let $A\in \CA^{ST}(SK)$, take a word $w$ such that $G_w \circ A$ is
well-defined, and take a word $w'$ such that $A \circ G_{w'}$ is well-defined.
Then
\[
G_w \circ A = A \circ G_{w'}.
\]
\end{lemma}

%
%
%
%
%
%
%
%
%
%
%
%
%
%
\begin{lemma}
$\KI^{ST}( A_T, B_T , w_1, w_2 )$ is independent of the
choice of presentation, and hence is an invariant of $T$.
\end{lemma}
\underline{Proof.}

It is clear that any two presentations can be related by a finite 
sequence of the following moves:
\[
\begin{array}{llll}
(1) & ( A \circ (1_{w_1} \otimes C), B, w_1, w_2 )
& \leftrightarrow & (A, (1_{w_1} \otimes C) \circ B, w_1, w_2'), \\
(2) & ( A \circ (C \otimes 1_{w_2}), B, w_1, w_2 )
& \leftrightarrow & (A, (C \otimes 1_{w_2}) \circ B, w_1', w_2).
\end{array}
\]
The lemma then follows from the functoriality of $\KIH$ and Lemma \ref{sweeplem}.
\rtb
Thus we revert to the notation $\KI^{ST}(T)$.

\

The normalisation of this invariant that is
appropriate for surgery considerations
is the following. This is the normalisation of \cite{LMMO}.
\begin{definition}\label{surgnorm}
Let $T$ be an $\mu$-string string link in a solid torus. Define:
\[
\KIC^{ST}(T) = \gamma((\underbrace{\nu \otimes \ldots \otimes \nu}_\mu) \circ \Delta^{\mu-1}(\nu)) \circ \KI^{ST}(T),
\]
in the space $\CA^{ST}(\uparrow^\mu)$, 
recalling that $\nu = \KIH(U) \in \CA(\uparrow)$.
\end{definition}

\subsection{The co-product}

We now equip $\CA^{ST}(K)$ with a co-product. The presence of winding
coupons does not affect the following familiar definition.

\begin{definition}
Take a diagram $D$ such that its dashed graph has connected components
indexed by the set $I$. If $J\subset I$ let $D_J$ indicate the diagram 
obtained by forgetting those components in the subset $J$.
Then, define the mapping $\Delta$ as the linear extension of
\[
\Delta(D) = \sum_{J\subset I} D_J \otimes D_{I-J}.
\]
\end{definition}

\begin{remark} 

\

\begin{enumerate}
\item{If $D$ has an empty dashed graph then this
operation is defined to be $\Delta(D) = D \otimes D$.}
\item{
This defines a co-product on the graded completions:
\[
\Delta : \CA^{ST}(K) \rightarrow \CA^{ST}(K)\, \hat{\otimes}\, \CA^{ST}(K).
\]}
\end{enumerate}
\end{remark}

\begin{lemma}
This is well-defined, co-commutative, co-associative, 
and commutes with compositions. 
\end{lemma}

To see that it is well-defined we must show that relations are mapped
to relations. The only novelty is a PUSH relation when two of the
involved edges are part of the skeleton; this relation is easily checked. 
Observe that it commutes with compositions by construction.
All other properties are standard.

\begin{lemma}\label{glprop}
For a string link in a solid torus $T$,
\[
\Delta( \KIC^{ST}(T) ) = \KIC^{ST}(T)\otimes \KIC^{ST}(T).
\]
\end{lemma}
\underline{Proof.} This follows for the usual reasons: that is, from the
corresponding property for $\hat{Z}$, the corresponding property for the
normalisation factors of Definition \ref{surgnorm},
from the obvious property that $\Delta(G_{w}) = G_{w}\otimes G_{w}$
and from the commutation of composition with the co-product. 
\rtb

\subsection{The Hopf algebra
$\CB^{ST}(X)$ }
\label{symsect}

\

We turn to the case of special string links in the solid torus. These are,
remember, string links in the solid torus such that a representative in general
position with respect to the meridional disc has algebraic intersection
zero with it.
For such an $\mu$-string string link in the solid torus $T$, $\KIC^{ST}(T)$
lies in a special subspace of $\CA^{ST}(\uparrow^\mu)$. 
\begin{definition}
Let $\CA^{ST,spec}(\uparrow^\mu)$ denote the subspace of $\CA^{ST}(\uparrow^\mu)$
spanned by diagrams with the property that the product of all
the labels on the winding coupons labelling 
some component of the skeleton is 1 (that is, using a factor of $t^{-1}$
if some coupon is oriented against the orientation of that component);
for each component.
\end{definition}

\begin{observation}
If $T$ is a special $\mu$-string string link in the solid torus
then
\[
\KIC(T) \in \CA^{ST,spec}(\uparrow^\mu).
\]
\end{observation}

\begin{remark}
If a diagram is in this subspace, then repeated 
applications of the $PUSH$ relation can be used to make all of the labels
on the skeleton 1 (say, by pushing all the labels to one end).
Then all the labels (all the {\it winding}) will be carried by the
dashed graph. 
\end{remark}

We now introduce an isomorphic description of this subspace. 
This is an enhancement of the familiar algebra $\CB$. Let
$X = \{x_1,\ldots,x_\mu\}$ denote a labelling set for the skeleton $\uparrow^\mu$.
\begin{definition}
Let a {\bf winding diagram on $X$}
be a graph with oriented trivalent vertices, winding coupons, and 
univalent vertices labelled from $X$.
\end{definition}

We may alternatively call this a {\bf symmetrised winding diagram} on 
$\uparrow^{\mu}$.

\begin{definition}
Define
\[
\CB_m^{ST}(X) = 
\frac{
\left\{
\begin{array}{c}
\mbox{Finite $\Qset$-linear combinations of degree $m$} \\
\mbox{winding diagrams on $X$}
\end{array}
\right\}}
{\Qset-\mbox{span of AS, IHX, OR, MULT and PUSH relations}}
\]
\end{definition}
Let $\CB^{ST}(X)$ denote the graded completion of $\oplus_{m=0}^{\infty}\CB^{ST}_m(X)$.
Equip this with the obvious analogs of the 
``disjoint-union'' product, the
``sum over partitions into two sets'' co-product, the ``empty set'' unit
and co-unit, and the ``$(-1)$ for every component'' antipode.

\begin{lemma}
$\CB^{ST}(X)$ is a commutative, co-commutative Hopf algebra.
\end{lemma}

\begin{definition}
Let $\chi : \CB^{ST}(X) \rightarrow \CA^{ST,spec}(\uparrow^\mu)$
be the operation defined on some symmetrised diagram of taking the average
of all diagrams obtained by locating all univalent vertices labelled with
$x_1$ on the first component, etc. ; linearly extended to each $\CB^{ST}_m$,
and to $\CB^{ST}$.
\end{definition}
\begin{lemma}
The mapping $\chi$ describes a $\Qset$-vector space isomorphism at each 
grade
\[
\CB^{ST}_m(X) \simeq \CA^{ST,spec}_m(\uparrow^\mu),
\]
commuting with coproducts
\[
(\chi \otimes \chi) \circ \Delta = \Delta \circ \chi.
\]
\end{lemma}
The inverse, $\sigma : \CA^{ST,spec}_m(\uparrow^\mu) \rightarrow 
\CB^{ST}_m(X)$, 
is obtained by first 
pushing 
all the winding coupons onto the dashed graph (say by pushing all coupons on
some component to the top of that component),
and then applying 
Bar-Natan's ``formal PBW map'' \cite{BN}. It is straightforward to check
(adapting \cite{BN})
that these maps are well-defined and inverses of each other.

The conclusion of this development follows:
\begin{lemma}\label{expform}
Let $T$ be a special $\mu$-string string link in the solid torus. Then
$\sigma( \KIC^{ST}(T) )$ is a group-like element in the Hopf algebra 
$\CB^{ST}(X)$. Thus it
is an exponential of a series of
connected diagrams, a finite $\Qset$-linear combination at each grade.
\end{lemma}

To see this, note that Lemma \ref{glprop} indicates that $\KIC^{ST}(T)$ is
group-like in $\CA^{ST,spec}$. Thus its image in $\CB^{ST}(X)$
is also group-like because of the commuting of the co-product with
the map $\sigma$. Thus it is an exponential of a primitive element
(for example, \cite{Qu}, Appendix A):
at each grade this will be a finite $\Qset$-linear combination of connected 
diagrams.

\subsection{The winding matrix}\label{windsect}

We now 
introduce $W(T,t) \in M_\mu( \Zset[t,t^{-1}] ),$ the 
{\bf winding matrix} of $T$, where $T$ is a $\mu$-string string link in
the solid torus.
Number the components of $T$ and choose a diagram for $T$ that is in general 
position with respect to the meridional disc.  We consider paths on this
diagram.
The ``algebraic
intersection of a path with the meridional disc'' is the sum over all crossings
of that path with the disc of: 
a plus one if the tangent vector of the path points in the 
direction of increasing $y$ at the intersection; 
and a minus one otherwise (according to the model of the solid
torus
described in Section \ref{outlineintro}).

\begin{definition}
For a crossing $c$ of strands $i$ and $j$, let $\epsilon(i,j,c)$ 
denote the algebraic intersection with the meridional disc of the path obtained
by travelling from the base along $i$ to $c$, crossing to $j$, and then
travelling to the top along $j$. (To be precise, if $i=j$ then change
strands the {\it first} time the crossing is encountered).
\end{definition}

\begin{definition}

For $T$, an $\mu$-string string link in the solid torus, choose a 
(blackboard-framed) diagram
for $T$ in general position with respect to the meridional disc. If
$c$ denotes a crossing then let $sgn(c)$ denote the sign of that crossing.
Let:

\[
{W}_{ij}(T,t) = 
\left\{
\begin{array}{ll}
\sum_{\mbox{\small c, a crossing of i and j}}\ \ \  \frac{1}{2}\mbox{sgn}(c)t^{\epsilon(i,j,c)}
& \mbox{if}\ i\neq j, \\
& \\
\sum_{\mbox{\small c, a self-crossing of $i$}}\ \ \  \frac{1}{2}\mbox{sgn}(c)
({t^{\epsilon(i,i,c)} + t^{-{\epsilon(i,i,c)}}})& \mbox{otherwise.}
\end{array}
\right.
\]
\end{definition}

\begin{remark}\label{Hermremark}

\

\begin{itemize}
\item{
The winding matrix is Hermitian:  ${W}_{ij}(T,t) = {W}_{ji}(T,t^{-1})$.}
\item{
The winding matrix specialises to the linking matrix of the underlying tangle: 
$W_{ij}(T,1) = \mbox{Lk}_{ij}(T)$.}
\end{itemize}
\end{remark}

\begin{lemma}
$W_{ij}(T,t)$ is an isotopy invariant of $T$, regarded
as a framed string link in a solid torus.
\end{lemma}
We will return to the topological interpretation of $W_{ij}$ 
at the end of this section when we make an
important connection with the Alexander polynomial. Besides, its isotopy
invariance follows from its appearance in the following theorem (the unique
Hermitian matrix with the following property).


Let $X=\{x_1,\ldots,x_\mu\}$ denote a labelling set for $T$.

\begin{theorem}\label{windfirst}
If $T$ is a string link in a solid torus then
\vspace{0.65cm}

\[
\sigma( \KIC^{ST}(T) ) = \mbox{exp}_{\sqcup}(\  \frac{1}{2}\sum_{i,j} \wcoup 
)\sqcup R,
\]
\vspace{0.65cm}

where $R$ is a series of $X$-substantial diagrams.
\end{theorem}

\underline{Proof.}

Lemma \ref{expform} indicated that 
$\sigma( \KIC^{ST}(T) )$ is of the form exp$_\sqcup(S)$ where $S$ is a series
of connected diagrams. Thus to prove the theorem it is
sufficient to calculate the degree one part of $\KIC^{ST}$. 

Take a diagram of $T$ in general position with respect to the meridional
disc, and an associated presentation of $T$. 
Examining the definition of $\KIC^{ST}$, 
we see that there will be a contribution of one chord from 
every crossing in this diagram.
The introduced complexity is that there will be a distribution of
winding coupons on the skeleton.
(The reader is invited to consider the example that follows
this proof).

Consider first a crossing between two different components, $i$ and
$j$. Use the relation $PUSH$ to push the coupons onto the chord 
in the following
manner: coupons that occur before the crossing as $i$ is traversed from the 
base should be pushed past the chord (where they will all cancel to 1);
coupons that occur after the crossing as $j$ is traversed to the top should
be pushed back past the chord (where they will all cancel to 1).

The reader can check that the label on the chord coming from that
crossing is precisely as follows. Note that in this section we will
use the notational convention described in Section \ref{CoupNot}.

\

\[
\Picture{
\put(-16,0){$
\frac{1}{2}\mbox{sgn}(c)\hspace{2cm}\crossingchord\hspace{1cm}$}
\put(-14,-10){$=\frac{1}{2}\left( \frac{1}{2}\mbox{sgn}(c)\hspace{2cm}
\crossingchordb\hspace{1cm}\ \ +\ \ \frac{1}{2}\mbox{sgn}(c)\hspace{2cm}
\crossingchordc\hspace{1cm} \right) $}}
\put(-16,-19){$
=
\chi(\hspace{0.15cm} 
\frac{1}{2}(\ \frac{1}{2}\,\mbox{sgn}(c)\hspace{0.5cm}\selfap \hspace{0.15cm}+\ \ \frac{1}{2}\,\mbox{sgn}(c)\hspace{0.5cm}\selfbp \hspace{0.25cm} ) )$}
\]

\vspace{0.85cm}

Now consider self-crossings of components. Choose some self-crossing 
$c$ of some component $x_i$. Its contribution is
as follows:

\

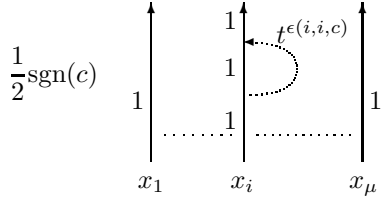
\[
\frac{1}{2}\mbox{sgn}(c)\ \ \ \ \ \ \ \ \ \ \ \ \ \ \ \ \ \ \crossingchordself\ \ \ \ \ \ \ \ \ \ 
\]

\vspace{1.5cm}

Now, use the following relation:

\

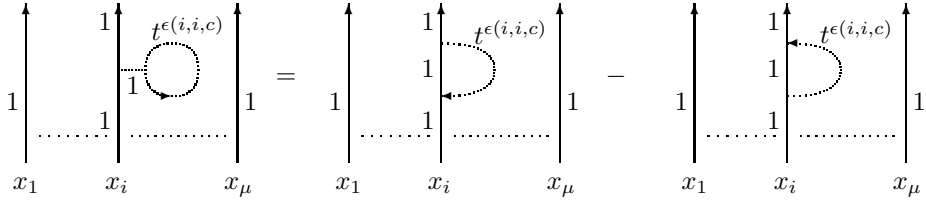
\[
\ \ \ \ \ \ \ \ \ \ \ \crossingchordselfloop\
\hspace{1cm}=\hspace{2cm}
\crossingchordselfb\hspace{1cm}\ \ -\ \ \hspace{2cm}
\crossingchordself 
\]

\vspace{1.6cm}

to write the contribution from this crossing as
follows:

\

\[
\chi(\hspace{0.25cm} 
\frac{1}{2}\mbox{sgn}(c)(\hspace{0.75cm}
\selfa \hspace{0.25cm}-\ \ \frac{1}{2}\hspace{0.75cm} \selfb
\hspace{0.5cm}))
\]
\vspace{0.5cm}

In other words:

\[
\chi\left(\hspace{0.15cm}\frac{1}{2}\,\left(\,
\frac{1}{2}\mbox{sgn}(c)
\left(\hspace{0.75cm} \selfa \hspace{0.25cm}+\ 
\hspace{0.75cm} \selfc\hspace{0.5cm}\right)\right) + r\ \ \right),
\]
\vspace{0.5cm}

where $r$ is an $X$-substantial diagram.
\rtb
\begin{example}
Consider the example associated with the presentation of the figure of 8
knot.
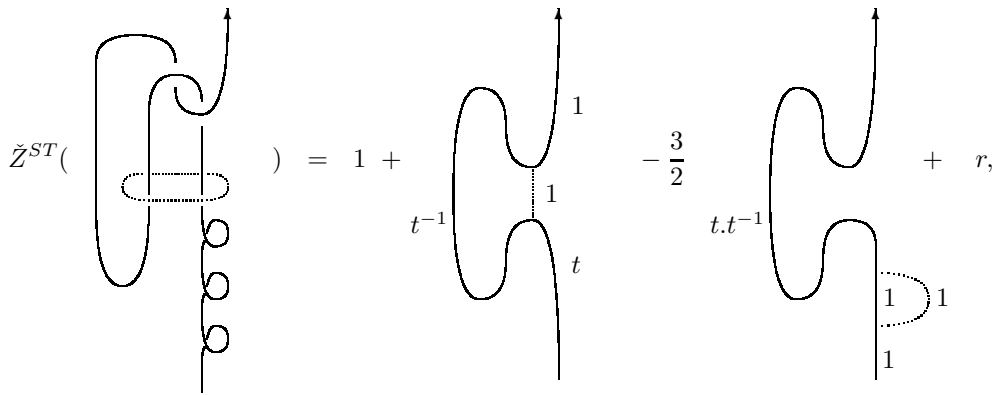
\vspace{1.5cm}

\[
\KIC^{ST}(\figateass
\hspace{2cm})\ \ =\ \ 1\ + \hspace{0.25cm} \figatediag\hspace{1.25cm}
\ \ \ \ \ \ \ -
\frac{3}{2}\hspace{0.75cm} \figatediagb\hspace{1.25cm}\ \ \ +\ \ \ r,\hspace{2cm}
\]

\vspace{5cm}

\[
\hspace{4cm} =1 \ \ +\hspace{0.5cm}\figatediagc\hspace{2.25cm}\ \ \ \ \ \ \ \ +\ \ \ r',\hspace{2cm}
\]
\end{example}

\vspace{1.5cm}

\hspace{-0.4cm}where $r$ and $r'$ are series of diagrams that are either of grade
greater than 1 or X-substantial.

\rtb

\vspace{0.5cm}

\subsection{A topological interpretation of $W(T,t)$ }

The attentive reader will have noticed the appearance of the 
Alexander polynomial for the figure of 8 knot in the previous calculation.
Let us examine the meaning of the matrix $W(T,t)$
more closely.

Choose a base point close to the base of the strings, 
and choose paths from that basepoint to the bases of the strings
(so that the bases, paths, and basepoint all lie
in some ball). Then close the string link on the left, obtaining
some link in the solid torus
$L$, say with components $\{K_1,\ldots,K_\mu\}$, with a path from
some basepoint to some point on each component.

Take the universal cyclic cover of the solid torus:
\[
p: \widetilde{ST} \rightarrow ST,
\]
and lift the link $\{K_1,\ldots,K_\mu\}$ to $\widetilde{ST}$. This can be done 
as we are restricting to special string links (that is, string links
whose algebraic intersection with any meridional disc is zero).
The group of translations is $\Zset$: choose an action such that a path
which starts at some point $p$; crosses the meridional disc in the direction
on increasing $y$; and the returns to the $p$ (without again crossing
the meridional disc) is lifted to a path starting at some $p$ and finishing
at some $tp$. The lifted link $\widetilde{L}$ can be
identified as the set of translates
\[
\{\ldots,t^{-1} \widetilde{K_1}, \widetilde{K_1}, t\widetilde{K_1},\ldots,
t^{-1}\widetilde{K_\mu},\widetilde{K_\mu},t \widetilde{K_\mu},\ldots \},
\]
where the components $\{ \widetilde{K_1},\ldots, \widetilde{K_\mu}\}$
can be fixed by choosing $\widetilde{K_1}$ and then following the
lifts of the arcs introduced when the closure was taken.

We define an invariant of $T$ using this arrangement, as follows.
\begin{definition}
Let $\widetilde{\mbox{Lk}}(T) \in M_\mu( \Zset[t,t^{-1}] )$ be
defined by 
\[
\widetilde{\mbox{Lk}}_{ij}(T) = \sum_{m=-\infty}^{\infty} t^{m}
lk( \widetilde{K}_i, t^m \widetilde{K}_j ).
\]
\end{definition}

\begin{lemma}
\[
W(T,t) = \widetilde{\mbox{Lk}}(T).
\]
\end{lemma}
\underline{Proof.}

Take a blackboard-framed
diagram for $T$ that is in general position with respect to the 
projection of the meridional disc.
A fundamental domain for a diagram for
the pair $(\widetilde{ST}, \widetilde{L})$ can be
obtained by cutting such a diagram along the projection of the meridional
disc, obtaining a rectangle, and then gluing a countable infinity of
copies of that rectangle in the appropriate way. 

Observe that a crossing between two different components, say $x_i$ and $x_j$,
in the diagram for $T$ 
lifts to a crossing between $K_i$ and some translate of $K_j$. Which 
translate is decided by counting intersections with the meridional disc
as in the definition of $W$.

The slightly different definition of $W$ along the diagonal is accounted
for by the fact that
self-crossings of a component $x_i$ in the diagram for $T$ will either
lift to a self-crossing of $K_i$, or to a {\it pair} of 
crossings, between $K_i$ and some $t^aK_i$ and between $K_i$ and $t^{-a}K_i$.
Observe that the factors give the right weights in both situations.
\rtb

\

\begin{example}
Continuing the example of the figure of 8 knot:
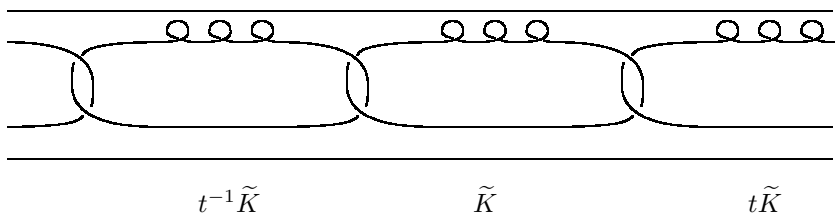
\vspace{0.5cm}

\setlength{\unitlength}{8pt}
\[
\figatesass
\]
\setlength{\unitlength}{10pt}
\end{example}
\vspace{2cm}

\rtb

\

Take $T$, a special string link in the solid torus, presenting some
pair $(M,K)$.
Remember that, according to Definition \ref{Alexdefn},  
$A_{(M,K)}(t)$ denotes the canonical Alexander polynomial of a knot in a
$\Zset HS^3$.

\begin{lemma}\label{alexconn}
\[
A_{(M,K)}(t) = \pm \mbox{det}( W(T,t) ).
\]
\end{lemma}
\underline{Proof.}

Surgery on the framed link $\widetilde{L}$ in $D^2\times \Rset$ 
recovers the universal cyclic cover of the complement of
$K_T$. 

The Mayer-Vietoris sequence then indicates that the matrix
$\widetilde{L}$, and hence $W(T,t)$, is a presentation matrix for
the $\Zset[t]$-module $H_1( \widetilde{M - K} ; \Zset )$.

The Alexander polynomial is defined as a generator of the order ideal
of that module. In the situation at hand, given a square presentation matrix,
this is calculated by the determinant of that matrix,
as in the statement of the lemma. Note that this only specifies the
polynomial up to multiplications of the form $\pm t^{n}$, and the
statement of the lemma asserts that the recovered polynomial 
is {\it symmetric} under $t\rightarrow t^{-1}$. 

To see that this is true note that $W(T,t)$ is
a Hermitian matrix, 
according to Remark \ref{Hermremark}.
\[
\mbox{det}(W(T,t))
=
\mbox{det}(W(T,t)^{Tr.}) = \mbox{det}(\overline{W(T,t)})
\]
\rtb

\begin{remark}
$W(T,t)$ is, presumeably, the matrix referred to in Exercise C.13 of
Rolfsen \cite{Rol}.
\end{remark}

\subsection{Alternative normalisations}\label{normsect}

We take this oppurtunity to draw attention to a certain subtle choice
of normalisation that has been made in the construction of
$\KI^{ST}$ given here.

Let $\alpha$ denote a group-like element of $\CA(\uparrow)$, the 
space of uni-trivalent diagrams on a single strand.
\begin{definition}
Let 
\[
\KI^{ST}[\alpha]:
\left\{
\begin{array}{c}
String\ links\ \\
in\ the\ solid\ torus.
\end{array}
\right\}
\rightarrow
\CA^{ST}(\uparrow^\mu).
\]
be defined in exactly the same way as given in Definition \ref{stinvdeaf},
except that in Equation \ref{invdeafeqn} $I_{\omega_1}\otimes G_{\omega_2}$
should be replaced by
\[
I_{w_1}\otimes ( \Delta_{\omega_2}(\alpha)\circ G_{\omega_2} ),
\]
where $\Delta_{\omega_2}(\alpha)$ is the paralleling operation
across the strands described by $\omega_2$, applied to $\alpha$.
Let $\KIC^{ST}[\alpha]$ denote the normalisation corresponding to 
Definition \ref{surgnorm}.
\end{definition}
\begin{remark}
Lemma \ref{expform} and Theorem \ref{windfirst} still hold if
$\KIC^{ST}$ is replaced by $\KIC^{ST}[\alpha]$.
\end{remark}

Certain choices of $\alpha$ will prove appropriate for certain applications.
For example, setting $\alpha=\nu^{-1}$ gives a normalisation that is better
adapted to questions involving covering spaces.

Applying $\int^{FG\, in\, ST}dX$, followed by $Thr^D$, 
presumeably gives knots invariants
which are different normalisations of the Kontsevich integral,
in some sense. We will return to this question in
the sequel.

\newpage

\section{Formal Gaussian integration}\label{intsect}
In this section we focus on the leftmost edge of Diagram \ref{mastercube}.
\[
\Picture{
\put(-2,-1){$
\Picture{
\put(-1,0){$\CB^{ST}(X)^{Int}$}
\put(1,-1){\vector(0,-1){4}}
\put(-3.5,-7){$\CB^{QST}(\phi)\times \Zset \times \Zset^1[t,t^{-1}]$}
\put(2,-3){$\int^{FG\,in\,ST}dX \times \sigma_+ \times \mbox{det}$}
}$}}
\]\vspace{2.5cm}

\begin{definition}\label{integdefn}
An element $S\in \CB^{ST}(X)$ is said to be {\bf integrable} if it
is of the following form.\vspace{0.1cm}

\begin{equation}\label{decompeqn}
S = \mbox{exp}_\sqcup \left( \frac{1}{2}\sum_{ij} \wcoup \right)\sqcup R,
\end{equation}
\vspace{0.5cm}

\begin{enumerate}
\item{$W_{ij}(t)$ is a Hermitian matrix ($W_{ij}(t) = W_{ji}(t^{-1})$)
wih the property that det$(W(1))=\pm 1$,}
\item{$R$ is a series of $X$-substantial diagrams: that is, a diagram
will have no chords (ignoring winding coupons).}
\end{enumerate}
\end{definition}

\begin{remark}
Observe that the above decomposition is unique. In particular, a
Hermitian matrix satisfying Equation \ref{decompeqn} (for some
given $S$) will be unique. The matrix will be called the {\bf Gaussian
matrix} of $S$.
\end{remark}

\begin{definition}
Let $\CB^{ST}(X)^{Int}$ denote the subspace of integrable elements
of $\CB^{ST}(X)$.
\end{definition}

\begin{remark}\label{matrixeval}
As the associated Hermitian matrix is uniquely specified, we will
freely apply any function of such matrices to the set $\CB^{ST}(X)^{Int}$
with the understanding that the function is to be applied to the 
Gaussian matrix of the element. An example of such use is the use of 
the functions $\sigma_+$ and det in the diagram above. These factors
will be taken up again in Section \ref{surgopsect}.
\end{remark}

\subsection{Rational winding diagrams}\label{rationaldefn}

The Aarhus calculation of the LMO invariant has an associated philosophy
of formal Gaussian integration. The idea is that given an integrable
element of $\CB(X)^{Int}$ one ``integrates'' by splitting off the quadratic
part, inverting it, and contracting the result with the remainder.

We would like to introduce an analog of this in the situation at hand,
that is, for the space $\CB^{ST}(X)^{Int}$. Whilst we have a clear
definition for an integrable element, it remains for us to
introduce a suitable space in which to ``invert'' the given
matrix of polynomials. 

The definition that follows generalises the 
definition of the space of winding diagrams, Definition \ref{winddefn}.
We will differentiate this space in our vocabulary by inserting the
adjective ``rational''. 

\begin{remark}

We should point out that for the proof of Theorem \ref{surgeryformula},
and hence Rozansky's conjecture, the following definition is unneccessary.
The reader who feels the following is an unneccessarily cumbersome space 
may happily evade this space by observing that, in this work,
every appearance of $\int^{FG\,in\,ST}dX$ is followed by a $Thr^D$.
(See Remark \ref{evasion}). 
The point of including such a definition here is 
Conjecture \ref{invcon}.
\end{remark}

\begin{definition}
A {\bf rational winding coupon} 
is a vertex of some even valency. The incoming edges
are ordered up to reorderings of the form $(\sigma,\sigma)$ for some
$\sigma \in \Sigma_m$ (if the valency is $2m$). 
A winding coupon is labelled from $\Qset^1(t),$ the ring of rational
functions in a single variable which are non-singular at 1.
\end{definition}

This is depicted as follows, the idea being that an edge leads through
the coupon to the opposite edge. 
The ellipsis used as below will 
indicate the possibility of
a number of other edges. The ordering of the edges is the bottom row, from
left to right, followed by the top row, from left to right.

\

\[
\qunit
\]
\vspace{0.25cm}

\begin{definition}
A {\bf rational 
winding diagram} labelled from the set $X$
is a graph with univalent vertices,
trivalent vertices and winding coupons, such that:
\begin{enumerate}
\item{(Univalent vertices.) Univalent vertices are labelled from $X$.}
\item{(Trivalent verices.) Trivalent vertices 
are vertex-oriented (have their incoming edges
cyclically ordered).}
\end{enumerate}
Two rational winding diagrams are identified if there is a graph isomporhism
between them respecting orientations at trivalent vertices, 
orientations and labels of winding 
coupons, and labels of univalent vertices.
\end{definition}

\begin{definition}
The grade of a rational 
winding diagram is half the number of trivalent vertices.
\end{definition}

The space in question will be a quotient of the space of finite $\Qset$-linear
combinations of rational winding diagrams of a fixed grade. The quotient will
be by the span of the following classes of vectors.
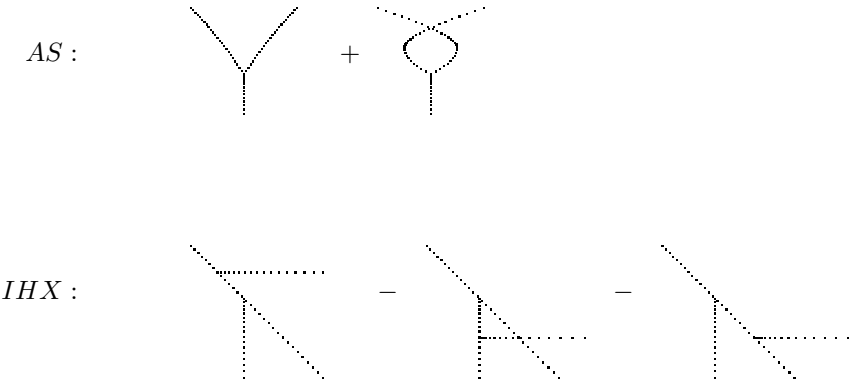
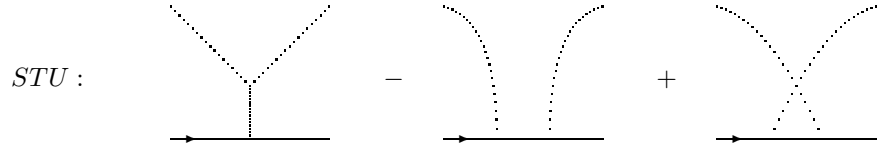
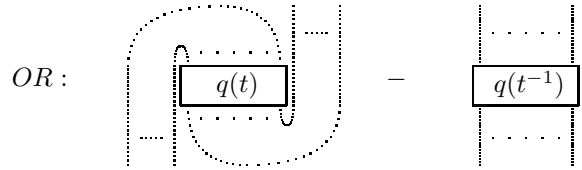
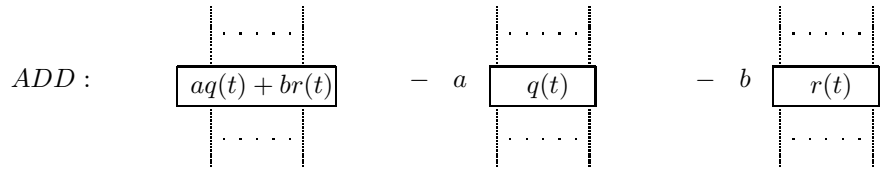
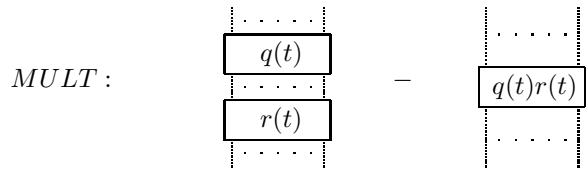
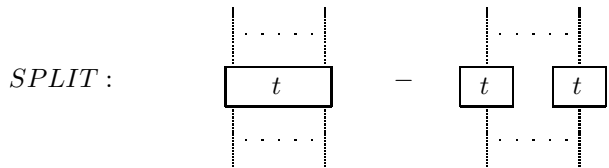
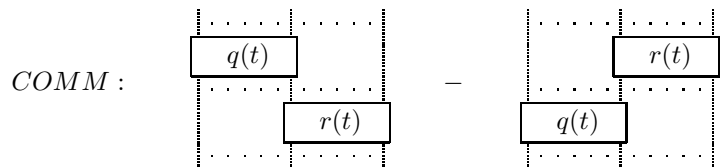
\vspace{0.75cm}

\begin{eqnarray*}
AS: & \hspace{1.5cm} & \ASR\hspace{0.5cm}+\hspace{0.85cm}\ASRB \\
& & \\
& & \\
& & \\
& & \\
& & \\
IHX: & \hspace{1.5cm} & \IHX\hspace{1cm}-\hspace{1cm}\IHXB\hspace{1cm}-
\hspace{1cm}\IHXC \\
& & \\
& & \\
& & \\
& & \\
& & \\
& & \\
& & \\
STU: & \hspace{1.5cm} & \STUA\hspace{1cm}-\hspace{1.5cm}\STUB\hspace{1cm}+\hspace{1.5cm}\STUC \\
& & \\
& & \\
& & \\
& & \\
& & \\
OR: & \hspace{1.5cm} &
\qrev \hspace{1.25cm} - \hspace{1.5cm}\invunit \\ 
& & \\
& & \\
& & \\
& & \\
& & \\
ADD: & \hspace{1.5cm} &
\addunit \hspace{1.25cm} - \ \ a \hspace{1cm}\qlong 
\hspace{1.25cm} -\ \ b \hspace{1cm}\rlong \\ 
& & \\
& & \\
& & \\
& & \\
& & \\
MULT: & \hspace{1.5cm} &
\qtimesr\hspace{0.75cm} - \hspace{1.5cm}\qrunit \\
& & \\
& & \\
& & \\
& & \\
& & \\
SPLIT: & \hspace{1.5cm} &
\tunit \hspace{0.75cm} - \hspace{1.5cm}\manyq \\ 
& & \\
& & \\
& & \\
& & \\
& & \\
COMM: & \hspace{1.5cm} &
\qcomm \hspace{1.25cm} - \hspace{2cm}\qcommb \\ 
& & \\
& & \\
& & \\
& & \\
& & \\
& & \\
& & \\
PUSH1: & \hspace{1.5cm} &
\qforka \hspace{0.75cm} - \hspace{1.5cm}\qforkb \\ 
& & \\
& & \\
& & \\
& & \\
& & \\
& & \\
PUSH2: & \hspace{1.5cm} &
\qloopa \hspace{0.75cm} - \hspace{1.5cm}\qloopb \\ 
& & \\
& & \\
& & \\
& & 
\end{eqnarray*}

\begin{definition}
Let $X$ be a set of labels. Let 
\[
\CB^{QST}_m(X) =
\frac{
\left\{
\begin{array}{l}
\mbox{Finite $\Qset$-linear combinations of degree $m$} \\
\mbox{rational winding diagrams labelled from $X$}
\end{array}
\right\}
}
{\Qset-\mbox{span of above relations}}
\]
Let $\CB^{QST}(X)$ denote the completion of $\oplus_{m=0}^{\infty} \CB^{QST
}(X)$ with respect to degree.
\end{definition}

\subsection{Formal Gaussian integration in the solid torus.}

We can now introduce the map
\[
\int^{FG\, in\, ST}dX : \CB^{ST}(X)^{Int} \rightarrow \CB^{QST}(\phi).
\]

\begin{definition}
If the unique decomposition of
an element $S\in \CB^{ST}(X)^{Int}$ is\vspace{0.3cm}

\begin{equation}\label{decomp}
S = \mbox{exp}_\sqcup \left( \frac{1}{2}\sum_{ij} \wcoup \right)\sqcup R,
\end{equation}
\vspace{0.5cm}

then 
\vspace{0.1cm}

\begin{equation}\label{fgdefeq}
\int^{FG\, in\, ST}dX\, S = 
\left<
\mbox{exp}_\sqcup \left( -\frac{1}{2}\sum_{ij} \wcoupinv \right) , R
\right>_X\ \in\ \CB^{QST}(\phi).
\end{equation}
\vspace{0.5cm}

\end{definition}

\newpage

\section{Threading}\label{threadsect}

In this section we will introduce $Thr^D$, the operation 
of threading (rational) winding diagrams.
This is used in the
following square from the master
diagram. In this section we will show it commutes.

\

\begin{diagram}\label{threaddiag}

\

\[
\Picture{
\put(-6,0){$
\Picture{
\put(-6,0){$
\left\{
\begin{array}{c}
Special\ string\ links \\
in\ the\ solid\ torus.
\end{array}
\right\}
$}
\put(6.3,0){\vector(1,0){5}}
\put(4.75,-8.5){\vector(1,0){7}}
\put(7.75,0.25){Thr}
\put(3.5,-10.75){$S \mapsto Thr^D(S)\sqcup \nu(k)$}
\put(12,0){$
\left\{
\begin{array}{l}
\mbox{Special} \\
\mbox{tangles.}
\end{array}
\right\} $}
\put(-4.5,-4.5){$\sigma \circ \KIC^{ST}$}
\put(11.5,-4.5){$\sigma \circ \KIC$}
\put(0,-2){\vector(0,-1){5}}
\put(15,-2){\vector(0,-1){5}}
\put(-2,-9){$\CB^{ST}(X)^{Int}$}
\put(13,-9){$\CB(X,\underline{k})$}
}$}}
\]
\end{diagram}
\vspace{4cm}

\begin{itemize}
\item{ ``Special string link in the solid torus'' is defined in Definition \ref{ssldefn}, and ``Special
tangle'' is defined in Definition \ref{sftdefn}. }
\item{ $Thr$ is the operation introduced in Definition \ref{thrdefn} which 
threads the hole
in the solid torus with a zero-framed unknot. }
\item{ The space $\CB^{ST}(X)^{Int}$ is defined in Definition 
\ref{integdefn}, and
$\CB(X,\underline{k})$ is recalled in Section \ref{notsect}.}
\item{ The operation $Thr^D$ will presently be introduced.}
\end{itemize}

This operation will be defined on any of the spaces introduced so 
far which involve (possibly rational) winding coupons. If the 
labelling information is denoted $L$ (so maybe a skeleton and a set
of labels, or possibly $\phi$) then the map will be between spaces as
follows. The $Q$ is in brackets because it may, or may not, be present.
\[
Thr^D : \CB^{(Q)ST}(L) \rightarrow \CB(L,\underline{k}).
\]

The definition of this map will be introduced via an intermediate
construction, a generating diagram. 
\begin{definition}
A {\bf generating diagram} is precisely
the same as a winding diagram, except that
its coupons are labelled 
with formal power series in a variable $k$. 
A generating diagram denotes a particular series of uni-trivalent diagrams.
\begin{enumerate}
\item{(One edge)
If there is only one edge going through the coupon, then the association is
as follows.

If $f(k)\in \Qset[[k]]$
is written $f_0 + f_1 k + f_2  k^2 + f_3 k^3 + \ldots$ then a coupon
labelled as follows, is to be expanded as shown. See below for the meaning
of the vertex.

\

\[
\fcoupon\hspace{0.05cm}= \hspace{0.05cm}f_0\hspace{0.25cm}
\fcoupa+\hspace{0.15cm}f_1\hspace{0.25cm} \fcoupbx
\hspace{1cm}+\hspace{0.15cm}f_2\hspace{0.25cm} \fcoupcx
\hspace{1cm}+\hspace{0.15cm}f_3\hspace{0.25cm} \fcoupdx \hspace{1cm} \ldots
\]\vspace{0.8cm}

}
\item{(More than one edge) In this case, one takes the above expansion
and then, for each diagram, takes the sum of diagrams obtained by lifting
each introduced leg to each edge going through the coupon.}
\end{enumerate}

This vertex depends on whether the edge is internal or part of the skeleton,
and then also on the orientation of the skeleton, as follows. 

\

\[
\fcoupbx \hspace{0.8cm} \rightarrow \hspace{0.25cm} \fcoupb
\hspace{0.8cm} ; \hspace{0.25cm} 
\fskx \hspace{0.8cm} \rightarrow \hspace{0.25cm} \fsk \hspace{0.8cm};
\hspace{0.25cm} 
\fskxr \hspace{0.8cm} \rightarrow\ (-1)\hspace{0.25cm} \fskr
\]\vspace{1cm}

\end{definition}

\begin{definition}\label{threaddefn}
The operation $Thr^D$ is defined by mapping a (rational) winding diagram
to the series of uni-trivalent diagrams represented by making the substitution
$t\mapsto e^k$ in the label of every winding coupon.
\end{definition}

\begin{lemma}
This is a well-defined operation.
\end{lemma}
For starters, the subsitution makes sense because we are restricting
labels
to the subspace of rational functions that are non-singular at $t=1$.
Furthermore, 
the relations involving winding coupons (OR, ADD, MULT, SPLIT, COMM, PUSH1
and PUSH2) are easily checked.

\begin{theorem}\label{threadmaintheorem}
Diagram \ref{threaddiag} commutes.
\end{theorem}

This theorem depends crucially on a certain consequence of the Kontsevich
integral proof of the ``Wheels Conjectures'' that has recently
been given by Bar-Natan, Le and Thurston \cite{TW}.
See the forthcoming paper
of Bar-Natan and Lawrence \cite{BNL} for 
the following calculation. (See that paper also for references to
other proofs of the ``Wheeling Conjecture'' that are in the literature.)

\begin{theorem}

\

\vspace{0.75cm}

\[
\KIH(\ \ \ \ \ \ \ \hopfundone\ \ )\  = \ \ \mbox{exp}_{\sqcup}(\  \ijchordWh
)\sqcup \nu(k)\ \in\ \CB(x,\underline{k}).
\]
\end{theorem}

\vspace{1.2cm}

We use the following corollary. The tangle below is
equipped with some choice of bracketting which is the same on 
both the top and the bottom boundary 
words.

\begin{corollary}

\

\vspace{1cm}

\[
\ \ \KIH(\hspace{1.9cm} \ \corolpic\hspace{1.4cm} ) =\hspace{2cm} \corolpict\hspace{2.5cm}\sqcup \nu(k) 
\]
\vspace{1.75cm}

in the space
$ 
\CA( \underbrace{\uparrow\ldots\uparrow}_r \underbrace{\downarrow\ldots\downarrow}_s ;  \underline{k}).$
\end{corollary}

\underline{Proof of corollary.}
This is proved with Le and Murakami's paralleling formula \cite{LM2}.
The small point to observe is that an application of the paralleling 
formula gives 

\[
\mbox{exp}_{\sqcup}( \ijchordWha )\sqcup \ldots \sqcup 
\mbox{exp}_{\sqcup}( \ijchordWhr ) 
\sqcup
\mbox{exp}_{\sqcup}( -\ijchordWhrp )
\sqcup \ldots \sqcup 
\mbox{exp}_{\sqcup}( -\ijchordWhrs ) \sqcup \nu(k). 
\]

\

This requires the {\it averaged sum of orderings} of legs on the 
components $\{x_i\}$, whereas the given statement requires the composition
exponential. But the other ends of these chords are unordered. 
So the corollary follows as stated, for $\KIH$. (Observe that there
are no problems with the choice of associator for this scenario).
\rtb

\underline{Proof of Theorem \ref{threadmaintheorem}}.
Take $T$, an $\mu$-string 
string link in the solid torus, given by some presentation
$(A_T,B_T,w_1,w_2)$. Without loss of generality we can assume that
$w_2$ is some bracketting of some word of 
the form $\underbrace{\uparrow\ldots\uparrow}_r \underbrace{\downarrow\ldots\downarrow}_s$.

Then, the functoriality of $\KIH$ indicates that $\KIC( Thr(T) )$ is equal
to the following expression, in the space $\CA({\uparrow^\mu},\underline{k}).$\vspace{0.2cm}

\[
\Picture{
\put(0,1){$\Picture{
\put(-16,-2){$ \left( (\underbrace{\nu \otimes \ldots \otimes \nu}_{\mu})\circ
\Delta^{\mu-1}(\nu)\circ \right.
$}
\put(-12,-8){$\left. \KIH(A_T) \circ \left(I_{w_1}\otimes\hspace{2cm} \corolpict
\hspace{2.5cm} \right) \circ \KIH(B_T)\right)\sqcup \nu(k)$}
}$}}
\]

\vspace{4cm}

Alternatively, $Thr^D(\KIC^{ST}(T))\sqcup \nu(k)$ may be calculated, as 
follows.
\[
Thr^D\left(
(\underbrace{\nu \otimes \ldots \otimes \nu}_{\mu}) \circ
\Delta^{\mu}(\nu)\circ \gamma(\KIH( A_T)) \circ (I_{w_1}\otimes G_{w_2}) \circ \gamma(\KIH( B_T ))\right) \sqcup \nu(k).
\]
This is exactly the same thing, in the space $\CA({\uparrow^\mu},\underline{k})$.
We have just proved the commutation of the top half of the 
following diagram.
The precise statement of the theorem
then follows from the commutation of the bottom half. 
\newpage

\vspace{1cm}

\[
\Picture{
\put(-12,0){$
\left\{
\begin{array}{c}
Special\ string\ links \\
in\ the\ solid\ torus.
\end{array}
\right\}
$}
\put(-4.5,-2){\vector(0,-1){2}}
\put(8.5,-2){\vector(0,-1){2}}
\put(1,0){$\stackrel{Thr}{\longrightarrow}$}
\put(5.6,0){$
\left\{
\begin{array}{c}
Special \\
tangles.
\end{array}
\right\}
$}
\put(-3.5,-3){$\KIC^{ST} $}
\put(9.5,-3){$\KIC $}
\put(0,-9){$
\Picture{
\put(-7,2){$\CA^{ST,spec}({\uparrow^\mu})$}
\put(1,2){$\stackrel{Thr^D}{\longrightarrow}$}
\put(7,2){$\CA({\uparrow^\mu};\underline{k})$}
\put(-4.5,-0.5){\vector(0,-1){2}}
\put(8.5,-0.5){\vector(0,-1){2}}
\put(-6,-5){$\CB^{ST}(X)$}
\put(1,-5){$\stackrel{Thr^D}{\longrightarrow}$}
\put(6.5,-5){$\CB(X;\underline{k})$}
\put(-3.5,-1.5){$\sigma $}
\put(9.5,-1.5){$\sigma $}
}$}}
\]\vspace{5.5cm}

\rtb

\subsection{Evading $\CB^{QST}$}\label{evasion}

The reader who wishes to evade the space $\CB^{QST}$ may do so 
by redefining the map $Thr^{D}\circ \int^{FG\, in\, ST}$
as follows :

\begin{altdefinition}
If the unique decomposition of
an element $S\in \CB^{ST}(X)^{Int}$ is\vspace{0.3cm}

\begin{equation}
S = \mbox{exp}_\sqcup \left( \frac{1}{2}\sum_{ij} \wcoup \right)\sqcup R,
\end{equation}
\vspace{0.5cm}

then 
\vspace{0.1cm}

\begin{equation}
Thr^{D}\left(\int^{FG\, in\, ST}dX\, (S)\right) = 
\left<
\mbox{exp}_\sqcup \left( -\frac{1}{2}\sum_{ij} \wcoupinvk \right) , Thr^{D}(R)
\right>_X\ \in\ \CB(\underline{k}).
\end{equation}
\vspace{0.5cm}

\end{altdefinition}

\newpage

\section{The LMO invariant}

The LMO invariant was introduced by Thang Le, Jun Murakami and Tomotada
Ohtsuki \cite{LMO}. We will specialise our presentation to the setting in
question, knots in integral homology three-spheres.

Let the pair of a knot $K$ in an integral homology three sphere $M$ be 
presented by some special framed tangle $T$. 
That is, $T$ has one closed component,
and forgetting that closed component leaves the tangle
a string link $T'$ whose components are to be surgered, after closure; the
determinant of the linking matrix of those components is $\pm 1$.
Let $X=\{x_1,\ldots,x_\mu \}$ be an index set for $T'$, 
and let $lk(T')$ denote the linking matrix of those components. 

The LMO invariant is constructed as a sequence
\[
\KI_n^{LMO}(M,K) \in \CB_{\leq n}( \underline{k} )
\]
with the property that
\begin{equation}\label{coheqn}
\mbox{Grad}_{\leq n}( \KI^{LMO}_{m}(M,K) ) = \KI_n^{LMO}(M,K) \in \CB_{\leq n}(\underline{k})
\end{equation}
when $m\geq n$. This sequence can thus be regarded as approximations to 
an invariant defined by
\begin{eqnarray}\label{LMOassem}
& & \\
\KI^{LMO}(M,K) & = & 1 + \mbox{Grad}_1(\KI^{LMO}_1(M,K))
+ \mbox{Grad}_2(\KI^{LMO}_2(M,K))
+ \ldots \in \CB(\underline{k}). \nonumber
\end{eqnarray}

It is convenient to introduce the operation that is the heart of this
definition over a more general space.

\begin{definition}\label{loopdefn}
For some labelling set $L$, the space $\CB^{o}(L)$ is defined
with exactly the same definition as $\CB(L)$, 
except that generating diagrams may also have a finite
number of closed dashed loops.
\end{definition}

\begin{remark}
The space $\CB(L)$ clearly injects into $\CB^{o}(L)$. 
\end{remark}

\begin{definition}
The mapping
\[
\int^{(n)} dX : \CB^{o}( X ,X',\underline{k})
\rightarrow \CB_{\leq n}( X',\underline{k} ),
\]
is defined on a diagram $D$ as
\[
\left< \prod_{i=1}^{\mu}\left( \frac{1}{n!} \left( \frac{1}{2} \iiichord \right)^{n} \right) , D\right>,
\]
followed by the exchange of each dashed loop component
(some extra may arise) for a multiplicative
factor of $-2n$.
\end{definition}

With this operation, $\KI_n^{LMO}(M,K)$ is defined as follows. Let $\sigma_{\pm}( lk(T') )$ denote
the number of positive (resp. negative) eigenvalues of $lk(T')$. 
Let $U_{\pm}$ be a $\pm 1$-framed unknot.

\begin{definition}\label{LMOdefn}
The invariant $\KI_n^{LMO}(M,K)$ is defined by
\[
\KI_n^{LMO}(M,K) = \frac
{ \int^{(n)}dX \sigma(\KIC( T )) }
{ \left(\int^{(n)}dU \sigma(\KIC(U_+))\right)^{\sigma_+( lk(T') )}
\left(\int^{(n)}dU \sigma(\KIC(U_-))\right)^{\sigma_-( lk(T') )}
},
\]
in the space $\CB_{\leq n}(\underline{k})$.
\end{definition}

\begin{notation}
Let the numerator in the expression above be denoted $\KI^{LMO,o}_n(T)$.
Note that this $o$ is of a different nature to the $o$ used in Definition
\ref{loopdefn} above.
\end{notation}

\begin{theorem}\cite{LMO,LeD}
This definition gives a well-defined invariant of knots in integral homology
three-spheres, assmebled via Formula \ref{LMOassem}, with the specialisation
\[
\KI^{LMO}(S^3,K) = \KIH(K).
\]
\end{theorem}

\begin{remark}
Observe that our normalisation of the non-surgered components does not
affect surgery formula.
\end{remark}

\newpage

\section{The surgery formula}\label{surgopsect}

This section is the crux of the paper. It will prove that the following
square from the master diagram commutes.

\begin{diagram}\label{surgdiag}\label{surgform}

\[
\Picture{
\put(-5,5){$
\Picture{
\put(13,-7){$\CB(X,\underline{k})$}
\put(-4,-7){$\CB^{ST}(X)^{Int}$}
\put(-2,-8.5){\vector(0,-1){9}}
\put(14.5,-8.5){\vector(0,-1){9}}
\put(-7,-20){$\CB^{QST}(\phi)\times \Zset \times \Zset^1[t,t^{-1}] $}
\put(-3,-22){$(S,\sigma,P(t)) \mapsto Thr^{D}(S)\sqcup (-1)^{n\sigma}
{Wh'}(P(t)) $}
\put(13,-20){$\CB_{\leq n}(\underline{k})$}
\put(-1,-12.5){$\int^{FG\,in\,ST} dX \times \sigma_+ \times \mbox{det}$}
\put(15,-12.5){$\int^{(n)}dX$}
\put(3,-7){\vector(1,0){6}}
\put(2,-6.25){$S\mapsto Thr^D(S)\sqcup \nu(k)$}
\put(5,-20){\vector(1,0){5.5}}
}$}}
\]
\end{diagram}
\vspace{6.5cm}

There are some components of this that are yet to be introduced.
\begin{definition}
Let $\Zset^1[t,t^{-1}]$ be the subring of $\Zset[t,t^{-1}]$ of polynomials
such that if $p(t)\in \Zset^{1}[t,t^{-1}]$ then
\begin{enumerate}
\item{$p(t)=p(t^{-1}),$}
\item{$p(1)=\pm 1.$}
\end{enumerate}
\end{definition}
Clearly $det$, evaluated on $\CB^{ST}(X)^{Int}$ (according to Remark \ref{matrixeval}) 
takes values in this ring.
\begin{definition}
The mapping $\sigma_{+}$ on $S\in \CB^{ST}(X)^{Int}$ 
is the number of positive eigenvalues
of $W(1)$, where $W(t)$ is the Gaussian matrix of $S$.
\end{definition}
\begin{definition}
The mapping 
\[
Wh' : \Zset^1[t,t^{-1}] \rightarrow \CB(\underline{k}),
\]
is defined by
\[
Wh'(P(t)) = \mbox{exp}_{\sqcup}
\left(\left. \left[-\frac{1}{2} \mbox{log}\left( \frac{P(e^{h})}{P(1)}\right) 
\right] \right|_{h^{2n}\rightarrow \omega_{2n}}\right) \sqcup \nu(k),
\]
where the operation indicated is to expand the term inside the square
brackets into a power series in $h$, and then to replace terms like
$ch^{2n}$ by $c \omega_{2n}$, in exactly the same fashion as 
Definition \ref{wheelsdefn}.
\end{definition}
Observe that the $P(1)$ factor just adjusts the sign of the polynomial.
\begin{remark}
\[
Wh(M,K) = Wh'( A_{(M,K)}(t) ).
\]
\end{remark}

\subsection{Translation by power series}\label{transsect}

\

Before we consider the details of this proof, we introduce a certain 
operation on diagrams.
For $F(x_1,\ldots,x_\mu,k) \in \CB(X,\underline{k})$,
the notation 
$F(x_1+z,x_2,\ldots,x_\mu,k)$, according to Aarhus,
denotes the series 
in $\CB(X,\underline{k},z)$ 
obtained by replacing every diagram with $l$ legs labelled by $x_1$ 
by the $2^l$ diagrams
obtained by relabelling each such leg with either $x_1$ or $z$.
This may be extended 
linearly to simultaneous ``translations'' of other variables.

Now we introduce something novel. If $f(k) \in \Qset[[k]]$ then 
$F(f(k)x_1,x_2,\ldots,x_\mu)$ denotes the element obtained by
replacing every diagram in the expansion of $F$ with a generating diagram
obtained by labelling (simultaneously) 
every leg marked $x_1$ as follows. 
Observe that the added coupon is oriented according to 
the position of the univalent vertex.
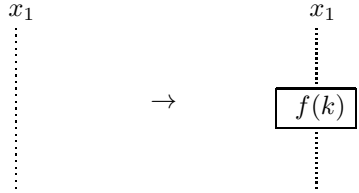
\vspace{1.3cm}

\[
\legf \hspace{1cm} \rightarrow \hspace{1.75cm} \legfb
\]

\vspace{1.6cm} 

We will want to combine these operations. For example, take some $f(k)\in \Qset[[k]]$:
\[
F(x_1+ f(k) z, x_2,\ldots, x_\mu, k) = 
F(x_1, x_2,\ldots, x_\mu, k) +
F(f(k) z, x_2,\ldots, x_\mu, k).
\]

\begin{notation}
Take a matrix $M(k) \in \mbox{M}_{\mu}( \Qset [[k]] )$. The notation
$F(\overline{x}+M(k)\overline{x}',k)$ represents the element
\[
F(x_1+\sum_{i_1} M_{1i_1}(k) x'_{i_1},
x_2+\sum_{i_2} M_{2i_2}(k)x'_{i_2},\,\ldots,x_\mu+\sum_{i_\mu} M_{\mu i_\mu}(k)x'_{i_\mu}, k) 
\]
in $\CB(X,X',\underline{k})$.
\end{notation}

\subsection{Translation invariance}

The core of our proof (our adaption of \cite{A3}) is the following
property.

\begin{theorem}\label{transinvprop}
Take some matrix of power series $M(k) \in {M}_{\mu}( \Qset [[k]] )$.
\[
\int^{(n)} dX F(\overline{x},{k}) =
\int^{(n)} dX F(\overline{x}+M(k)\overline{x}', k ) \ \ \in\ \
\CB_{\leq n}
(X',\underline{k}).
\]
\end{theorem}

The proof of this lemma requires a certain subspace of $\CB^{o}(X,X',\underline{k})$. This
is the subspace generated by $C_n$ vectors, which were introduced
in \cite{A3}. 

\begin{definition}[$C_n$ vectors]
Consider some uni-trivalent diagram drawn to include some box with $n$ 
dashed edges
attached (in an ordered fashion) to its top,
and $n$ dashed edges attached (ditto.) to its base. 
If the box is labelled by some permutation $\sigma \in \Sigma_n$, 
then this diagram 
is defined to represent the diagram obtained by joining up the edges
according to $\sigma$. A $C_n$ relation is a linear combination of diagrams
obtained from 
such a diagram by summing over all the diagrams
obtained by labelling that box with all possible permutations.
\end{definition}

\begin{lemma}[\cite{LMO,A3}]

A $C_{m}$ vector is in the kernel of $\int^{(n)}dX$, if $m\geq 2n+1$.

\end{lemma}

\begin{remark}
This is a slightly different viewpoint than is taken in the literature.
LMO works with a different relation $P_{n+1}$. This was shown to be
equivalent to $C_{2n+1}$ in \cite{A3}.
LMO introduces the relation $P_{n+1}$ and then shows that it is ``redundant''
in the image (that is, the quotient by diagrams of grade greater than $n$). 
Here we are in 
a more general situation than is explicitly found in the literature:
namely, we are allowing non-surgered univalent vertices. The proof
that $P_{n+1}$ is ``redundant'', which is accessibly described in
Section 2.5.4 of \cite{LeGr}, adapts immediately to this generality.
Following the discussion there, one sees that every leg on a $P_{n+1}$
must still meet a seperate vertex.

\end{remark}

\begin{flushleft}
\underline{Proof of Theorem \ref{transinvprop}.}
\end{flushleft}

Consider some diagram $D(\overline{x},{k})$ appearing
in the expansion of $F(\overline{x},{k})$. There are three
cases which cover the possibilities.
{\bf Case 1}: there is some component $x_i$ which has less than
$2n$ legs labelled with it. In this case both $D(\overline{x},k)$ and 
$D(\overline{x}+M\overline{x}',k)$ are in the kernel of $\int^{(n)}dX$.
{\bf Case 2}: each component of $X$ has exactly $2n$ legs labelled with it.
In this case there is exactly one contributing diagram on the right
hand side, because if any leg $x_i$ is relabelled $M_{ij}x_j'$ then 
there will then be less than $2n$ legs labelled with that component (such
a diagram will
be in the kernel.) The one contributing term is the same as that obtained
from the left hand side.
{\bf Case 3}: some component has more than $2n$ legs labelled with it.
This vanishes on the left hand side by definition.
This also vanishes on the right hand side because this is expressible
as a series of
$C_m$ vectors, for some $m\geq 2n+1$.
\rtb

\subsection{Diagram \ref{surgdiag} commutes}

Take some element $S(\overline{x})$ of $\CB^{ST}(X)^{Int}$, with canonical decomposition:
\vspace{0.1cm}

\[
S(\overline{x}) = \mbox{exp}_\sqcup \left( \frac{1}{2}\sum_{ij} \wcoup \right)\sqcup R(\overline{x}).
\]\vspace{0.4cm}

We shall calculate $\int^{(n)}dX( Thr^D( S(\overline{x})) \sqcup \nu(k) )$.
All expressions in the following are to 
be evaluated in the space $\CB_{\leq n}(X)$.

The first step is to split off the remainder, $Thr^D( R(\overline{x}) )$.
Denote this $R'(\overline{x},k)$.
We can do the following
because $R'(\overline{x},k)$ is a series of X-substantial diagrams.

\vspace{3cm}

\[
\Picture{
\put(-15,6){$
\int^{(n)}dX\ \left( \mbox{exp}_{\sqcup}\left( \frac{1}{2}\sum_{i,j} \wcoupk \right) \sqcup 
R'(\overline{x},k) \sqcup \nu(k) \right) $}
\put(-14,-1){$
=
\left< R'(\overline{x}',k),\int^{(n)} dX \mbox{exp}_{\sqcup}\left( \frac{1}{2}\sum_{i,j}
\wcoupk + \sum_{i} \wcoupii \right)\right>_{X'}\sqcup\, \nu(k).$}}
\] 

\vspace{5cm}

Theorem \ref{transinvprop} tells us that making the ``translation''
\[
x_i \rightarrow x_i - \sum_{j} W^{-1}_{ij}(e^{k}) x_{j}'
\]
inside the integrand will not affect the result.
[The matrix $W(T,e^k)$ is invertible in $M_\mu( \Qset[[k]] )$
by assumption (Definition \ref{integdefn}.) ]

Making that substitution transforms the integral into the following.
[Note that this requires the property that $W^{-1}(T,e^{k}$) is 
Hermitian, which is again by assumption (Definition \ref{integdefn}.)]
\vspace{0.25cm}

\[
\int^{(n)}dX \mbox{exp}_{\sqcup}\left( \frac{1}{2}\sum_{i,j} \wcoupk
- \frac{1}{2} \sum_{i,j} \wcoupinvkp \,\right) .
\]\vspace{0.75cm}

Substituting this expression back into the above gives the following.
\vspace{0.75cm}


\begin{equation}\label{finaleqn}
\Picture{
\put(-10,0){$
\left(
\nu(k)\sqcup 
\int^{(n)}dX \mbox{exp}_{\sqcup}\left( \frac{1}{2}\sum_{i,j} \wcoupk 
\right)\right) $}
\put(-2,-7){$
\sqcup \left< \mbox{exp}_{\sqcup}\left( -\frac{1}{2}\sum_{i,j}\ \wcoupinvk \right)\ ,\ 
R'(\overline{x},k) \right> \ \ \ \ \  \, $}}
\end{equation}\vspace{3.25cm}

\hspace{-0.5cm} evaluated 
in $\CB_{\leq n}(\underline{k})$. The theorem follows from the
calculation of the leading term that is performed in the following
section.
\rtb

\newpage

\subsection{The wheels line}\label{wheelsect}

Examining the structure of Equation \ref{finaleqn}, we see
that the projection onto the wheels line of the
element $\int^{(n)}dX(Thr^D(S(\overline{x}))\sqcup \nu(k))$
is precisely the following function of the entries of the 
matrix $W(t)$.\vspace{0.2cm}

It is calculated as follows.
\begin{theorem}
Let $W(t)$ be a Hermitian matrix satisfying det$(W(1))=\pm 1$. Then\vspace{0.25cm}

\[\left(
\nu(k)\sqcup 
\int^{(n)}dX \mbox{exp}_{\sqcup}\left( \frac{1}{2}\sum_{i,j} \wcoupk 
\right)\right) = (-1)^{n\sigma_+(W(1))} Wh'(\mbox{det}(W(t))).
\]\
\end{theorem}
\vspace{0.5cm}

\underline{Proof}.

We begin by observing that every such matrix $W(t)$ can be realised
as the winding matrix of some special string link in the solid torus $T$,
presenting some pair $(M_T,K_T)$.
The calculation is performed by calculating the wheels line
of $\KI^{LMO}(M_T,K_T)$ in two different ways.

On the one hand, 
this ``wheels line''  has already been calculated by other means.
Many authors contributed to this result.
Let us highlight the original conjecture of Melvin--Morton \cite{MM} and
the Kontsevich integral proof given by   
Bar-Natan--Garoufalidis
\cite{BNG} (acknowledging other contemporaneous proofs \cite{RMM,KM}).
See Garoufalidis--Habegger \cite{GH} for
the following formula in the setting of $\Zset HS^3$s (recently
extended to null-homologous knots in $\Qset HS^3s$, \cite{L}).

\begin{theorem}\label{wheelthm}
Let $T$ be a $\mu$-string string link in the solid torus presenting
some pair $(M_T,K_T)$.
Then,
\[
\KI^{LMO}(M_T,K_T) = Wh'(\mbox{det}(W(T,t))) \sqcup (1+R),
\]
where $R$ is a series of diagrams whose dashed graphs have 
negative Euler charcteristic.
\end{theorem}

On the other hand, we have just seen (Equation \ref{finaleqn})
that $\KI^{LMO}(M_T,K_T)$ may be calculated
 
\[
\frac{
[Wheels\ bit] \sqcup [The\ rest]}
{ \left(\int^{(n)}dU \sigma(\KIC(U_+))\right)^
{\sigma_+(W(1))}
\left(\int^{(n)}dU \sigma(\KIC(U_-))\right)^{\sigma_-(W(1))}}
\]
in the space $\CB_{\leq n}(\underline{k})$, where \vspace{0.5cm}

\begin{eqnarray*}
\ [Wheels\ bit] & = & \left(
\nu(k)\sqcup 
\int^{(n)}dX \mbox{exp}_{\sqcup}\left( \frac{1}{2}\sum_{i,j} \wcoupk 
\right)\right), \\
& & \\
& & \\
& & \\
& & \\
& & \\
\ [The\ rest] & = & \left< \mbox{exp}_{\sqcup}\left( -\frac{1}{2}\sum_{i,j}\ \wcoupinvk \right)\ ,\ 
R'(\overline{x},k) \right>. \\
& & \\
& & 
\end{eqnarray*}
Given that \cite{LMO}
\[
\int^{(n)}dU \sigma(\KIC(U_{\pm 1})) = 
(\mp 1)^{n} + \mbox{terms of negative Euler characteristic,}
\]
the theorem follows by equating wheels lines. 
\rtb

\begin{exercize}
Prove Theorem \ref{wheelthm} directly. That is, calculate\vspace{0.5cm}

\[
\int^{(n)}dX 
\mbox{exp}_{\sqcup}\left( \frac{1}{2}\sum_{i,j} \wcoupk \right).
\]
\end{exercize}\vspace{0.75cm}

This is a recommended, non-trivial, exercise; 
it provides some perspective on what lies behind the commutativity
of Diagram \ref{surgform}
and gives a particular 
sense in which wheels tangibly become cycles of something.
Hint: Step A is the Aarhus formula; Step B is 
an appropriate identity from determinant theory, being
careful with combinatorial factors.

\section{Rozansky's rationality conjecture}\label{mainthsect}

We start by recalling the surgery formula. Just say the
pair $(M,K)$ is presented by
some $T$, a $\mu$-string string link in the solid torus. If the
associated decomposition is denoted\vspace{0.3cm}

\[
\sigma(\KIC^{ST}(T)) = \mbox{exp}_{\sqcup}\left(\ \frac{1}{2}\sum_{i,j} \wcoup
\ \right)\sqcup R(\overline{x}) 
\in \CB^{ST}(X),
\]\vspace{0.5cm}

where $R(\overline{x})$ is a series of $X$-substantial terms,
then $\KI^{LMO}_n(M,K)$ is equal to \vspace{0.5cm}

\[
\frac{
\begin{array}{l}
Wh(M,K) \sqcup \left< \mbox{exp}_{\sqcup}\left( -\frac{1}{2}\sum_{i,j}\ \wcoupinvk \ \ \right)\ ,\ 
Thr^D(R(\overline{x})) \right> \\ \\ \\ \, \end{array}}
{ \left((-1)^n\int^{(n)}dU \sigma(\KIC(U_+))\right)^
{\sigma_+(W(1))}
\left(\int^{(n)}dU \sigma(\KIC(U_-))\right)^{\sigma_-(W(1))}}
\ \ \in\ \ \CB_{\leq n}(\underline{k}).
\]

We already know that this is a group-like element. In fact, we know that
\[
\KI^{LMO}(M,K) = Wh(M,K) \sqcup \mbox{exp}(q),
\]
where $q$ is a series of connected diagrams of degree greater than one 
whose dashed graphs have
Euler characteristic less than zero. Thus we can read $q$ directly from
the expression above.
Write $q=\sum_{i=1}^{\infty} q^{(i)}$ where $q^{(i)}$ 
is a series of diagrams
of Euler characteristic minus $i$.
Comparing the two expressions, it is clear that  
$q^{(i)}$ is calculated from precisely those
diagrams appearing 
in $R(\overline{x})$ which have $2i$ trivalent vertices, 
plus some finite 
contribution from the signature corrections. 

In other words, the term $q^{(i)}$ 
is a sum over those ways of gluing labelled chords into pairs of legs of
some finite combination of polynomial generating diagrams
which result in connected diagrams. 

The correspondence is completed by noting that the determinant of the matrix
$W(T,e^{k})$ is $\pm A_{(M,K)}(e^{k})$ 
so that the entries of the matrix
$W^{-1}(T,e^{k})$ lie in $L_{(M,K)}$. 
Thus the edge-labels on the chords 
gluing the legs fall into that subspace, and the labels on the edges 
of the diagrams that are assembled according to the gluing
formula also lie in this subspace.
\rtb

It is an informative exercise to prove this theorem
without appealing to the prior-known fact that the result is group-like.

\bibliographystyle{amsalpha}

\end{document}